\documentclass[12pt]{amsart}
\usepackage{amssymb}
\usepackage[mathcal]{eucal}
\usepackage{lmodern}
\usepackage[T1]{fontenc}
\usepackage{hyperref}
\usepackage{enumitem}
\usepackage{nicematrix}
\usepackage{tikz-cd}
\usepackage[nobysame]{amsrefs}
\usepackage[total={6.25truein,9truein},centering]{geometry}
\usepackage{mathdots}

\allowdisplaybreaks

\NiceMatrixOptions{
  xdots/line-style = {dotted, thick},
  xdots/shorten-start = 5pt,
  xdots/shorten-end = 5pt
}

\numberwithin{equation}{section}

\newtheorem{theorem}[equation]{Theorem}
\newtheorem{lemma}[equation]{Lemma}
\newtheorem{proposition}[equation]{Proposition}
\newtheorem{corollary}[equation]{Corollary}

\newtheorem{theoremintro}{Theorem}

\theoremstyle{definition}

\theoremstyle{remark}

\newtheorem{remark}[equation]{Remark}

\newtheorem*{acknowledgments}{Acknowledgments}

\newenvironment{subsec}[1]{\vspace{\topsep}
\refstepcounter{equation}\theequation.\ \textit{#1.}}{\vspace{\topsep}}

\newenvironment{alphenumerate}{\begin{enumerate}[label=\textup{(\alph*)}]}{\end{enumerate}}

\newcommand{\FF}{\mathbb{F}}
\newcommand{\ZZ}{\mathbb{Z}}
\newcommand{\QQ}{\mathbb{Q}}

\newcommand{\CC}{\mathbb{C}}
\newcommand{\PP}{\mathbb{P}}
\newcommand{\NN}{\mathbb{N}}
\newcommand{\TT}{\mathbb{T}}

\newcommand{\HH}{\mathbb{H}}

\newcommand{\SSS}{\mathbb{S}}

\newcommand{\bB}{\mathbf{B}}
\newcommand{\bC}{\mathbf{C}}
\newcommand{\bD}{\mathbf{D}}
\newcommand{\be}{\mathbf{e}}

\newcommand{\bP}{\mathbf{P}}
\newcommand{\bQ}{\mathbf{Q}}

\newcommand{\bV}{\mathbf{V}}
\newcommand{\bW}{\mathbf{W}}

\newcommand{\bz}{\mathbf{z}}

\newcommand{\cB}{\mathcal{B}}
\newcommand{\cC}{\mathcal{C}}

\newcommand{\cE}{\mathcal{E}}
\newcommand{\cF}{\mathcal{F}}
\newcommand{\cG}{\mathcal{G}}

\newcommand{\cO}{\mathcal{O}}

\newcommand{\cS}{\mathcal{S}}
\newcommand{\cT}{\mathcal{T}}

\newcommand{\fG}{\mathfrak{G}}
\newcommand{\fp}{\mathfrak{p}}

\newcommand{\rN}{\mathrm{N}}

\newcommand{\rR}{\mathrm{R}}

\newcommand{\rV}{\mathrm{V}}

\newcommand{\sC}{\mathsf{C}}

\newcommand{\bschi}{\boldsymbol{\chi}}

\newcommand{\bszeta}{\boldsymbol{\zeta}}
\newcommand{\bsxi}{\boldsymbol{\xi}}
\newcommand{\bsomega}{\boldsymbol{\omega}}

\DeclareMathOperator{\diag}{diag}

\DeclareMathOperator{\id}{id}

\DeclareMathOperator{\GL}{GL}

\DeclareMathOperator{\Mat}{Mat}

\DeclareMathOperator{\End}{End}
\DeclareMathOperator{\res}{res}
\DeclareMathOperator{\Sf}{Sf}
\DeclareMathOperator{\Span}{Span}

\DeclareMathOperator{\Gal}{Gal}
\DeclareMathOperator{\Hom}{Hom}

\DeclareMathOperator{\Res}{Res}
\DeclareMathOperator{\sgn}{sgn}
\DeclareMathOperator{\ord}{ord}

\DeclareMathOperator{\im}{im}

\newcommand{\oK}{\mkern2.5mu\overline{\mkern-2.5mu K}}

\newcommand{\hM}{\hat{M}}
\newcommand{\hN}{\hat{N}}

\newcommand{\cont}{\mathrm{cont}}

\newcommand{\tors}{\mathrm{tors}}

\newcommand{\tbe}{\widetilde{\be}}

\newcommand{\teta}{\tilde{\eta}}

\newcommand{\tpi}{\tilde{\pi}}

\newcommand{\C}{\CC_{\infty}}

\newcommand{\oFF}{\overline{\FF}}

\newcommand{\iso}{\stackrel{\sim}{\longrightarrow}}

\newcommand{\power}[2]{{#1 [\![ #2 ]\!]}}
\newcommand{\laurent}[2]{{#1 (\!( #2 )\!)}}
\newcommand{\brac}[2]{\genfrac{\{}{\}}{0pt}{}{#1}{#2}}

\newcommand{\dnorm}[1]{\lVert #1 \rVert}

\newcommand{\inorm}[1]{{\lvert #1 \rvert}_{\infty}}

\newcommand{\hotimes}{\mathbin{\hat{\otimes}}}
\newcommand{\tr}{{\mathsf{T}}}
\newcommand{\sumprime}[1]{\sideset{}{'}\sum_{#1}}
\newcommand{\assign}{\mathrel{\vcenter{\baselineskip0.5ex \lineskiplimit0pt
                     \hbox{\scriptsize.}\hbox{\scriptsize.}}}%
                     =}
\newcommand{\rassign}{=%
                     \mathrel{\vcenter{\baselineskip0.5ex \lineskiplimit0pt
                     \hbox{\scriptsize.}\hbox{\scriptsize.}}}%
                     }

\begin{document}

\title{Dual Gauss-Thakur sums for Drinfeld modules}

\author{Athena Bassil}
\address{Department of Mathematics, Texas A{\&}M University, College Station, TX 77843, U.S.A.}
\email[A.\ Bassil]{athenabassil@tamu.edu}

\author{Matthew Papanikolas}
\email[M.\ Papanikolas]{papanikolas@tamu.edu}


\date{November 18, 2026}

\begin{abstract}
We define dual Gauss-Thakur sums for Drinfeld modules, building on constructions of Gauss-Thakur sums due to Thakur, in the case of the Carlitz module, and Gazda and Maurischat, for general Anderson $t$-modules. Using results of Ferraro, we demonstrate that these sums can be obtained by evaluations of dual special functions at elements of finite fields. Our main results are matrix identities assembled from Gauss-Thakur sums and dual Gauss-Thakur sums in terms of multiplication coefficients for the Drinfeld module and its adjoint. We further prove determinant identities and explore connections with special Goss $L$-values and Drinfeld modular forms.
\end{abstract}


\maketitle

\tableofcontents

\section{Introduction} \label{S:Intro}

Thakur defined function field analogues of Gauss sums in~\cite{Thakur88} and established their fundamental properties, demonstrating a striking parallel with the classical theory. Then in the intervening years Gauss-Thakur sums have appeared as central quantities in various contexts (e.g., see \cites{AnglesPellarin14, AnglesPellarin15, AT17, Ferraro22, Ferraro25, Thakur91a, Thakur91b, Thakur93a, Thakur93b, Thakur}). Gazda and Maurischat~\cite{GazdaMaurischat21} expanded the definition of Gauss-Thakur sums, built from torsion on an arbitrary $t$-module, and made a crucial connection with modules of special functions, augmenting work of Angl\`es and Pellarin~\cite{AnglesPellarin15} for the Carlitz module. In the present paper we construct and investigate dual Gauss-Thakur sums associated to the a Drinfeld module and its adjoint, rooted in the Drinfeld module's rigid analytic trivializations and built in part on work of Ferraro~\cite{Ferraro25}.

\begin{subsec}{Drinfeld modules}
For a finite field $\FF_q$ with $q = p^m$ elements, let $A = \FF_q[\theta]$ denote the ring of polynomials over $\FF_q$ in a variable $\theta$. Let $k = \FF_q(\theta)$ be its fraction field, let $k_{\infty} = \laurent{\FF_q}{\theta^{-1}}$ be the completion of $k$ at its infinite place, and let $\C$ be the completion of an algebraic closure of~$k_{\infty}$. We let $\deg = -\ord_{\infty} : \C \to \QQ \cup \{ -\infty\}$, so for $x \in \C$, $\inorm{x} = q^{\deg x}$, and in particular $\inorm{\theta} = q$. We let~$A_+$ be the monic elements of~$A$.

Let $\phi : A \to K[\tau]$ be a Drinfeld module of rank~$r$ defined over an intermediate field $k \subseteq K \subseteq \C$, where $\tau$ is the $q$-th power Frobenius endomorphism,
\begin{equation} \label{E:rhodef}
\phi_\theta = \theta + b_1 \tau + \dots + b_r \tau^r, \quad b_r \neq 0.
\end{equation}
We construct also its adjoint $\phi^* : A \to \oK[\tau^{-1}]$ defined over the algebraic closure $\oK$ of $K$ (see \S\ref{SS:Drinfeldmodules}, \S\ref{SS:adjoints}). Both $\phi$ and $\phi^*$ can be extended in natural ways (see \S\ref{SS:formalphi}, \S\ref{SS:phistarext}) to $\FF_q$-algebra maps $\phi : k_{\infty} \to \laurent{\C}{\tau}$ and $\phi^* : k_{\infty} \to \laurent{\C}{\tau}$, which for $\alpha \in k_{\infty}$ we express
\[
\phi_\alpha = \sum_n \brac{\alpha}{n}_{\phi} \tau^n, \quad
\phi^*_{\alpha} = \sum_n \tau^n \brac{\alpha}{n}_{\phi^*}.
\]

To $\phi$ we associate two matrices $\Upsilon$, $\Psi \in \GL_r(\TT)$, over the Tate algebra $\TT \subseteq \power{\C}{t}$ of functions converging on the closed unit disk of $\C$, which are rigid analytic trivializations for the $t$-motive and $t$-comotive of $\phi$ respectively as in~\cites{ABP04, Chang20, CP12, ChenYT24, HartlJuschka20, KP23, NP24, P08, Pellarin08}. Specializations of these matrices at $t=\theta$ produce periods and quasi-periods of $\phi$, and here we consider also specializations at $t=\varepsilon \in \oFF_q$,  connected to Gauss-Thakur sums as in~\cites{AnglesPellarin15, GazdaMaurischat21} (see \eqref{E:gievalintro}, \eqref{E:gi*evalintro}, Theorem~\ref{T:DualGaussThakurSums/Zetas}). As part of these formulations we obtained new power series identities for $\Upsilon$ and $\Psi$ in terms of the multiplication coefficients of $\phi$ and $\phi^*$. As a higher rank version of \cite{ABP04}*{Lem.~5.4.1}, we prove (Proposition~\ref{P:PsiNUpsilon11}) for $N \in \ZZ$,
\begin{equation} \label{E:PsiNUpsilon11intro}
\bigl[ \Psi^{(N)} \Upsilon \bigr]_{11} = \begin{cases}
{\displaystyle \sum_{m=0}^{\infty} \brac{\theta^{-m-1}}{N}_{\phi} \, t^m,} & \textup{if $N \geq 0$,} \\[15pt]
{\displaystyle -\sum_{m=0}^{\infty} \brac{\theta^{-m-1}}{-N}_{\phi^*}\, t^m,} & \textup{if $N < 0$,}
\end{cases}
\end{equation}
where for a matrix $M$, we let $M_{ij} = [M]_{ij}$ denote its $(i,j)$-entry. This identity underlies our main duality result (Theorem~\ref{T:GaussSumduality1Intro}) outlined below.

The first column of~$\Upsilon$ and the first row of $\Psi$,
\[
\bsomega_{\phi} = (\omega_1, \dots, \omega_r)^{\tr}, \quad
\bszeta_{\phi} = (\zeta_1, \dots, \zeta_r),
\]
consist of special functions (or Anderson generating functions) and dual special functions as defined in~\cites{AndThak90, EP14, GazdaMaurischat21, Goss94, Ferraro22, Ferraro25, Maurischat22, Pellarin08}, which satisfy certain Frobenius difference equations. As presented here these functions depend on a choice of $A$-basis $\pi_1, \dots, \pi_r$ for the period lattice $\Lambda_{\phi}$. See~\S\ref{SS:RATs} for details.

While $\bsomega_{\phi}$ and $\bszeta_{\phi}$ had appeared in the literature in various guises for a while, Ferraro~\cites{Ferraro22, Ferraro25} showed that they represent universal objects in a fundamental pairing between special functions and dual special functions (see Theorems~\ref{T:FerraroUniv} and~\ref{T:Ferraropairing}). Moreover, after comparing previous work in~\cites{CP12, HartlJuschka20, NP24}, it follows from Ferraro's results (see Corollary~\ref{C:samezetas}) that each $\zeta_i$, previously arising through Anderson generating functions~\cites{CP12,HartlJuschka20}, has a lattice sum incarnation,
\begin{equation}
\zeta_i = -\sumprime{(a_1, \dots, a_r) \in A^r} \frac{a_i(t)}{a_1 \pi_1 + \cdots + a_r\pi_r},
\end{equation}
where $\Sigma'$ denotes that the sum omits the $0$-tuple.
\end{subsec}

\begin{subsec}{Gauss-Thakur sums and dual Gauss-Thakur sums} For each $1\leq i \leq r$, we define continuous $\FF_q$-linear functions $\be_i$, $\be_i^* : k_{\infty} \to \C$ so that
\[
\omega_i = \sum_{m=0}^{\infty} \be_i(\theta^{-m-1}) t^m,
\quad \zeta_i = \sum_{m=0}^{\infty} \be_i^*(\theta^{-m-1}) t^m.
\]
We have $\be_i(x) = \exp_{\phi}(\pi_i x)$, and $\be_i^*(x)$ is similar to functions defined in~\cite{ABP04}*{\S 5.2} (though note Remark~\ref{R:oldestar}). As a matter of course, for monic irreducible $\wp \in A$ and $\FF_{\wp} = A/\wp$, we have induced injections of $A$-modules,
\[
z \mapsto \be_i(z/\wp) : \FF_{\wp} \to \phi[\wp], \quad
z \mapsto \be_i^*(z/\wp) : \FF_{\wp} \to \phi^*[\wp],
\]
where $\phi[\wp]$ and $\phi^*[\wp]$ are the $\wp$-torsion submodules of $\phi$ and $\phi^*$ respectively. Furthermore,
\[
\phi[\wp] = \bigoplus_{i=1}^r \be_i(\wp^{-1} A), \quad
\phi^*[\wp] = \bigoplus_{i=1}^r \be_i^*(\wp^{-1} A),
\quad \textup{(see Lemma~\ref{L:e*propertiesfortorsion}).}
\]

Gazda and Maurischat~\cite{GazdaMaurischat21} defined Gauss-Thakur sums $g_i(\chi;a)$ and we similarly define their duals $g_i^*(\chi;a)$ in the following way. Letting $\chi \in \widehat{\FF}_q^{\times}$ be a multiplicative character and $a \in A$, for $1 \leq i \leq r$,
\begin{equation}
g_i(\chi;a) = - \sum_{z \in \FF_{\wp}^{\times}} \chi(z)^{-1} \be_i(z/\wp), \quad
g_i^*(\chi;a) = - \sum_{z \in \FF_{\wp}^{\times}} \chi(z)^{-1} \be_i^*(z/\wp).
\end{equation}
When $\phi$ is the Carlitz module, one recovers Thakur's Gauss sums from~\cite{Thakur88}. These Gauss-Thakur sums and their duals satisfy a number of properties that parallel classical Gauss sums (see~\S\ref{SS:Gauss-ThakurSums}), particularly that
\[
g_i(\chi;a) = \chi(a) g_i(\chi), \quad
g_i^*(\chi;a) = \chi(a) g_i^*(\chi),
\]
where $g_i(\chi) = g_i(\chi;1)$ and $g_i^*(\chi) = g_i^*(\chi;1)$. Gazda and Maurischat define also a tensor product version, which characterizes Gauss-Thakur sums as a class of special functions for~$\phi$, and then dual sums arise as a class of dual special functions (see Lemma~\ref{L:GTspaces}).
\end{subsec}

\begin{subsec}{Special function evaluations}
For $d=\deg \wp$ and a root $\varepsilon \in \FF_{q^d} \subseteq \C$ of $\wp$, the map $\chi_{\varepsilon} : A \to \FF_{q^d}$ defined by $\chi_{\varepsilon}(a) = a(\varepsilon)$ induces an evaluative Dirichlet character modulo~$\wp$. Thakur proved in the Carlitz case, and Gazda and Maurischat in general, that $g_i(\chi;a) = 0$ if $\chi$ is not evaluative. The same holds for $g_i^*(\chi;a)$, so we primarily consider only $g_i(\chi_k;a)$ and $g_i^*(\chi_k;a)$, where $\chi_k = \chi_{\varepsilon^{q^k}}$.

One of the main theorems of Gazda and Maurischat~\cite{GazdaMaurischat21}*{Thm.~5.2} identifies Gauss-Thakur sums with evaluations of $\bsomega_{\phi}$,
\begin{equation} \label{E:gievalintro}
g_i(\chi_{\varepsilon}; \wp'(\theta)) = \omega_i(\varepsilon), \quad 1 \leq i \leq r,
\end{equation}
where $\wp'(\theta)$ is the derivative of $\wp$. This result extends a theorem of Angl\`es and Pellarin~\cite{AnglesPellarin15}*{Thm.~2.9} for the Carlitz module and the Anderson-Thakur function~$\omega_{\sC}$ (see~\eqref{E:Carlitzomegaevals}). We prove a version of this result (Theorem~\ref{T:DualGaussThakurSums/Zetas}) for dual Gauss-Thakur sums in terms of~$\bszeta_{\phi}$,
\begin{equation} \label{E:gi*evalintro}
g_i^*(\chi_{\varepsilon}; \wp'(\theta)) = \zeta_i(\varepsilon), \quad 1 \leq i \leq r.
\end{equation}
For $k \in \ZZ$, we also twist Gauss-Thakur sums,
\[
g_i^{[k]}(\chi;a) = - \sum_{z \in \FF_{\wp}^{\times}} \chi(z)^{-1} \be_i(z/\wp)^{q^k}, \quad
g_i^{*[k]}(\chi;a) = - \sum_{z \in \FF_{\wp}^{\times}} \chi(z)^{-1} \be_i^*(z/\wp)^{q^k},
\]
and determine (Corollary~\ref{C:GaussSumsTwistsAndZetas}) that
\[
g_i^{[k]}(\chi_{\varepsilon}; \wp'(\theta)) = \omega_i^{(k)}(\varepsilon), \quad
g_i^{*[k]}(\chi_{\varepsilon}; \wp'(\theta)) = \zeta_i^{(k)}(\varepsilon), \quad 1 \leq i \leq r,
\]
i.e., that they are evaluations of Frobenius twists (see \S\ref{S:RAT}) of~$\bsomega_{\phi}$ and~$\bszeta_{\phi}$.
\end{subsec}

\begin{subsec}{Matrices of Gauss-Thakur sums and dual Gauss-Thakur sums} \label{SS:GTSMatrices}
Our main theorem below (Theorem~\ref{T:GaussSumduality1Intro}, and Theorem~\ref{T:GaussSumduality1} later) is a matrix identity for Gauss-Thakur sums, their duals, and their twists. For $\chi \in \widehat{\FF}_q^{\times}$ and $a \in A$, define
\begin{equation} \label{E:frakG}
\fG(\chi;a) \assign \begin{pNiceMatrix}
g_1(\chi;a) & g_1^{[1]}(\chi;a)& \Cdots & g_1^{[r-1]}(\chi;a)\\
g_2(\chi;a) & g_2^{[1]}(\chi;a) & \Cdots & g_2^{[r-1]}(\chi;a)\\
\Vdots & \Vdots &  & \Vdots & \\
g_r(\chi;a) & g_r^{[1]}(\chi;a) & \Cdots & g_r^{[r-1]}(\chi;a) \\
\end{pNiceMatrix}
\quad \in \Mat_{r\times r}(\oK)
\end{equation}
and
\begin{equation} \label{E:frakG*}
\fG^*(\chi;a) \assign \begin{pNiceMatrix}
g_1^*(\chi;a) & g_2^*(\chi;a) & \Cdots & g_r^*(\chi;a) \\
g_1^{*[1]}(\chi;a) & g_2^{*[1]}(\chi;a) & \Cdots & g_r^{*[1]}(\chi;a) \\
\Vdots & \Vdots &  & \Vdots & \\
g_1^{*[r-1]}(\chi;a) & g_2^{*[r-1]}(\chi;a) & \Cdots & g_r^{*[r-1]}(\chi;a) \\
\end{pNiceMatrix}
\quad \in \Mat_{r\times r}(\oK).
\end{equation}
Keeping with previous conventions, we set $\fG(\chi) \assign \fG(\chi;1)$ and $\fG^*(\chi) \assign \fG^*(\chi;1)$. We find that just as above, these matrices are obtained by specializing our rigid analytic trivializations $\Upsilon$ and $\Psi$ (see~\eqref{E:matrixGaussZetas}),
\[
\fG(\chi_{\varepsilon};\wp'(\theta)) = \wp'(\varepsilon)\fG(\chi_{\varepsilon}) = \Upsilon(\varepsilon), \quad
\fG^*(\chi_{\varepsilon};\wp'(\theta)) = \wp'(\varepsilon) \fG^*(\chi_{\varepsilon}) = J\Psi^{(r-1)}(\varepsilon),
\]
where
\[
J = \left( \begin{smallmatrix} & & 1 \\ & \iddots & \\ 1 & &\end{smallmatrix} \right) \in \Mat_{r \times r}(\FF_q).
\]
Now for $a \in A$, define the $rd \times \NN$ matrix,
\begin{equation} \label{E:Gammadef}
\Gamma(\phi;a) \assign \begin{pNiceMatrix}
\fG(\chi_0;a) & \fG^{[r]}(\chi_0; a) & \fG^{[2r]}(\chi_0; a) & \Cdots[shorten-end=10pt] \\
\fG(\chi_1;a) & \fG^{[r]}(\chi_1; a) & \fG^{[2r]}(\chi_1;a) & \Cdots[shorten-end=10pt]  \\
\Vdots & \Vdots & \Vdots  &  \\
\fG(\chi_{d-1};a) & \fG^{[r]}(\chi_{d-1}; a) & \fG^{[2r]}(\chi_{d-1}; a) & \Cdots[shorten-end=10pt] 
\end{pNiceMatrix}.
\end{equation}
Define also the $\NN \times rd$ matrix,
\begin{equation} \label{E:Gamma*def}
\Gamma^*(\phi^*;a) \assign \begin{pNiceMatrix}
\fG^*(\chi_0;a) & \fG(\chi_1;a) & \Cdots & \fG^*(\chi_{d-1};a) \\
\fG^{*[r]}(\chi_0;a) & \fG^{[r]}(\chi_1;a) & \Cdots & \fG^{*[r]}(\chi_{d-1};a) \\
\fG^{*[2r]}(\chi_0;a) & \fG^{*[2r]}(\chi_1;a) & \Cdots & \fG^{*[2r]}(\chi_{d-1};a) \\
\Vdots[shorten-end=5pt] & \Vdots[shorten-end=5pt] &  & \Vdots[shorten-end=5pt]
\end{pNiceMatrix}.
\end{equation}
Here $\fG^{[k]}(\chi;a)$ and $\fG^{*[k]}(\chi;a)$ are obtained by twisting each entry of $\fG(\chi;a)$ and $\fG^{*}(\chi;a)$ $k$~times respectively as above. For $\ell \in \NN$, we let $\Gamma(\phi;a) \in \Mat_{rd\times r\ell}(\oK)$ be the truncation of $\Gamma(\phi;a)$ after $r\ell$ column blocks, and let $\Gamma_{\ell}^*(\phi^*;a) \in \Mat_{r\ell \times rd}(\oK)$ be the truncation of $\Gamma^*(\phi^*;a)$ after $r\ell$ row blocks. We also set $\Gamma(\phi) \assign \Gamma(\phi;1)$, $\Gamma^*(\phi^*) \assign \Gamma^*(\phi^*;1)$, etc.

Now for $\beta \in k$ that is regular at $\varepsilon$, define the block diagonal matrix,
\begin{equation} \label{E:bQintro}
\bQ(\beta)\assign \diag \bigl( \beta(\varepsilon)I_r, \dots, \beta(\varepsilon^{q^{d-1}})I_r \bigr)
\quad \in \Mat_{rd\times rd}(\FF_{q^d}),
\end{equation}
where $I_r$ is the $r \times r$ identity matrix.  For $\beta \in k$ define
\begin{equation} \label{E:bPintro}
\bP(\beta) \assign \begin{pNiceMatrix}
\beta & -\brac{\beta}{1}_{\phi^*}^q & -\brac{\beta}{2}_{\phi^*}^{q^2} & -\brac{\beta}{3}_{\phi^*}^{q^3} & \Cdots[shorten-end=10pt] & \Hspace*{5pt} \\
\brac{\beta}{1}_{\phi} & \beta^q & -\brac{\beta}{1}_{\phi^*}^{q^2} & -\brac{\beta}{2}_{\phi^*}^{q^3} & \Ddots[shorten-end=15pt] & \Hspace*{5pt} \\
\brac{\beta}{2}_{\phi} & \brac{\beta}{1}_{\phi}^q & \beta^{q^2} & -\brac{\beta}{1}_{\phi^*}^{q^3} & \Ddots[shorten-end=15pt] & \Hspace*{5pt} \\
\brac{\beta}{3}_{\phi} & \brac{\beta}{2}_{\phi}^q & \brac{\beta}{1}_{\phi}^{q^2} & \beta^{q^3} & \Ddots[shorten-end=15pt] & \Hspace*{5pt} \\
\Vdots[shorten-end=3pt] & \Ddots[shorten-end=30pt] & \Ddots[shorten-end=40pt] & \Ddots[shorten-end=30pt] & \Ddots[shorten-end=20pt] & \Hspace*{5pt} \\
\end{pNiceMatrix}
\quad \in \Mat_{\NN \times \NN}(\oK).
\end{equation}
Similarly, let $\bP_{\ell}(\beta) \in \Mat_{r\ell \times r\ell}(\oK)$ be the truncation of $\bP(\beta)$ after $r\ell$ rows and columns. For $a \in A$, set $\langle a/\wp \rangle = (a \bmod \wp)/\wp \in k$.  We then obtain the following matrix duality between Gauss-Thakur sums and their duals (Theorem~\ref{T:GaussSumduality1}).
\end{subsec}

\begin{theoremintro} \label{T:GaussSumduality1Intro}
With the matrices defined as above, we have the following identities.
\begin{alphenumerate}
\item For $a$, $b \in A$,
\[
\Gamma^*(\phi^*;a) \cdot \bQ(\wp') \cdot \Gamma(\phi;b)=\bP (\langle ab/\wp \rangle ).
\]
\item For $a$, $b \in A$ and $\ell \in \NN$,
\[
\Gamma_\ell^*(\phi^*;a) \cdot \bQ(\wp') \cdot \Gamma_\ell (\phi;b) = \bP_\ell (\langle ab/\wp \rangle ).
\]
\item In particular,
\[
\Gamma^*(\phi^*) \cdot \bQ(\wp') \cdot \Gamma(\phi) = \bP(1/\wp).
\]
\end{alphenumerate}
\end{theoremintro}

\begin{subsec}{Carlitz module}
In \S\ref{S:CarlitzApp}, we investigate the case of the Carlitz module $\sC$, which is the rank~$1$ Drinfeld module defined by $\sC_{\theta} = \theta + \tau$. The Gauss-Thakur sums $g_{\sC}(\chi)$ for $\sC$ are then precisely the ones designed by Thakur~\cite{Thakur88}. Thakur proved a number of results about $g_{\sC}(\chi)$, including identities and $\wp$-divisibility properties \cite{Thakur88}*{Thms.~I--IV}, which we summarize in Theorem~\ref{T:ThakurI-IV}. One consequence of Theorem~\ref{T:GaussSumduality1Intro} for $\sC$ is that
\begin{equation} \label{E:ggstaridentityintro}
g_{\sC}(\chi_{\varepsilon})g_{\sC}^*(\chi_{\varepsilon}) = \frac{1}{(\theta - \varepsilon) \wp'(\varepsilon)^2},
\end{equation}
and using this identity we can transform Thakur's results into corresponding results for Carlitz dual Gauss-Thakur sums (Theorem~\ref{T:dualThakurI-IV}). Thakur also proved result on the $\infty$-adic sign of $g_{\sC}(\chi)$ in~\cite{Thakur93a}*{Thm.~3.2}, and we prove a corresponding result for $g_{\sC}^*(\chi)$ (Proposition~\ref{P:gg*signs}). We prove the following application of Theorem~\ref{T:GaussSumduality1Intro} (Theorem~\ref{T:detGamma.Carlitz}).
\end{subsec}

\begin{theoremintro} \label{T:detGamma.CarlitzIntro}
Let $\wp \in A_+$ be irreducible of degree $d$, and fix a root $\varepsilon \in \FF_{q^d}$ of $\wp$. Let $K = k(\be_{\sC}(1/\wp))$ and $H=\FF_{q^d}K$. 
Define $\Gamma_d(\sC)$, $\Gamma_d^*(\sC^*) \in \GL_d(H)$ as above.
\begin{alphenumerate}
\item We have in $H$,
\[
\det \Gamma_d(\sC)^{q-1} = -\wp \quad \text{and} \quad \det \Gamma_d^*(\sC^*)^{q-1} = -\frac{1}{\wp^{q^d}}.
\]
\item In the completion $H_{\fp_0}$, where $\fp_0$ is the unique prime of $H$ above $\theta-\varepsilon \in \FF_{q^d}[\theta]$, by choosing a root $(-\wp)^{1/(q^d-1)}$ as in Theorem~\ref{T:ThakurI-IV}(c),
\[
\det \Gamma_d(\sC) = (-\wp)^{1/(q-1)} \quad \text{and} \quad
\det \Gamma_d^*(\sC^*) = \frac{(-1)^d}{\rN(\wp'(\varepsilon)) \cdot (-\wp)^{q^d/(q-1)}},
\]
where $\rN : \FF_{q^d}^{\times} \to \FF_q^{\times}$ is the norm.
\end{alphenumerate}
\end{theoremintro}

We prove $\infty$-adic identities for $\det \Gamma_d(\sC)$ and $\det \Gamma_d^*(\sC^*)$ in Theorem~\ref{T:detGamma.Carlitz}(c), which we have not repeated here for brevity, but can be obtained from Theorem~\ref{T:detGamma.GeneralIntro} below as well.

As another application of Theorem~\ref{T:GaussSumduality1Intro} for the Carlitz module, we obtain identities for Goss-Dirichlet $L$-functions,
\[
L(\chi_k,s) = \sum_{a \in A_+} \frac{\chi_k(a)}{a}
\]
by combining Theorem~\ref{T:GaussSumduality1Intro} with Pellarin's special $L$-value formula~\cite{Pellarin12}*{Thm.~1} and an identity of Gazda and Maurischat~\cite{GazdaMaurischat21}*{Thm.~6.1}. We obtain a summation identity (Proposition~\ref{P:Lgtrace}) between $L$-values, twists of Gauss-Thakur sums, the Carlitz period~$\tpi$, and Carlitz multiplication coefficients,
\begin{equation}
\sum_{k=0}^{d-1} L(\chi_k,q^{i}) g^{[j]}(\chi_k;b) =
\begin{cases}
{\displaystyle \tpi^{q^{i}}\brac{ \langle b/\wp \rangle }{i-j}_{\sC}^{q^{j}}} & \textup{if $i-j \geq 0$,} \\[10pt]
{\displaystyle -\tpi^{q^{i}}\brac{ \langle b/\wp \rangle }{j-i}_{\sC^*}^{q^{j}}} & \textup{if $i-j < 0$.}
\end{cases}
\end{equation}
We explore also connections with Anderson's log-algebraicity~\cite{And96} in Corollary~\ref{C:logalg}.

\begin{subsec}{Determinants in higher rank}
In \S\ref{S:DrinfeldApps}, we investigate the determinants of $\Gamma_d(\phi)$ and $\Gamma^*(\phi^*)$ for general~$\phi$. As an intermediate result we prove (Proposition~\ref{P:detfGq-1}) that
\[
\prod_{j=0}^{d-1}\det \fG_\phi(\chi_j)^{q-1} = \frac{(-1)^{rd}\wp}{b_r^d}
\quad \text{and} \quad
\prod_{j=0}^{d-1}\det \fG^*_{\phi^*}(\chi_j)^{q-1} = \frac{(-1)^{rd}b_r^{d}}{\wp^{q^r}},
\]
which serve as higher rank versions of a result of Thakur~\cite{Thakur88}*{Thm.~II} (see Theorem~\ref{T:ThakurI-IV}(b)). The full result is the following (Theorem~\ref{T:detGamma.General}).
\end{subsec}

\begin{theoremintro} \label{T:detGamma.GeneralIntro}
Let $\wp \in A_+$ be irreducible of degree $d$, and fix a root $\varepsilon \in \FF_{q^d}$ of $\wp$. Define $\Gamma_d(\phi)$, $\Gamma_d^*(\phi^*) \in \GL_d(\oK)$ as in~\eqref{E:Gammadef}--\eqref{E:Gamma*def}.
\begin{alphenumerate}
\item We have
\[
\det \Gamma_d(\phi)^{q-1} = \frac{(-1)^{r}\wp}{b_r^{(q^{rd}-1)/(q^r-1)}}
\quad \text{and} \quad
\det \Gamma_d^*(\phi^*)^{q-1} = \frac{(-1)^{r}b_r^{(q^{rd}-1)/(q^r-1)}}{\wp^{q^{rd}}}.
\]
\item By choosing $(-1)^{1/(q-1)}$ to have the same sign as $\tpi$ and $\wp^{1/(q-1)}$ to have sign~$1$,
\[
\det \Gamma_d(\phi) = \frac{(-1)^{rd(d-1)/2} \cdot (-1)^{d/(q-1)} \cdot \mu^{(q^{rd}-1)/(q^r-1)}}{\rV(\varepsilon)^r} \cdot \wp^{1/(q-1)}
\]
and
\[
\det \Gamma_d^*(\phi^*) = \frac{1}{(-1)^{d/(q-1)} \cdot \mu^{(q^{rd}-1)/(q^r-1)} \cdot \rV(\varepsilon)^r}\cdot \frac{1}{\wp^{q^{rd}/(q-1)}},
\]
where $\mu \in \oK$ satisfies $\det \Upsilon = \mu \omega_{\sC}$ and $\rV(\varepsilon)$ is the Vandermonde determinant $V(\varepsilon, \dots, \varepsilon^{q^{d-1}})$.
\end{alphenumerate}
\end{theoremintro}

\begin{subsec}{Drinfeld modular forms}
We apply Theorem~\ref{T:detGamma.GeneralIntro} to Drinfeld modular forms of higher rank with values in Tate algebras, building on identities of Gekeler~\cite{Gekeler17}, Chen and Gezmi\c{s}~\cites{ChenGezmis22, ChenGezmis26}, Pellarin~\cite{Pellarin12}, and Perkins~\cite{Perkins14}. In particular, we show (Corollary~\ref{C:DMFs}) that Gekeler's modular forms $h_r(\bz)$ and $\Delta_r(\bz)$, for $\bz$ in the Drinfeld upper half space $\HH^r$, are expressible as determinants of dual Gauss-Thakur sums. For example, we find
\[
\Delta_r(\bz)^{(q^{rd}-1)/(q^r-1)} = (-1)^r \wp^{q^{rd}} \cdot \det \Gamma_d^*(\bz)^{q-1},
\]
where $\Gamma_d^*(\bz) = \Gamma_d^*(\phi^{\bz})$ is associated to the Drinfeld module $\phi^{\bz}$. See \S\ref{SS:DMFspecialvalues} for details.
\end{subsec}

\begin{acknowledgments}
The authors thank G.~H.\ Ferraro for many helpful discussions on the contents of this paper.
\end{acknowledgments}

\section{Drinfeld modules and their adjoints} \label{S:Drinfeld}

\begin{subsec}{Rigid analytic functions} \label{SS:RAfns}
Let $t$ be a variable independent of $\theta$. For $\alpha \in \C$, $\alpha \neq 0$, the \emph{Tate algebra}
\[
\TT_{\alpha} \assign \biggl\{ \sum c_m t^m \in \power{\C}{t} : \lim_{m \to \infty} \inorm{\alpha}^m \inorm{c_m} = 0 \biggr\},
\]
consists of power series that converge on the closed disk of radius $\inorm{\alpha}$. The Gauss norm on~$\TT_{\alpha}$, defined as $\dnorm{\sum c_m t^m} = \max_m \{ \inorm{\alpha}^m \inorm{c_m} \}$, makes $\TT_{\alpha}$ into a complete $\C[t]$-algebra. We let $\TT \assign \TT_1$. We further define rings of double-tailed Laurent series,
\[
\SSS_{\alpha} \assign \biggl\{ \sum_{m \in \ZZ} c_m t^m \in \power{\C}{t,t^{-1}} : \lim_{m \to \pm \infty} \inorm{\alpha}^m \inorm{c_m} = 0 \biggr\},
\]
which provide functions that converge on the annular strip $\{ x \in \C : \inorm{x} = \inorm{\alpha} \}$ (an affinoid ring domain in the sense of~\cite{FresnelvdPut}*{\S 2.2}).
\end{subsec}

\begin{subsec}{Twisted polynomials}
We let $\tau : \C \to \C$ denote the $q$-th power Frobenius endomorphism, which we extend to matrices $\tau : \Mat_{m \times n}(\C) \to \Mat_{m\times n}(\C)$ entry-wise. For an $\FF_q$-algebra $R \subseteq \Mat_{n\times n}(\C)$, we let $R[\tau]$ be the ring of twisted polynomials in $\tau$ over $R$, i.e., the $\FF_q$-subalgebra of $\End_{\FF_q}(\C)$ generated by $R$ and $\tau$. Likewise, if $R \subseteq \Mat_{n\times n}(\C)$ is a perfect $\FF_q$-subalgebra, we similarly define twisted polynomials $R[\tau^{-1}]$ in $\End_{\FF_q}(\C)$. We construct the ring of Laurent polynomials $R[\tau,\tau^{-1}]$ as well.
\end{subsec}

\begin{subsec}{The $*$-operator}
If $R$ is a perfect field, then the $\FF_q$-linear involution $*: R[\tau,\tau^{-1}] \to R[\tau,\tau^{-1}]$ defined by
\begin{equation}
\biggl( \sum b_i \tau^i \biggr)^* \assign \sum b_i^{1/q^i} \tau^{-i} = \sum \tau^{-i} b_i,
\end{equation}
satisfies the anti-commutative property, $(\alpha \beta)^* = \beta^* \alpha^*$, for $\alpha, \beta \in R[\tau,\tau^{-1}]$.
As such~$*$ induces anti-isomorphisms, $* : R[\tau] \iso R[\tau^{-1}]$ and $* : R[\tau^{-1}] \iso R[\tau]$.
We also extend~$*$ to Laurent series rings $* : \laurent{R}{\tau} \to \laurent{R}{\tau^{-1}}$ and $* : \laurent{R}{\tau^{-1}} \to \laurent{R}{\tau}$ in the natural way.
For more information, see for example \citelist{\cite{Goss}*{\S 4.14} \cite{Poonen96} \cite{Thakur}*{\S 2.10}}.
\end{subsec}

\begin{subsec}{Drinfeld modules} \label{SS:Drinfeldmodules}
For an intermediate field $k \subseteq K \subseteq \C$, a \emph{Drinfeld module over~$K$} is defined by an $\FF_q$-algebra homomorphism $\phi : A \to K[\tau]$ such that, as in~\eqref{E:rhodef}, for $r \geq 1$, we have $\phi_{\theta} = \theta + b_1\tau + \cdots + b_r \tau^r$, $b_r \neq 0$.
In this case $r$ is the \emph{rank} of $\phi$. See \cites{Goss, Papikian, Thakur} for fundamental properties of Drinfeld modules. For $a \in A$, we write
\[
\phi_a = \sum_{n=0}^{r \deg a} \brac{a}{n}_{\phi} \tau^n, \quad \brac{a}{n}_{\phi} \in K,
\]
Then $\phi$ induces an $A$-module structure on $\C$ such that for $a \in A$, $x \in \C$, we have $\phi_a(x) = \sum_{n=0}^{r \deg a} \brac{a}{n}_{\phi} x^{q^n}$.
The exponential $\exp_\phi \in 1 + \power{K}{\tau} \tau$ and logarithm $\log_{\phi} \in 1 + \power{K}{\tau}\tau$ are the unique series such that
\begin{gather} \label{E:explog}
{\exp_{\phi}} \cdot {\log_\phi} = {\log_\phi} \cdot {\exp_\phi} = 1, \\
{\exp_\phi} \cdot {a} = {\phi_a} \cdot {\exp_{\phi}}, \quad
{\log_\phi} \cdot {\phi_a} = {a} \cdot {\log_{\phi}}, \quad \forall\,a \in A. \notag
\end{gather}
Both $\exp_\phi$ and $\log_\phi$ induce functions on $z \in \C$, where $\exp_\phi(z)$ is an entire function and $\log_\phi(z)$ has a finite radius of convergence (e.g., see \citelist{\cite{Goss}*{Prop.~4.14.2} \cite{KP23}*{Cor.~4.5} \cite{Papikian}*{Lem.~5.1.5}}). If $\exp_{\phi} = \sum_{i \geq 0} \alpha_i \tau^i$ and $\log_{\phi} = \sum_{i \geq 0} \beta_i \tau^i$, then for a variable $x$,
\begin{equation} \label{E:bracxdef}
{\exp_{\phi}} \cdot {x} \cdot {\log_{\phi}} = \sum_{k \geq 0} \brac{x}{k}_{\phi} \tau^k
\quad \Rightarrow \quad
\brac{x}{n}_{\phi} = \sum_{i=0}^n \alpha_i \beta_{n-i}^{q^i}\, x^{q^i} \in K[x].
\end{equation}
We note that if $a \in A$, then as one would expect, $\brac{x}{n}_{\phi}\big|_{x=a} = \brac{a}{n}_{\phi}$.
\end{subsec}

\begin{subsec}{Adjoint Drinfeld modules} \label{SS:adjoints}
When $K$ is perfect, the \emph{adjoint} of $\phi$ is the $\FF_q$-algebra homomorphism $\phi^* : A \to K[\tau^{-1}]$ defined by
\[
\phi^*_a \assign (\phi_a)^* = \sum_{n=0}^{r \deg a} \brac{a}{n}_{\phi}^{1/q^n} \tau^{-n} \rassign \sum_{n=0}^{r \deg a} \tau^{-n} \brac{a}{-n}_{\phi^*}.
\]
Note that $\brac{a}{n}_{\phi} = \brac{a}{-n}_{\phi^*}$ for $a \in A$ and $n \in \ZZ$. Similar to the Drinfeld module case, $\phi^*$ induces an $A$-module structure on $\C$ such that for $a \in A$ and $x \in \C$,
\[
\phi^*_a(x) = \sum_{n=0}^{r \deg a} \brac{a}{n}_{\phi}^{1/q^n} x^{1/q^n} = \sum_{n=0}^{r \deg a} \brac{a}{-n}_{\phi^*}^{1/q^n} x^{1/q^n}.
\]
For basic properties of adjoint modules, see \citelist{\cite{Ferraro25} \cite{Goss}*{\S 4.14} \cite{HuangP25}*{\S 2.2, \S 3.2} \cite{Poonen96}}. By applying the $*$ operator to~\eqref{E:explog}, we obtain the companion formulas for $\exp_{\phi}^*$, $\log_{\phi}^* \in 1 + \power{K}{\tau^{-1}} \tau^{-1}$,
\begin{gather} \label{E:explogstar}
{\log_{\phi}^*} \cdot {\exp_\phi^*} = {\exp_\phi^*} \cdot {\log_\phi^*} = 1, \\
{\log_\phi^*} \cdot {a} = {\phi_a^*} \cdot {\log_{\phi}^*}, \quad
{\exp_\phi^*} \cdot {\phi_a^*} = {a} \cdot {\exp_{\phi}^*}, \quad \forall\,a \in A. \notag
\end{gather}
For $z\in \C$, $\exp_{\phi}^*(z)$ converges everywhere, but $\log_{\phi}^*(z)$ converges only at $z=0$.
\end{subsec}

\begin{subsec}{Torsion submodules and exponentiation} \label{SS:torsexp}
For $b \in A$, $b \neq 0$, we have
\[
\phi[b] = \{ x \in \oK \mid \phi_b(x) = 0 \} \cong (A/bA)^r, \quad
\phi^*[b] = \{ x \in \oK \mid \phi^*(x) = 0 \} \cong (A/bA)^r,
\]
(see \citelist{\cite{Goss}*{\S 4.14} \cite{Poonen96}*{\S 8}}).
The kernel $\Lambda_{\phi}$ of $\exp_{\phi}$ is a discrete free $A$-submodule of $\C$ of rank~$r$, and we have an induced isomorphism of $A$-modules, $\exp_\phi : \C/\Lambda_\phi \iso \phi(\C)$, where $\phi(\C)$ denotes $\C$ with $A$-module structure via $\phi$. This implies that
\begin{equation} \label{E:torsionexp}
\phi[b] = \exp_{\phi} \bigl( b^{-1} \Lambda_{\phi} / \Lambda_{\phi} \bigr).
\end{equation}
The corresponding map $\exp_{\phi}^* : \C \to \phi^*(\C)$ is a bit more subtle but has been studied extensively by Ferraro~\cite{Ferraro25} and Poonen~\cite{Poonen96}. We see from~\eqref{E:explogstar} that $\phi^*_{\tors} \subseteq \ker \exp_{\phi}^*$, where $\phi^*_{\tors}$ is the $A$-torsion submodule of $\phi^*(\C)$, but it is not an equality. Indeed $\ker \exp_{\phi}^*$ is compact and contains $\phi^*_{\tors}$ as a dense subset~\cite{Poonen96}*{Prop.~20}.

The analogous identity to~\eqref{E:torsionexp} for $\phi^*[b]$ would more naturally, given the identities in~\eqref{E:explogstar}, arise from $\log_{\phi}^*$. However, since $\log_{\phi}^*$ converges only at $0$, this does not produce results. Nevertheless, we develop a stand-in for \eqref{E:torsionexp} using techniques from~\cite{ABP04}*{\S 5}, which we explore in \S\ref{S:emaps} (see Lemma~\ref{L:e*propertiesfortorsion}).
\end{subsec}

\begin{subsec}{Extension $\phi : k_{\infty} \to \laurent{\C}{\tau}$} \label{SS:formalphi}
Using \eqref{E:bracxdef} we extend $\phi$ to an injective $\FF_q$-algebra homomorphism $\phi: k_{\infty} \to \laurent{\C}{\tau}$, setting
\[
\phi_{\alpha} \assign \sum_{n=0}^{\infty} \brac{\alpha}{n}_{\phi} \tau^n, \quad \forall\,\alpha \in k_{\infty}.
\]
In this way we obtain a formal Drinfeld $k_{\infty}$-module \citelist{\cite{Goss}*{\S 4.6} \cite{Papikian}*{\S 6.5} \cite{Rosen03}}. In particular, $\brac{\alpha}{0}_{\phi} = \alpha$, for all $\alpha \in k_{\infty}$. For convenience we set $\brac{\alpha}{n}_{\phi}=0$ if $n < 0$. We note that $\phi_{\theta\alpha} = \phi_{\theta}\phi_{\alpha} = \phi_{\alpha}\phi_{\theta}$ implies for $\alpha \in k_{\infty}$ and $n \in \ZZ$ (using $b_0=\theta$),
\begin{equation} \label{E:phibracrec}
\brac{\theta\alpha}{n}_{\phi}
= \sum_{j=0}^r b_j \brac{\alpha}{n-j}_{\phi}^{q^j}
= \sum_{j=0}^r b_j^{q^{n-j}} \brac{\alpha}{n-j}_{\phi}.
\end{equation}
As $\brac{x}{n} \in K[x]$ by \eqref{E:bracxdef}, each coefficient function of $\phi : k_{\infty} \to \laurent{\C}{\tau}$ is continuous.
\end{subsec}

\begin{subsec}{Extension $\phi^* : k_{\infty} \to \laurent{\C}{\tau}$} \label{SS:phistarext}
As defined by Ferraro~\cite{Ferraro25}*{Prop.~7.16}, we extend $\phi^*$ to $\phi^*: k_{\infty} \to \laurent{\C}{\tau}$. The definition here matches Ferraro's $\hat{\Phi}$, though as we are working with $A=\FF_q[\theta]$, we obtain explicit identities as well. For $a \in A$ of degree $d$,
\begin{equation} \label{E:InvertingPrep}
 \tau^{rd}\phi_{a}^* = \sum_{n=0}^{rd} \brac{a}{-n}_{\phi^*}^{q^{rd-n}} \tau^{rd-n} = \sum_{n=0}^{rd} \brac{a}{n-rd}_{\phi^*}^{q^n} \tau^n,   
\end{equation}
so then $\tau^{rd} \phi_{a}^*$ is invertible in $\power{\C}{\tau}$, as $\brac{a}{-rd}_{\phi^*} \neq 0$. Thus in $\laurent{\C}{\tau}$ we find
\begin{equation} \label{E:SettingUpInverse}
\bigl(\phi_{a}^*\bigr)^{-1} = \brac{a}{-rd}_{\phi^*}^{-1} \tau^{rd} - \brac{a}{-rd}_{\phi^*}^{-q-1} \brac{a}{-rd+1}_{\phi^*}^{q} \tau^{rd+1} + \cdots.
\end{equation}
Setting
\[
\phi_{1/a}^* \assign \bigl(\phi_{a}^*\bigr)^{-1} \quad \in \laurent{\C}{\tau},
\]
we obtain a well-defined $\FF_q$-algebra homomorphism $\phi^* : k \to \laurent{\C}{\tau}$, such that $\phi_{a/b}^* = \phi_a^* \cdot \phi_{1/b}^*$, for $a$, $b \in A$, $b \neq 0$ (cf.~\cite{Goss}*{Prop.~5.5.3}).
Furthermore, since for $\alpha \in k$ we have $\phi_{\alpha}^* \in \power{\C}{\tau} \tau^{-r\deg \alpha}$, it follows that the coefficient functions of $\phi^*$ are continuous. We thus define $\brac{\alpha}{n}_{\phi^*}$ for $\alpha \in k_{\infty}$ and $n \in \ZZ$ by
\[
\phi_{\alpha}^* \rassign \sum_{n=-r\deg \alpha}^\infty \brac{\alpha}{n}_{\phi^*}^{q^n} \tau^n = \sum_{n=-r \deg \alpha}^{\infty} \tau^n \brac{\alpha}{n}_{\phi^*}.
\]
By convention we set $\brac{\alpha}{n}_{\phi^*}=0$ for $n < -r \deg \alpha$. 

Since $\phi_{\theta\alpha}^* = \phi_{\alpha}^*\phi_\theta^* = \phi_{\theta}^*\phi_{\alpha}^*$ for $\alpha \in k_{\infty}$, similar to~\eqref{E:phibracrec} we have for $n \in \ZZ$,
\begin{equation} \label{E:phistarbracrec}
\brac{\theta\alpha}{n}_{\phi^*}
= \sum_{j=0}^r b_j \brac{\alpha}{n+j}_{\phi^*}^{q^j}
= \sum_{j=0}^r b_j^{q^{-n-j}} \brac{\alpha}{n+j}_{\phi^*}.
\end{equation}
We note that if $\alpha = \sum_{m\geq m_0} c_m \theta^{-m} \in k_{\infty}$, $c_m \in \FF_q$, then
\[
  \phi_{\alpha}^* = \sum_{m=m_0}^{\infty} c_m \phi_{\theta^{-m}}^* = \sum_{m=m_0}^\infty c_m \sum_{n=rm}^{\infty}  \tau^n \brac{\theta^{-m}}{n}_{\phi^*}.  
\]
A word of caution is that generally $\phi_{\alpha}^* \neq ( \phi_{\alpha} )^*$ for $\alpha \in k_{\infty}$, unless $\alpha \in A$.
If we had extended $\phi^*$ to a homomorphism $k_{\infty} \to \laurent{\C}{\tau^{-1}}$, arguing as in \S\ref{SS:formalphi}, such an identity would hold for all $\alpha \in k_{\infty}$. However, the definition of $\phi^* : k_{\infty} \to \laurent{\C}{\tau}$ here is the right one for our future considerations (e.g., see Proposition~\ref{P:PsiNUpsilon11} and Theorem~\ref{T:GaussSumduality1}).
\end{subsec}

\begin{lemma} \label{L:qPowerBracketStarPoly}
Let $a$, $b \in A$, with $b \neq 0$. Then $\brac{a/b}{m}_{\phi^*}^{q^m}\in A[b_1,\dots,b_r,b_r^{-1}]$ for all $m \geq 0$.
\end{lemma}

\begin{proof}
Let $R = A[b_1, \dots, b_r,b_r^{-1}]$. 
For $b \in A$, recall that for $n \in \ZZ$, $\brac{b}{-n}_{\phi^*} = \brac{b}{n}_{\phi}$, and so $\brac{b}{-n}_{\phi^*} \in A[b_1, \dots, b_r]$. Furthermore, if $\deg b = d$, then $\brac{b}{-rd}_{\phi^*} = b_r^{1+q+\dots+q^{\deg a-1}} \in R^{\times}$. Therefore by \eqref{E:InvertingPrep}, we have $\tau^{rd} \phi_b^* \in \power{R}{\tau}^{\times}$, and so $\phi_{1/b}^* = (\tau^{rd} \phi_b^*)^{-1} \tau^{rd} \in \power{R}{\tau}$. Since $\phi_{1/b}^* = \sum_{m \geq 0} \tau^m \brac{1/b}{m}_{\phi^*}$, the result holds for $1/b$. In general,
\[
\phi_{a/b}^* = \sum_{k= r(\deg b - \deg a)}^{\infty} \tau^k \brac{a/b}{k}_{\phi^*} = \phi_{1/b}^* \cdot \phi_a^*
= \sum_{m=r \deg b}^{\infty} \sum_{n=0}^{r \deg a} \tau^{m-n} \brac{1/b}{m}_{\phi^*}^{q^n} \brac{a}{-n}_{\phi^*}.
\]
For $m-n \geq 0$, when we raise a given coefficient in the final sum to the $q^{m-n}$-th power, the result follows from the result for $1/b$ and the fact that $\brac{a}{-n}_{\phi^*} \in A[b_1, \dots, b_r]$.
\end{proof}

\section{Rigid analytic trivializations} \label{S:RAT}

For $i \in \ZZ$, we define Frobenius twisting
\[
f = \sum c_m t^m \mapsto f^{(i)} \assign \sum c_m^{q^i} t^m \;\; : \;\; \laurent{\C}{t} \to \laurent{\C}{t},
\]
which we extend to matrices in $\laurent{\C}{t}$ entry-wise. Frobenius twisting induces automorphisms of the Tate algebra $\TT$. For an $\FF_q[t]$-subalgebra $S$ of $\laurent{\C}{t}$, we define twisted polynomials $S[\tau, \tau^{-1}]$ subject to the relation, $\tau f = f^{(1)} \tau$, $f \in S$. In this way, for $\eta = \sum B_i \tau^i \in \Mat_{d}(S)[\tau,\tau^{-1}]$ and $F \in \Mat_{d \times n}(\laurent{\C}{t})$, we have the operation
\[
\eta(F) \assign \sum_i B_i F^{(i)} \quad \in \Mat_{d \times n}\bigl( \laurent{\C}{t} \bigr).
\]

\begin{subsec}{Special functions} \label{SS:special}
Fix a Drinfeld module $\phi$ of rank $r$ defined over $\C$ as in~\eqref{E:rhodef}. A \emph{special function} and a \emph{dual special function} for $\phi$ are $f$ and $g \in \TT$ respectively satisfying
\[
\phi_a(f) = a(t) f, \quad \phi_a^*(g) = a(t)g, \quad \forall\,a \in A.
\]
Special functions for Drinfeld modules and $t$-modules have been studied extensively in various contexts, especially related to Anderson generating functions (e.g., see \cites{CP12, EP14, Ferraro22, Ferraro25, GazdaMaurischat21, Goss94, KP23, Maurischat22, NP24, Pellarin08}). Dual special functions arise via Frobenius difference equations with respect to $\tau^{-1}$ (e.g., see \cites{ABP04, CP12, Pellarin08, NP24}), but more recently Ferraro~\cite{Ferraro25} has provided new insights and directions for them (e.g., see Theorems~\ref{T:FerraroUniv} and~\ref{T:Ferraropairing}).
\end{subsec}

\begin{subsec}{Anderson generating functions} \label{SS:AGF}
For $u \in \C$, we let
\[
\omega_u \assign \omega_u(t) \assign \sum_{m=0}^{\infty} \exp_{\phi} \biggl( \frac{u}{\theta^{m+1}} \biggr) t^m \quad \in \TT.
\]
Defined originally by Anderson (see~\cite{Goss94}), the function~$\omega_u$ satisfies the principal functional equation $\phi_t( \omega_u) = t\omega_u + \exp_{\phi}(u)$. For $\lambda \in \Lambda_{\phi}$, we then find that
\begin{equation} \label{E:AGFpereq}
\phi_a(\omega_{\lambda}) = a(t) \omega_{\lambda}, \quad \forall\,a \in A.
\end{equation}
Thus $\omega_{\lambda}$ is a special function for $\phi$ as defined above, and moreover, all special functions for $\phi$ arise in this way~\cite{Pellarin08}*{\S 4.2} (see also~\cite{EP14}*{Rem.~6.3}). Pellarin~\cite{Pellarin08}*{\S 4.2} showed further that $\omega_u$ has a partial fraction decomposition,
\[
\omega_u = \sum_{n=0}^{\infty} \frac{\alpha_n u^{q^n}}{\theta^{q^n}-t}, \quad u \in \C,
\]
where $\exp_{\phi} = \sum_{n\geq 0} \alpha_n \tau^n$, converging in $\TT$ with respect to the Gauss norm. And so $\res_{t=\theta} (\omega_u\, dt) = -u$. For extensions to higher dimensional $t$-modules see \citelist{\cite{Maurischat22} \cite{NP24}*{\S 4}}.
\end{subsec}

\begin{subsec}{Rigid analytic trivializations} \label{SS:RATs}
We set
\[
\Theta \assign \begin{pNiceMatrix}
0 & 1 & \Cdots[shorten=20pt] & 0 \\
\Vdots & \Vdots & \Ddots[shorten=20pt] & \Vdots \\
0 & 0 & \Cdots[shorten=20pt] & 1 \\
(t-\theta)/b_r & -b_1/b_r & \Cdots & -b_{r-1}/b_r
\end{pNiceMatrix}
\]
and
\[
\Phi \assign\begin{pNiceMatrix}
0 & 1 & \Cdots[shorten=30pt] & 0 \\
\Vdots & \Vdots & \Ddots[shorten=30pt] & \Vdots \\
0 & 0 & \Cdots[shorten=30pt] & 1 \\
(t-\theta)/b_r^{(-r)} & -b_1^{(-1)}/b_r^{(-r)} & \Cdots & -b_{r-1}^{(-r+1)}/b_r^{(-r)}
\end{pNiceMatrix},
\]
which represent multiplication by $\tau$ (respectively by $\tau^{-1}$) on the $t$-motive of $\phi$ (respectively the $t$-comotive of $\phi$). Though we will not pursue the motivic aspects in this paper directly, the reader is directed to \cites{CP12, GazdaMaurischat25, HartlJuschka20, NP24, Pellarin08} for more information.

We now fix an $A$-basis $\pi_1, \dots, \pi_r$ of $\Lambda_{\phi}$, and let $\omega_i \assign \omega_{\pi_i}$ for $i = 1, \dots, r$. Setting
\begin{equation} \label{E:Upsilondef}
\Upsilon \assign \begin{pNiceMatrix}
\omega_1 & \omega_1^{(1)} & \Cdots & \omega_1^{(r-1)}, \\
\Vdots & \Vdots & & \Vdots \\
\omega_r & \omega_r^{(1)} & \Cdots & \omega_r^{(r-1)}
\end{pNiceMatrix},
\end{equation}
Pellarin~\cite{Pellarin08}*{\S 4.2} proved that $\Upsilon \in \GL_r(\TT)$, and moreover we observe that
\begin{equation} \label{E:Upsilontwist}
\Upsilon^{(1)} = \Upsilon \Theta^{\tr},
\end{equation}
where $\Theta^{\tr}$ is the transpose of $\Theta$.
Otherwise, as in \citelist{\cite{CP12}*{\S 3.4} \cite{HartlJuschka20}*{\S 2.5.2} \cite{Pellarin08}}, if we let
\begin{equation} \label{E:Vdef}
V \assign \begin{pNiceMatrix}
b_1 & b_2^{(-1)} & \Cdots & b_{r-1}^{(-r+2)} & b_r^{(-r+1)} \\
b_2 & b_3^{(-1)} & \Cdots & b_r^{(-r+2)} & \\
\Vdots & \Vdots & \Iddots & & \\
b_{r-1} & b_r^{(-1)} & & &\\
b_r & & & & 
\end{pNiceMatrix},
\end{equation}
then it follows that $V^{(-1)} \Phi = \Theta^{\tr} V$. Setting
\begin{equation} \label{E:Psidef}
\Psi \assign -V^{-1} \bigl( \Upsilon^{(1)} \bigr)^{-1},
\end{equation}
we see that
\begin{equation} \label{E:Psitwist}
\Psi^{(-1)} = \Phi \Psi.
\end{equation}
It follows from \cite{ABP04}*{Prop.~3.1.3} that the entries of $\Psi$ are entire functions on $\C$, and we see that $\Psi$ has the form
\begin{equation} \label{E:zetadef}
\Psi = \begin{pNiceMatrix}
\zeta_1 & \Cdots & \zeta_r \\
\zeta_1^{(-1)} & \Cdots & \zeta_r^{(-1)} \\
\Vdots & & \Vdots \\
\zeta_1^{(-r+1)} & \Cdots & \zeta_r^{(-r+1)}
\end{pNiceMatrix},
\end{equation}
thus defining $\zeta_1, \dots, \zeta_r \in \TT$. Now \eqref{E:Psitwist} implies $\zeta_1, \dots, \zeta_r$ are dual special functions for~$\phi$, and moreover they form an $\FF_q[t]$-basis of all dual special functions~\cite{P08}*{Prop.~3.4.7(b)}.
\end{subsec}

\begin{subsec}{The matrix $V$}
The matrix in \eqref{E:Vdef} represents the Hartl-Juschka pairing between the $t$-comotive of $\phi$ and its module of biderivations~\cite{HartlJuschka20}*{\S 2.5.2} (see also \cite{NP24}*{\S 4} for additional constructions). We find through calculation~\citelist{ \cite{HartlJuschka20}*{Eq.~(2.5.12)} \cite{NP24}*{Ex.~4.117}},
\begin{equation} \label{E:VtwistPhi}
V^{(-1)} \Phi = \Theta^{\tr} V = \begin{pNiceMatrix}
t-\theta & 0 & \Cdots & & 0 \\
0 & b_2^{(-1)} & \Cdots &  & b_r^{(-r+1)} \\
 & &  &  & 0\\
\Vdots & \Vdots & \Iddots & \Iddots[draw-first] & \Vdots \\
0 & b_r^{(-1)} & 0 & \Cdots & 0
\end{pNiceMatrix}.
\end{equation}
From \eqref{E:Upsilontwist} and~\eqref{E:Psidef}, we see that $\Psi\Upsilon = -V^{-1} (\Theta^{\tr})^{-1}$. Combining this with the preceding identity yields the following.
\end{subsec}

\begin{lemma} \label{L:PsiUpsilon}
We have
\[
\Psi \Upsilon = -\bigl( \Theta^{\tr} V \bigr)^{-1} = -\bigl( V^{(-1)} \Phi \bigr)^{-1} =
\begin{pNiceMatrix}
\dfrac{1}{\theta - t} & 0 & \Cdots & 0 \\
0 & 0 & \Cdots & \dfrac{-1}{b_r^{(-1)}} \\
\Vdots & \Vdots & \Iddots & \Vdots \\
0 & \dfrac{-1}{b_r^{(-r+1)}} & \Cdots & *
\end{pNiceMatrix}.
\]
\end{lemma}

\begin{remark}
(1) All the identities to this point have depended on our choices of the defining homomorphism for $\phi$ and then the $A$-basis $\pi_1, \dots, \pi_r$ for $\Lambda_\phi$. However, these are the only significant choices we make, and we design our future constructions to be consistent with them. (2) The inclusion of the negative sign in \eqref{E:Psidef} makes $\Psi$ differ slightly from previous incarnations, such as in~\cites{CP12, NP24}, but it aligns better with Ferraro~\cite{Ferraro25}.
\end{remark}

For $N \in \ZZ$, using~\eqref{E:Psitwist} we note that $\Psi^{(N)}\Psi^{-1} = (\Phi^{-1})^{(N)} \cdots (\Phi^{-1})^{(1)}$ when $N \geq 0$, and $\Psi^{(N)}\Psi^{-1} = \Phi^{(N+1)} \cdots \Phi^{(-1)} \Phi$ when $N < 0$.
Then applying Lemma~\ref{L:PsiUpsilon},
\begin{equation} \label{E:PsiNUpsilon11}
\bigl[ \Psi^{(N)} \Upsilon \bigr]_{11} = \frac{1}{\theta-t} \cdot \begin{cases}
\Bigl[ \bigl( \Phi^{-1} \bigr)^{(N)} \cdots \bigl( \Phi^{-1} \bigr)^{(1)} \Bigr]_{11}, & \textup{if $N \geq 0$,} \\[10pt]
\Bigl[ \Phi^{(N+1)} \cdots \Phi^{(-1)} \Phi \Bigr]_{11}, & \textup{if $N < 0$.}
\end{cases}
\end{equation} 
This yields the following proposition, which serves as a generalization of~\cite{ABP04}*{Lem.~5.4.1} to arbitrary Drinfeld $A$-modules. One should also compare with~\cite{Ferraro25}*{Thm.~7.18}.

\begin{proposition} \label{P:PsiNUpsilon11}
For $N \in \ZZ$,
\[
\bigl[ \Psi^{(N)} \Upsilon \bigr]_{11} = \begin{cases}
{\displaystyle \sum_{m=0}^{\infty} \brac{\theta^{-m-1}}{N}_{\phi} \, t^m,} & \textup{if $N \geq 0$,} \\[15pt]
{\displaystyle -\sum_{m=0}^{\infty} \brac{\theta^{-m-1}}{-N}_{\phi^*}\, t^m,} & \textup{if $N < 0$.}
\end{cases}
\]
\end{proposition}

\begin{proof}
Let $\cS_N(t)$ denote the left-hand side of the identity in the proposition, and likewise let $\cT_N(t)$ denote the right-hand side. We let $\bszeta_{\phi} \assign (\zeta_1, \dots, \zeta_r) \in \Mat_{1 \times r}(\TT)$ and $\bsomega_{\phi} \assign (\omega_1, \dots, \omega_r)^{\tr} \in \Mat_{r \times 1}(\TT)$, which are the first row of $\Psi$ and first column of $\Upsilon$ respectively.
We note that the $r$ cases where $-r+1 \leq N \leq 0$ are precisely obtained by considering the first column of $\Psi\Upsilon$ in Lemma~\ref{L:PsiUpsilon}. Indeed in these cases,
\begin{equation} \label{E:SNbase}
\cS_N(t) = [ \Psi\Upsilon]_{-N+1,1} = \begin{cases}
{\displaystyle \frac{1}{\theta-t},} & \textup{if $N=0$,} \\[10pt]
0 & \textup{if $-r+1 \leq N \leq -1$.}
\end{cases}
\end{equation}
When $N=0$, we have $\cT_0(t) = \sum_{m=0}^{\infty} \theta^{-m-1} t^m = 1/(\theta-t) = \cS_0(t)$, with this equality occurring in $\TT$. Recalling from~\S\ref{SS:phistarext} that $\brac{\theta^{-m-1}}{-N}_{\phi^*} = 0$ for $N > -r(m+1)$, i.e., for $m > -N/r - 1$, it follows that $\cT_N(t) = 0 = \cS_N(t)$ for $-r+1 \leq N \leq -1$.

Let $N \in \ZZ$. We claim that $\cS_N$ satisfies two recursions,
\begin{equation} \label{E:SNrecs}
  \cS_N=\frac{1}{t-\theta}\sum_{i=1}^{r}b_i \cS_{N-i}^{(i)} = \frac{1}{t-\theta^{(N)}} \sum_{i=1}^r b_i^{(N-i)} \cS_{N-i}.
\end{equation}
Since $\phi_\theta(\bsomega_\phi)=t\bsomega_\phi$, we have $\bsomega_\phi= (t-\theta)^{-1} \sum_{i=1}^{r}b_i \bsomega_\phi^{(i)}$,
from which
\[
\cS_N = \bszeta_{\phi}^{(N)} \bsomega_\phi=\frac{1}{t-\theta}\sum_{i=1}^{r}b_i \bszeta_{\phi}^{(N)} \bsomega_\phi^{(i)} = \frac{1}{t-\theta} \sum_{i=1}^r b_i \cS_{N-i}^{(i)},
\]
the first recursion in~\eqref{E:SNrecs}. Likewise, $\phi_{\theta}^*(\bszeta_{\phi}) = t\bszeta_\phi$ implies $\bszeta_\phi = (t-\theta)^{-1} \sum_{i=1}^{r}b_i^{(-i)} \bszeta_\phi^{(-i)}$. Since $\cS_N = \bszeta_{\phi}^{(N)} \bsomega$, we obtain the second recursion in~\eqref{E:SNrecs}.

We observe that either recursion in~\eqref{E:SNrecs} uniquely determines $\{\cS_N \}_{N \in \ZZ}$ from the base information in~\eqref{E:SNbase}. Indeed for $N > 0$, $\cS_N$ is determined uniquely from $\cS_{N-1}, \dots, \cS_{N-r}$ using either recursion. By rearranging and reindexing the second sum in~\eqref{E:SNrecs},
\[
\cS_N = \frac{t-\theta^{(N+r)}}{b_r^{(N)}} \cS_{N+r} - \frac{1}{b_r^{(N)}} \sum_{i=1}^{r-1} b_i^{(N+r-i)} \cS_{N+r-i}, \quad N \in \ZZ.
\]
One can make a similar rearrangement using the first sum in~\eqref{E:SNrecs}, which we omit. Nevertheless for $N \leq -r$, we have that either recursion in~\eqref{E:SNrecs} defines $\cS_N$ in terms of $\cS_{N+1}, \dots, \cS_{N+r}$, starting with the base cases in~\eqref{E:SNbase}.

Let $N > 0$. Noting that $\brac{1}{N}_\phi=0$ for $N>0$, it follows from~\eqref{E:phibracrec} and the base cases for $-r+1 \leq N \leq -1$ that $\cT_N$ satisfies the two recursions,
\begin{equation} \label{E:posTNrecursion}
  \cT_N=\frac{1}{t-\theta}\sum_{i=1}^{r}b_i \cT_{N-i}^{(i)} = \frac{1}{t-\theta^{(N)}} \sum_{i=1}^r b_i^{(N-i)} \cT_{N-i}.
\end{equation}
Now let $N<0$. Noting that $\brac{1}{N}_{\phi^*}=\brac{1}{-N}_{\phi}=0$ for $N<0$, it follows from substituting $n=-N$ in ~\eqref{E:phistarbracrec} and the base cases for $-r+1 \leq N \leq -1$ that $\cT_N$, for $N < 0$, also satisfies~\eqref{E:posTNrecursion}. Therefore, $\{ \cS_N \}_{N\in \ZZ}$ and $\{ \cT_N \}_{N \in \ZZ}$ satisfy the same recursions with the same initial conditions, so as noted above we must have $\cS_N = \cT_N$ for all $N \in \ZZ$.
\end{proof}

\begin{remark} \label{R:PsiUpsilonij}
For all $M$, $N \in \ZZ$ and $1 \leq i, j \leq r$, Proposition~\ref{P:PsiNUpsilon11} permits calculation of $[\Psi^{(N)} \Upsilon^{(M)}]_{ij}$. Indeed, using~\eqref{E:Upsilondef} and~\eqref{E:zetadef}, we have
\[
\bigl[ \Psi^{(N)} \Upsilon^{(M)} \bigr]_{ij}
= \bigl[ \Psi^{(N-i+1)} \Upsilon^{(M+j-1)} \bigr]_{11}
= \bigl[ \Psi^{(N-M-i-j+2)} \Upsilon \bigr]_{11}^{(M+j-1)}.
\]
\end{remark}

\begin{remark}
When considering the $N\geq 0$ case of Proposition~\ref{P:PsiNUpsilon11}, we note that Chen used a similar recursion to~\eqref{E:SNrecs} in~\cite{ChenYT24}*{Lem.~3.1} to prove that $\cS_N= \cB_N/(\theta-t)$, where $\{ \cB_N \}$ is a sequence of rational functions defined in~\cite{EP14} using the combinatorics of shadowed partitions. A consequence of Proposition~\ref{P:PsiNUpsilon11} is that
\[
\dfrac{\cB_N}{\theta-t} = \sum_{m=0}^{\infty} \brac{\theta^{-m-1}}{N}_{\phi} t^m, \quad N \geq 0,
\]
which gives a new description of $\cB_N$ as a power series. Indeed the recursion in~\eqref{E:SNrecs} and the corresponding one in~\cite{ChenYT24}*{Lem.~3.1} are the characterizing recursions for $\cB_N$ in~\cite{EP14}*{Lem.~6.12}, and one can check using~\eqref{E:bracxdef} that $(\theta-t)\cS_N|_{t=\theta} = \beta_N$ as in~\cite{EP14}*{Eq.~(6.5)}.
\end{remark}

\section{Universal Anderson eigenvectors} \label{S:Aeigenvectors}

We now consider results of Ferraro~\cite{Ferraro25}, which provide remarkable extensions of work of Poonen~\cite{Poonen96}, particularly focusing on universal Anderson eigenvectors and dual Anderson eigenvectors. We compare them to the rigid analytic trivializations reviewed in~\S\ref{SS:RATs}, leading to new identities for $\zeta_1, \dots, \zeta_r \in \TT$ from~\eqref{E:zetadef} (see Corollary~\ref{C:samezetas}).

\begin{subsec}{Modules of eigenvectors}
Because we are considering the case $A=\FF_q[\theta]$, as opposed to general base rings~$A$, several aspects of~\cite{Ferraro25} simplify (see also \cite{Ferraro25}*{\S 7.5}). For an $A$-module~$M$, we let $\hM \assign \Hom_{\FF_q}^{\cont}(M,\FF_q)$ denote its Pontryagin dual, and for a pair of $A$-modules $M$, $N$, with $N$ locally compact, their topological tensor product is $M \hotimes N \assign \Hom_{\FF_q}^{\cont}(\hN,M)$. See \citelist{\cite{Ferraro25}*{\S 2} \cite{Poonen96}*{\S 7}}. Under these definitions it follows that
\begin{equation} \label{E:completeTate}
\C \hotimes A \cong \TT, \quad (\theta \otimes 1 \mapsto \theta,\ 1 \otimes \theta \mapsto t),
\end{equation}
as topological $A \otimes_{\FF_q} A$-modules. Now for a discrete $A$-module $M$, we set
\begin{align*}
\Sf_{\phi}(M) &\assign \Hom_{A}^{\cont}(\hM,\phi(\C))
\cong \{ \alpha \in \phi(\C) \hotimes M \mid (\phi_{\theta} \otimes 1)\alpha = (1 \otimes \theta)\alpha \}, \\
\Sf_{\phi^*}(M) &\assign \Hom_{A}^{\cont}(\hM,\phi^*(\C))
\cong \{ \alpha \in \phi^*(\C) \hotimes M \mid (\phi_{\theta}^* \otimes 1)\alpha = (1 \otimes \theta)\alpha \}.
\end{align*}
to be the $A$-modules of \emph{Anderson eigenvectors} and \emph{dual Anderson eigenvectors of $M$} (see \cite{Ferraro25}*{\S 3, \S 6}). And indeed we see that $\Sf_{\phi}(A)$ and $\Sf_{\phi^*}(A)$ are precisely the modules of special functions and dual special functions of $\phi$ contained in $\TT$, as defined in~\S\ref{SS:special} (see also~\cite{GazdaMaurischat21}). Considering $\Sf_{\phi}$ and $\Sf_{\phi^*}$ as functors from $A$-modules to $A$-modules, one of the main results of~\cite{Ferraro25} is the following. We recall the convention that a sum $\sum'_{a \in G} s_a$ over an additive group $G$ denotes that index $a=0$ is omitted.
\end{subsec}

\begin{theorem}[{Ferraro \cite{Ferraro25}*{Thm.~3.10 \& Thm.~6.4}}] \label{T:FerraroUniv}
The functor $\Sf_{\phi}$ is isomorphic to $\Hom_{A}(\Hom_A(\Lambda_{\phi},A\,d\theta),-)$, and $\Sf_{\phi^*}$ is isomorphic to $\Hom_A(\Lambda_\phi,-)$.
\begin{alphenumerate}
\item Let $\eta_1, \dots, \eta_r \in \Hom_A(\Lambda_{\phi},A)$ denote the dual $A$-basis of  $\pi_1, \dots, \pi_r$. Then the universal element in $\Sf_{\phi}(\Hom_{A}(\Lambda_{\phi},A\,d\theta))$ is
\[
\omega_{\phi} \assign \sum_{i=1}^r \sum_{m=0}^{\infty} \exp_{\phi} \biggl( \frac{\pi_i}{\theta^{m+1}} \biggr) \otimes \theta^m \eta_i\,d\theta \quad \in \phi(\C) \hotimes \Hom_A(\Lambda_{\phi},A\,d\theta).
\]
\item The universal element in $\Sf_{\phi^*}(\Lambda_{\phi})$ is
\[
\zeta_{\phi} \assign - \sumprime{\lambda \in \Lambda_{\phi}} \lambda^{-1} \otimes \lambda \quad \in \phi^*(\C) \hotimes \Lambda_{\phi}.
\]
\end{alphenumerate}
\end{theorem}

\begin{subsec}{Anderson eigenvectors}
Theorem~\ref{T:FerraroUniv}(a) implies that for a discrete $A$-module~$M$,
\[
\Sf_{\phi}(M) = \biggl\{ \sum_{i=0}^r \sum_{m=0}^{\infty} \exp_{\phi} \biggl( \frac{\pi}{\theta^{m+1}} \biggr) \otimes \theta^m \mu_i \biggm| \mu_1, \dots, \mu_r \in M \biggr\}.
\]
In particular by~\eqref{E:completeTate}, $\omega_1, \dots, \omega_r$ form an $\FF_q[t]$-basis of $\Sf_{\phi}(A)$, as noted in~\S\ref{SS:AGF}.
\end{subsec}

\begin{subsec}{Dual Anderson eigenvectors}\label{SS:DualAndersonEigenvectors}
Theorem~\ref{T:FerraroUniv}(b) implies for a discrete $A$-module~$M$,
\[
\Sf_{\phi^*}(M) = \biggl\{ -\sumprime{\lambda \in \Lambda_{\phi}} \lambda^{-1} \otimes \beta(\lambda) \biggm| \beta \in \Hom_A(\Lambda_{\phi},M) \biggr\}.
\]
Therefore, with $\eta_i : \Lambda_{\phi} \to A$ as above, we define for $i=1, \dots, r$,
\[
\xi_i \assign -\sumprime{\lambda \in \Lambda_{\phi}} \lambda^{-1} \otimes \eta_i(\lambda)  \quad \in \TT,
\]
and then $\xi_1, \dots, \xi_r$ form an $\FF_q[t]$-basis of $\Sf_{\phi^*}(A) \subseteq \TT$. If we consider the dual basis relation of $\eta_i(\pi_j) = \delta_{ij}$ for $1\leq i,j \leq r$, then we can rewrite the $\xi_i$ as 
\begin{equation} \label{E:xisum}
\xi_i = -\sumprime{(a_1, \dots, a_r) \in A^r} \frac{a_i(t)}{a_1 \pi_1 + \cdots + a_r\pi_r}, \quad 1 \leq i \leq r.
\end{equation}
\end{subsec}

\begin{theorem}[Ferraro~\cite{Ferraro25}*{\S 7}] \label{T:Ferraropairing}
The canonical $A$-bilinear pairing, $(\lambda,\alpha) \mapsto \alpha(\lambda) : \Lambda_{\phi} \times \Hom_A( \Lambda_{\phi},A \,d\theta) \to A\,d\theta$, induces a $\C\otimes_{\FF_q} A$-bilinear pairing
\begin{equation} \label{E:Ferraropairing}
(\alpha,\beta) \mapsto \alpha \cdot \beta : \C \hotimes \Lambda_{\phi} \times \C \hotimes \Hom_A(\Lambda_{\phi},A\,d\theta) \to (\C \hotimes A)\,dt,
\end{equation}
under which
\[
\zeta_\phi \cdot \omega_{\phi}^{(k)} = \begin{cases}
\dfrac{d(1\otimes \theta)}{\theta\otimes 1- 1 \otimes \theta} & \text{if $k=0$,} \\[10pt]
0 & \text{if $1 \leq k \leq r-1$,}
\end{cases}
\]
where $\omega_{\phi}^{(k)}$ is the $k$-th twist under $\tau \otimes \id$.
\end{theorem}

As in the proof of Proposition~\ref{P:PsiNUpsilon11}, we recall $\bszeta_{\phi} =(\zeta_1, \dots, \zeta_r) \in \Mat_{1 \times r}(\TT)$ and $\bsomega_{\phi} =(\omega_1, \dots, \omega_r)^{\tr} \in \Mat_{r \times 1}(\TT)$, and we further let $\bsxi_{\phi} \assign (\xi_1, \dots, \xi_r) \in \Mat_{1\times r}(\TT)$. Now under our choice of $A$-basis $\pi_1, \dots, \pi_r$ for $\Lambda_{\phi}$, the $\C[t]$-bilinear pairing induced by \eqref{E:Ferraropairing},
$\Mat_{1\times r}(\TT) \times \Mat_{r\times 1}(\TT) \to \TT$,
is simply the dot product (canonically identifying $\TT \cong \TT\,dt$). In particular, Theorem~\ref{T:Ferraropairing} implies
\begin{equation} \label{E:bszetabsomega}
\bsxi_{\phi} \cdot \bsomega_{\phi}^{(k)} = \begin{cases}
\dfrac{1}{\theta - t} & \text{if $k=0$,} \\[10pt]
0 & \text{if $1 \leq k \leq r-1$.}
\end{cases}
\end{equation}
The following corollary shows that our dual special functions $\zeta_1, \dots \zeta_r$ from \eqref{E:zetadef} are in fact the same as $\xi_1, \dots, \xi_r$, i.e., $\bszeta_{\phi} = \bsxi_{\phi}$.

\begin{corollary} \label{C:samezetas}
Under our choice of $A$-basis $\pi_1, \dots, \pi_r$ for $\Lambda_{\phi}$, we have
\[
\zeta_i = \xi_i = - \sumprime{(a_1, \dots, a_r) \in A^r} \frac{a_i(t)}{a_1 \pi_1 + \cdots + a_r\pi_r}, \quad \forall\, i = 1, \dots, r.
\]
\end{corollary}

\begin{proof}
We note that~\eqref{E:bszetabsomega} can be translated into the matrix form, $(\xi_1, \dots, \xi_r) \cdot \Upsilon = ( 1/(\theta-t), 0, \dots, 0)$. Since $\Upsilon \in \GL_r(\TT)$, we have
\[
(\xi_1, \dots, \xi_r)
= \biggl( \frac{1}{\theta-t}, 0, \dots, 0 \biggr) \Upsilon^{-1}
= \textup{first row of $(\Psi \Upsilon)\Upsilon^{-1}$}
= (\zeta_1, \dots, \zeta_r),
\]
where the second equality follows from Lemma~\ref{L:PsiUpsilon}.
We are done by~\eqref{E:xisum}.
\end{proof}

\section{The maps \texorpdfstring{$\be_{\lambda}$}{e-lambda} and \texorpdfstring{$\be_{\eta}^*$}{e-eta-star}} \label{S:emaps}

\begin{subsec}{Residue maps}
For a field $F$ and a double-tailed Laurent series
\[
x = \sum_{j \in \ZZ} c_j z^{-j} \in \power{F}{z,z^{-1}},
\]
we have the $F$-linear map $\rR_z : \power{F}{z,z^{-1}} \to F$ given by
\[
\rR_z(x) \assign c_1.
\]
We note that if $t$ is the canonical uniformizing parameter on $\PP^1/F$ vanishing at $0$ and if, for $f \in F(t)$, the differential $f\,dt$ has expansion at $\infty$,
$f\,dt = \sum_{j \geq j_0} c_j t^{-j} \, dt$, then
\begin{equation} \label{E:residue}
c_1 = \rR_t \biggl( \sum_{j \geq j_0} c_j t^{-j} \biggr) = -\Res_{\infty} ( f\,dt ),
\end{equation}
where the latter is the standard residue of a differential form on an algebraic curve.
\end{subsec}

\begin{subsec}{The maps $\be_{\lambda}$} \label{SS:elambda}
For $\lambda \in \Lambda_{\phi}$, define the $A$-module homomorphism
\begin{equation} \label{E:elambdadef}
\be_{\lambda} \assign \bigl( x \mapsto \exp_{\phi}(\lambda x) \bigr) : k_{\infty} \to \phi(\C).
\end{equation}
As such, $\ker \be_\lambda \supseteq A$ for each $\lambda \in \Lambda_{\phi} \setminus \{ 0 \}$. Moreover, the induced map
\begin{equation} \label{E:elambdaHom}
\lambda \mapsto \be_{\lambda} : \Lambda_{\phi} \to \Hom_A(k_{\infty}/A,\phi(\C))
\end{equation}
is injective and $A$-linear. For our choice of $A$-basis $\pi_1, \dots, \pi_r$ of $\Lambda_{\phi}$, we set
\[
\be_i(x) \assign \be_{\pi_i}(x), \quad 1 \leq i \leq r,
\]
for which $\ker \be_i = A$. In particular for each $b \in A$, $b \neq 0$, it follows from~\eqref{E:torsionexp} that
\begin{equation} \label{E:eitorsion}
\phi[b] = \bigoplus_{i=1}^r \be_i \bigl( b^{-1} A\bigr).
\end{equation}

Now for $x \in k_{\infty}$, we can form the product $x(t) \omega_\lambda(t) \in \power{\C}{t,t^{-1}}$. Indeed, for example by the argument in~\cite{NP24}*{Rem.~4.29}, for $\alpha \in \C$ with $1 < \inorm{\alpha}< \inorm{\theta}$, we have $x(t)$, $\omega_\lambda(t) \in \SSS_{\alpha}$ as in \S\ref{SS:RAfns}, and so $x(t)\omega_\lambda(t) \in \SSS_{\alpha} \subseteq \power{\C}{t, t^{-1}}$. We then have the following.
\end{subsec}

\begin{lemma} \label{L:elambdaRes}
Let $\lambda \in \Lambda_{\phi}$. For $x \in k_{\infty}$, we have $\be_\lambda(x) = \rR_t ( x(t) \omega_\lambda(t) )$.
\end{lemma}

\begin{proof}
Let $x = \sum_{j\geq j_0} c_j \theta^{-j}$. From the definition of $\omega_\lambda$ in~\S\ref{SS:AGF}, we have
\begin{equation} \label{E:omegalambdaelambda}
\omega_\lambda(t) = \sum_{m \geq 0} \be_\lambda(1/\theta^{m+1}) t^m.
\end{equation}
Thus, the coefficient of $t^{-1}$ in $x(t)\omega_\lambda(t)$ is
\[
\sum_{m=0}^{\infty} c_{m+1} \be_{\lambda} (1/\theta^{m+1}) = \be_{\lambda} \biggl( \sum_{j\geq 1} c_j \theta^{-j} \biggr) = \be_\lambda(x),
\]
the first equality by the continuity and $\FF_q$-linearity of $\be_\lambda$, the second since $\be_{\lambda}(A)=0$.
\end{proof}

\begin{subsec}{The maps $\be_\eta^*$} \label{SS:eetastar}
As in \S\ref{SS:elambda}, we would like to define $A$-module homomorphisms
\[
\be_\eta^* : k_{\infty} \to \phi^*(\C), \quad \eta \in \Hom_A(\Lambda_\phi,A),
\]
that satisfy similar properties such as \eqref{E:elambdadef}--\eqref{E:omegalambdaelambda}, and particularly Lemma~\ref{L:elambdaRes}. However, as indicated in~\S\ref{SS:torsexp}, the series $\exp_{\phi}^*$ and $\log_{\phi}^*$ are inadequate to the task.

Instead, let $\eta_1, \dots, \eta_r \in \Hom_A(\Lambda_{\phi},A)$ denote the dual basis of $\pi_1, \dots, \pi_r \in \Lambda_{\phi}$.
For $\eta \in \Hom_A(\Lambda_\phi,A)$, writing $\eta = s_1 \eta_1 + \cdots + s_r \eta_r$ for $s_1 \dots, s_r \in A$, we set
\begin{equation} \label{E:zetaetadef}
\zeta_{\eta}(t) \assign \bszeta_{\phi} \cdot (s_1(t), \dots, s_r(t))^{\tr} = s_1(t) \zeta_1(t) + \dots + s_r(t)\zeta_r(t).
\end{equation}
We note that $\zeta_i(t) =\zeta_{\eta_i}(t)$, for $1 \leq i \leq r$.
The reader will recognize that Theorem~\ref{T:FerraroUniv} of Ferraro necessitates indexing $\zeta_{\eta}(t)$ by $\Hom_A(\Lambda_\phi,A)$, rather than simply by~$\Lambda_\phi$. It is also clear from the following lemma.
\end{subsec}

\begin{lemma} \label{L:zetaeta}
For $\eta \in \Hom_A(\Lambda_\phi,A)$, the function $\zeta_{\eta}(t) \in \TT$ is independent of the choice of $A$-basis $\pi_1, \dots, \pi_r$.
\end{lemma}

\begin{proof}
Suppose $\pi_1', \dots, \pi_r'$ is an $A$-basis of $\Lambda_{\phi}$ with corresponding dual basis $\eta_1', \dots, \eta_r'$, and let $\zeta_{\eta}'(t)$ denote \eqref{E:zetaetadef} with respect to this basis. Then for some $B \in \GL_r(A)$,
\[
\bigl( \pi_1', \dots, \pi_r' \bigr)^{\tr}
= B ( \pi_1, \dots, \pi_r )^{\tr},
\quad
B^{\tr} \bigl( \eta_1', \dots, \eta_r' \bigr)^{\tr}
= ( \eta_1, \dots, \eta_r )^{\tr}.
\]
Let $\Upsilon' \in \GL_r(\TT)$ be defined as in \eqref{E:Upsilondef} with respect to $\pi_1', \dots, \pi_r'$. Taking $B(t) = B|_{\theta=t} \in \GL_r(\FF_q[t])$, it follows from \eqref{E:AGFpereq} that $\Upsilon' = B(t) \Upsilon$. The corresponding rigid analytic trivialization $\Psi'$ as in~\eqref{E:Psidef} is then
\[
\Psi' = -V^{-1} \bigl( \Upsilon'{}^{(1)} \bigr)^{-1} = -V^{-1} \bigl( \Gamma^{(1)} \bigr)^{-1} B(t)^{-1} = \Psi B(t)^{-1},
\]
from which we see that for $\bszeta_{\phi}' = (\zeta_1', \dots, \zeta_r')$, the first row of $\Psi'$, we have $\bszeta_{\phi} = \bszeta_{\phi}' \cdot B(t)$. Thus, if $\eta = (s_1, \dots, s_r) \cdot (\eta_1, \dots, \eta_r)^{\tr}$, then by definition, $\zeta_{\eta}(t) = \bszeta_{\phi} \cdot (s_1(t), \dots, s_r(t))^{\tr}$. On the other hand, $\eta = (s_1, \dots, s_r) B^{\tr} (\eta_1', \dots, \eta_r')^{\tr}$,
which implies
\[
\zeta_{\eta}'(t) = \bszeta_\phi' \cdot B(t) \cdot (s_1(t), \dots, s_r(t))^{\tr} = \bszeta_{\phi} \cdot (s_1(t), \dots, s_r(t))^{\tr} = \zeta_{\eta}(t).
\qedhere
\]
\end{proof}

Now taking inspiration from~\cite{ABP04}*{\S 5.2.5}, we define for $\eta \in \Hom_A(\Lambda_{\phi},A)$ and $x \in k_{\infty}$,
\begin{equation} \label{E:eetastardef}
\be_\eta^*(x) \assign \rR_t\bigl( x(t)\zeta_{\eta}(t) \bigr),
\end{equation}
noting as in \S\ref{SS:elambda} that $x(t)\zeta_\eta(t) \in \SSS_{\alpha} \subseteq \power{\C}{t,t^{-1}}$ for any $\alpha \in \C$ with $1 < \inorm{\alpha} < \inorm{\theta}$. Thus by definition $\be_\eta^*$ satisfies the conclusion of Lemma~\ref{L:elambdaRes}, and by Lemma~\ref{L:zetaeta} its definition does not depend on the representation of $\eta$ with respect to our $A$-basis of~$\Lambda_{\phi}$. Because $\zeta_{\eta}(t) \in \power{\C}{t}$, the definition of $\rR_t$ implies $\ker \be_{\eta}^* \supseteq A$, and as in~\eqref{E:elambdaHom}, the map
\[
\eta \mapsto \be_{\eta}^* : \Hom_A(\Lambda_{\phi},A) \to \Hom_A(k_{\infty}/A,\phi^*(\C))
\]
is injective and $A$-linear. By~\eqref{E:eetastardef}, taking $x = \theta^{-m-1}$ yields 
\begin{equation}
\zeta_{\eta}(t) = \sum_{m=0}^\infty \be_{\eta}^*(1/\theta^{m+1}) t^m,
\end{equation}
which gives the companion identity to~\eqref{E:omegalambdaelambda}. For an $A$-basis $\pi_1, \dots \pi_r \in \Lambda_{\phi}$ and corresponding dual basis $\eta_1, \dots, \eta_r \in \Hom_A(\Lambda_\phi,A)$, we set
\[
\be_i^* \assign \be_{\eta_i}^* : k_{\infty} \to \phi^*(\C), \quad 1 \leq i \leq r,
\]
and so if $\eta = s_1 \eta_1 + \cdots + s_r \eta_r$ with $s_1, \dots, s_r \in A$, then $\be_{\eta}^* = \phi_{s_1}^* \circ \be_1^* + \dots + \phi_{s_r}^* \circ \be_r^*$.

\begin{lemma} \label{L:eetaresidue}
Let $\eta \in \Hom_A(\Lambda_{\phi},A)$.
For $x \in k_{\infty}$,
\[
\be_{\eta}^*(x) = -\sumprime{\lambda \in \Lambda_{\phi}} \frac{\rR_{\theta} (x \eta(\lambda))}{\lambda}.
\]
\end{lemma}

\begin{proof}
By forgetting $A$-linearity, we have an inclusion $\Hom_A(\Lambda_{\phi},A) \hookrightarrow \Hom_{\FF_q}(\Lambda_\phi,A)$, and by picking an $\FF_q$-basis of $A$, we have $\Hom_{\FF_q}(\Lambda_\phi,A) \cong \oplus_{j \geq 0} \Hom_{\FF_q}(\Lambda_{\phi},\FF_q)$. Taking $\{ \theta^j \}$ as this basis, the induced isomorphism is given by
\[
\eta \mapsto \Bigl( \rR_{\theta} \bigl( \theta^{-1} \eta(\,\cdot\,) \bigr),\, \rR_{\theta} \bigl( \theta^{-2} \eta(\,\cdot\,) \bigr), \dots \Bigr).
\]
Since $\rR_{\theta}(x \eta(\lambda)) \in \FF_q$ and $\Lambda$ is a finitely generated lattice over $A$ in $\C$, the sum
\[
\sumprime{\lambda \in \Lambda_{\phi}} \frac{\rR_{\theta}(x \eta(\lambda))}{\lambda}, \quad \eta \in \Hom_A(\Lambda_{\phi},A),\ x \in k_\infty,
\]
converges in $\C$. Now if $\eta = s_1 \eta_1 + \cdots + s_r \eta_r$ with $s_1, \dots, s_r \in A$, then for $j \geq 0$,
\begin{align*}
-\sumprime{\lambda \in \Lambda_{\phi}} \frac{\rR_{\theta}(\theta^{-j-1} \eta(\lambda))}{\lambda}
&= - \sum_{i=1}^r \sumprime{(a_1, \dots, a_r) \in A^r} \frac{\rR_{t}(t^{-j-1} s_i(t) a_i(t))}{a_1 \pi_1 + \dots + a_r \pi_r} && \\
&= \sum_{i=1}^r \rR_t \bigl( t^{-j-1}s_i(t) \zeta_i(t) \bigr) && \textup{(Corollary~\ref{C:samezetas})} \\
&=\rR_t \bigl( t^{-j-1} \zeta_{\eta}(t) \bigr) && \textup{(Equation~\eqref{E:zetaetadef})} \\
&= \be_{\eta}^*(1/\theta^{j+1}) && \textup{(Equation~\eqref{E:eetastardef})}.
\end{align*}
The result then follows by extending $\FF_q$-linearly and continuously to $k_{\infty}$.
\end{proof}

\begin{subsec}{Determining $\ker \be_{\eta}^*$}
We recall results of Poonen~\cite{Poonen96}*{Thm.~10} and Ferraro~\cite{Ferraro25}*{Prop.~5.18}. As mentioned in \S\ref{SS:torsexp}, we know from \cite{Poonen96}*{Prop.~20} that $\phi_{\tors}^*$ is a dense subset of the compact space $\ker \exp_{\phi}^*$. By~\cite{Ferraro25}*{Lem.~5.3}, for $\beta \in \ker \exp_{\phi}^*$ there is a unique $h_{\beta} \in \power{\C}{\tau}$ that converges on all of $\C$ such that $(1-\tau)h_{\beta} = \beta \exp_{\phi}$. Poonen~\cite{Poonen96}*{Thm.~10} then proved (see also \cite{Ferraro25}*{Thm.~5.5}) that the map
\begin{equation} \label{E:Pooneniso}
\beta \mapsto h_{\beta}|_{\Lambda_{\phi}} : \ker \exp_{\phi}^* \iso \hat{\Lambda}_{\phi}
\end{equation}
is an $A$-linear homeomorphism. Furthermore, Ferraro~\cite{Ferraro25}*{Prop.~5.18} proved that
\begin{equation} \label{E:betarecover}
\beta = - \sumprime{\lambda \in \Lambda_{\phi}} \frac{h_{\beta}(\lambda)}{\lambda},
\quad \forall\, \beta \in \ker \exp^*.
\end{equation}

Fixing $\eta \in \Hom_A(\Lambda_{\phi},A)$, for $x \in k_{\infty}$ the $\FF_q$-linear map
\[
\gamma_x = \bigl( \lambda \mapsto \rR_{\theta}(x\eta(\lambda)) : \Lambda \to \FF_q \bigr)
\]
is an element of $\hat{\Lambda}_{\phi}$. Thus by~\eqref{E:Pooneniso} there is a unique $\beta_x \in \ker \exp_{\phi}^*$ so that $\gamma_x = h_{\beta_x}|_{\Lambda_{\phi}}$. Then by Lemma~\ref{L:eetaresidue} and~\eqref{E:betarecover}, $\be_{\eta}^*(x) = \beta_x$. Therefore,
\[
\ker \be_{\eta}^* = \{ x \in k_{\infty} \mid \rR_{\theta}(x\eta(\lambda)) = 0, \ \forall\,\lambda \in \Lambda_{\phi} \}.
\]
Now let us assume that $\eta \neq 0$. In this case $\im \eta = bA \subseteq A$ for some $b \in A$, $b \neq 0$. Then
\[
x \in \ker \be_{\eta}^* \quad \Leftrightarrow \quad \rR_{\theta}(xba) = 0, \quad \forall\, a \in A.
\]
We claim this implies that
\begin{equation} \label{E:kereetastar}
\ker \be_{\eta}^* = b^{-1} A.
\end{equation}
Indeed, $b^{-1}A \subseteq \ker \be_{\eta}^*$ is clear. On the other hand, if $bx \notin A$, for some $x \in k_{\infty}$, then $bx = a + \sum_{j\geq 1} c_j \theta^{-j}$ for some $a \in A$ and with $c_{j_0} \neq 0$ for some $j_0 \geq 1$. However, then $\rR_{\theta}(bx\theta^{j_0-1}) = c_{j_0} \neq 0$, whence $bx \notin \ker \be_{\eta}^*$.
\end{subsec}

\begin{lemma} \label{L:e*propertiesfortorsion}
For an $A$-basis $\pi_1, \dots, \pi_r \in \Lambda_{\phi}$ and corresponding dual basis $\eta_1, \dots, \eta_r \in \Hom_A(\Lambda_{\phi},A)$, the following hold.
\begin{alphenumerate}
\item For each $i$, $\ker \be_i^* = A$.
\item For $b \in A$, $b \neq 0$,
\[
\phi^*[b] = \bigoplus_{i=1}^r \be_i^* \bigl( b^{-1} A \bigr).
\]
\end{alphenumerate}
\end{lemma}

\begin{proof}
Part (a) follows from~\eqref{E:kereetastar}, since $\eta_i(\Lambda) = A$ for each $i$. For~(b), the case that $\deg b = 0$ is clear, so we then assume $\deg b \geq 1$. We first note that $\be_i^*(b^{-1}A) \cong A/bA$ as $A$-modules by part~(a). Second, we let $N=\sum_{i=1}^r \be_i^*(b^{-1}A)$, and we show this sum is direct. Suppose that $\be_1^*(b^{-1}s_1) + \dots + \be_r^*(b^{-1}s_r) = 0$ for some $s_1, \dots, s_r \in A$ with $\gcd(s_1, \dots, s_r)=s$. Therefore, if we let $\eta = s_1\eta_1 + \dots + s_r\eta_r$, then
\[
0 = \phi_{s_1}^* \be_1^*(1/b) + \dots + \phi_{s_r}^* \be_r^*(1/b)  = \be_{\eta}^*(1/b).
\]
This implies that $b^{-1}A \subseteq \ker \be_{\eta}^*$, whereas by~\eqref{E:kereetastar}, $\ker \be_{\eta}^* = s^{-1}A$. Hence $b \mid s_i$ for each~$i$, and so $\be_i^*(b^{-1}s_i)=0$ for each $i$, as desired.
\end{proof}

\section{Gauss-Thakur sums and dual Gauss-Thakur sums} \label{S:GTsums}

We review definitions of Gauss-Thakur sums of Gazda and Maurischat~\cite{GazdaMaurischat21} for Drinfeld modules, which extend Thakur's original constructions~\citelist{\cite{Thakur88} \cite{Thakur93b} \cite{Thakur}*{Ch.~4}}. At the same time we define their dual counterparts, built from torsion on the adjoint Drinfeld module.

\begin{subsec}{Evaluative characters} \label{SS:evalchars}
We fix a monic irreducible polynomial $\wp \in A$ of degree~$d$, and we let $\FF_{\wp} \assign A/\wp$. A multiplicative character $\chi : \FF_{\wp}^{\times} \to \oFF_q^{\times}$ extends in the usual way to a Dirichlet character modulo $\wp$, $\chi : A \to \oFF_q$,
where $\chi(a)=0$ if $\wp \mid a$. By selecting $\varepsilon \in \FF_{q^d} \subseteq \oFF_q$ that is a root of $\wp$, we define the \emph{evaluative character},
\[
\chi_{\varepsilon} \assign \bigl( a \mapsto a(\varepsilon) \bigr) : A \to \FF_{q^d},
\]
and we find that $\widehat{\FF}_{\wp}^{\times} = \langle \chi_{\varepsilon} \rangle$. If $\varepsilon$ has been fixed, then to ease notation, for $k \in \ZZ$ we will sometimes write $\chi_k \assign \chi_{\varepsilon^{q^k}}$.
\end{subsec}

\begin{subsec}{Gauss-Thakur sums} \label{SS:Gauss-ThakurSums}
Let $\phi$ be a Drinfeld module of rank~$r$ as in~\eqref{E:rhodef}. For a Dirichlet character $\chi : A \to \FF_{q^d}$ modulo $\wp$, $A$-module homomorphisms $\psi : A \to \phi[\wp]$ and $\psi^*: A \to \phi^*[\wp]$, and $a \in A$, the Gauss-Thakur sum of~\cite{GazdaMaurischat21} and its dual sum are
\begin{align*}
G(\chi,\psi;a) &\assign -\sum_{z \in \FF_{\wp}^{\times}} \chi(z)^{-1} \otimes \psi(az) \quad \in \FF_{q^d} \otimes_{\FF_q} \phi(\C), \\
G^*(\chi,\psi^*;a) &\assign -\sum_{z \in \FF_{\wp}^{\times}} \chi(z)^{-1} \otimes \psi^*(az) \quad \in \FF_{q^d} \otimes_{\FF_q} \phi^*(\C).
\end{align*}
We set $G(\chi,\psi) \assign G(\chi,\psi;1)$ and $G^*(\chi,\psi^*) \assign G(\chi,\psi^*;1)$, and we have for $a$, $b \in A$,
\begin{align} \label{E:Gcharab}
G(\chi,\psi;ab) &= G(\chi,\phi_a \circ \psi;b) = (\chi(a) \otimes 1) G(\chi,\psi;b), \\
G^*(\chi,\psi^*;ab) &= G^*(\chi,\phi_a^*\circ \psi^*;b) = (\chi(a) \otimes 1) G^*(\chi,\psi^*;b). \notag
\end{align}
Composing with the product map $\FF_{q^d} \otimes_{\FF_q} \C \to \C$, we also have the tensorless sums (as named in~\cite{GazdaMaurischat21}*{Rem.~4.3}),
\[
g(\chi,\psi;a) \assign -\sum_{z \in \FF_{\wp}^{\times}} \chi(z)^{-1} \psi(az), \quad
g^*(\chi,\psi^*;a) \assign -\sum_{z \in \FF_{\wp}^{\times}} \chi(z)^{-1} \psi^*(az),
\]
which take values in $\oK$. We also set $g(\chi,\psi) \assign g(\chi,\psi;1)$ and $g^*(\chi,\psi^*) = g^*(\chi,\psi;1)$, which for $a$, $b \in A$ satisfy,
\begin{align} \label{E:gcharab}
g(\chi,\psi;ab) &= g(\chi,\phi_a\circ \psi;b) = \chi(a) g(\chi,\psi;b), \\
g^*(\chi,\psi^*;ab) &= g^*(\chi,\phi_a^* \circ \psi^*;b) =  \chi(a) g^*(\chi,\psi^*;b). \notag
\end{align}

As in~\cite{GazdaMaurischat21}*{\S 4}, we further have the $\FF_{q^d}$-linear vector spaces of special elements,
\begin{align} \label{E:GTsumspaces}
\cG(\phi,\chi) &\assign \{ G \in \FF_{q^d} \otimes_{\FF_q} \phi(\C) \mid (1 \otimes a) G = (\chi(a) \otimes 1) G,\ \forall\, a \in A \}, \\
\cG^*(\phi^*,\chi) &\assign \{ G^* \in \FF_{q^d} \otimes_{\FF_q} \phi^*(\C) \mid (1 \otimes a) G^* = (\chi(a) \otimes 1) G^*,\ \forall\, a \in A \}. \notag
\end{align}
We see in Lemma~\ref{L:GTspaces} that $\cG(\phi,\chi)$ and $\cG^*(\phi^*,\chi)$ are really spaces of Gauss-Thakur sums and dual Gauss-Thakur sums.

As in \cite{GazdaMaurischat21}*{Ex.~4.4}, all additive characters $\psi : A \to \phi[\wp]$, $\psi^* : A \to \phi^*[\wp]$, have the form
\[
\psi_{\gamma} \assign \bigl( a \mapsto \phi_a(\gamma)\bigr) : A \to \phi[\wp], \ \gamma \in \phi[\wp]; \quad
\psi^*_{\delta} \assign \bigl( a \mapsto \phi^*_a(\delta)\bigr) : A \to \phi^*[\wp], \ \delta \in \phi^*[\wp].
\]
Also for a multiplicative character $\chi \in \widehat{\FF}_{\wp}^{\times}$, define $\FF_{q^d}$-linear maps, $\eta_\chi : \FF_{q^d} \otimes_{\FF_q} \phi[\wp] \to \FF_{q^d} \otimes_{\FF_q} \phi[\wp]$ and $\eta_\chi^* : \FF_{q^d} \otimes_{\FF_q} \phi^*[\wp] \to \FF_{q^d} \otimes_{\FF_q} \phi^*[\wp]$ by
\[
\eta_\chi(\alpha) \assign -\sum_{z \in \FF_{\wp}^\times} \bigl( \chi(z)^{-1} \otimes z \bigr) \alpha, \quad
\eta_\chi^*(\beta) \assign -\sum_{z \in \FF_{\wp}^\times} \bigl( \chi(z)^{-1} \otimes z \bigr) \beta.
\]
We then have the following lemma of Gazda and Maurischat, with added dual Gauss-Thakur sum variant. The proof for dual Gauss-Thakur sums is exactly the same with minor changes, so we omit it.
\end{subsec}

\begin{lemma}[{Gazda-Maurischat~\cite{GazdaMaurischat21}*{Lem.~4.5}}] \label{L:GTspaces}
Fix a multiplicative character $\chi \in \widehat{\FF}_{\wp}^{\times}$. With notation as above the following hold.
\begin{alphenumerate}
\item The map $\eta_\chi$ is an $\FF_{q^d}$-linear projection onto $\cG(\phi,\chi)$, and $\eta_{\chi}^*$ is an $\FF_{q^d}$-linear projection onto $\cG^*(\phi^*,\chi)$.
\item For $\gamma \in \phi[\wp]$, $\delta \in \phi^*[\wp]$, we have
$\eta_{\chi} ( 1 \otimes \gamma) = G(\chi, \psi_{\gamma} )$ or
$\eta_{\chi}^* (1 \otimes \delta) = G^*(\chi,\psi_{\delta}^*)$.
\item Every Gauss-Thakur sum $G(\chi,\psi)$ is in $\cG(\phi,\chi)$, and $\cG(\phi,\chi)$ is generated as an $\FF_{q^d}$-vector space by the Gauss-Thakur sums.
\item Likewise, every dual Gauss-Thakur sum $G^*(\chi,\psi^*)$ is in $\cG^*(\phi^*,\chi)$, and $\cG^*(\phi^*,\chi)$ is generated as an $\FF_{q^d}$-vector space by the dual Gauss-Thakur sums.
\end{alphenumerate}
\end{lemma}

We also have the following vanishing lemma for non-evaluative characters, whose proof for dual Gauss-Thakur sums is identical to the one in~\cite{GazdaMaurischat21}.

\begin{proposition}[{Gazda-Maurischat~\cite{GazdaMaurischat21}*{Prop.~4.8(1)}}] \label{P:vanishing}
Fix a multiplicative character $\chi \in \widehat{\FF}_{\wp}^{\times}$. If $\chi$ is not evaluative, then $\cG(\phi,\chi) = \{0\}$ and $\cG^*(\phi^*,\chi) = \{0\}$; in particular, for additive characters $\psi:A \to \phi[\wp]$ or $\psi^* : A \to \phi^*[\wp]$, we have
\[
G(\chi,\psi) = G^*(\chi,\psi^*) = 0, \quad
g(\chi,\psi) = g^*(\chi,\psi^*) = 0.
\]
\end{proposition}

\begin{subsec}{Gauss-Thakur sums with respect to a basis} \label{Ss:Gauss-ThakurSumBasis}
As in \S\S\ref{S:RAT}--\ref{S:emaps}, we choose an $A$-basis $\pi_1, \dots, \pi_r$ of $\Lambda_{\phi}$. Then the maps $\psi_i : A \to \phi[\wp]$, $\psi_i^* : A \to \phi^*[\wp]$ defined by
\[
\psi_i(z) \assign \be_i ( z/\wp ), \quad \psi_i^*(z) \assign \be_i^*(z/\wp), \quad 1 \leq i \leq r,
\]
are $A$-module homomorphisms each with kernel $\wp A$ (for the latter by Lemma~\ref{L:e*propertiesfortorsion}(a)). If $\lambda = a_1 \pi_1 + \dots +a_r \pi_r$ or $\eta = a_1 \eta_1 + \dots + a_r \eta_r$, for $a_1, \dots, a_r \in A$, we similarly have $\psi_{\lambda} : A \to \phi[\wp]$ and $\psi_{\eta}^*: A \to \phi^*[\wp]$ defined by
\[
\psi_{\lambda}(z) \assign \be_{\lambda}(z/\wp) = \sum_{i=1}^r \psi_i(a_i z), \quad
\psi_{\eta}^*(z) \assign \be_{\eta}^*(z/\wp) = \sum_{i=1}^r \psi_i^*(a_i z).
\]
We note that $\psi_{\lambda}$ depends on the congruence class of $\lambda$ modulo $\wp \Lambda_{\phi}$, and likewise $\psi_{\eta}^*$ depends on the class of $\eta$ modulo $\wp \Hom_A(\Lambda_{\phi},A)$.
Moreover, by \eqref{E:eitorsion} and Lemma~\ref{L:e*propertiesfortorsion},
\[
\Hom_A(A,\phi[\wp]) = \Span_A(\psi_1, \dots, \psi_r),\quad
\Hom_A(A,\phi^*[\wp]) = \Span_A(\psi_1^*,\dots, \psi_r^*),
\]
and so every $\psi \in \Hom_A(A,\phi[\wp])$ (respectively, $\psi^* \in \Hom_A(A,\phi^*[\wp])$) is of the form~$\psi_{\lambda}$ for some $\lambda \in \Lambda_{\phi}$ (respectively, $\psi_{\eta}^*$ for some $\eta \in \Hom_{A}(\Lambda_{\phi},A)$.

Thus, as a matter of notation, we define the Gauss-Thakur sums for $\lambda \in \Lambda_{\phi}$ and $a \in A$,
\begin{align}
G_{\lambda}(\chi; a) &\assign G(\chi,\psi_{\lambda}; a) = -\sum_{z \in \FF_{\wp}^{\times}} \chi(z)^{-1} \otimes \be_{\lambda}(az/\wp), &G_\lambda(\chi) &\assign G(\chi, \psi_\lambda;1), \\
g_{\lambda}(\chi; a) &\assign g(\chi, \psi_{\lambda}; a) = -\sum_{z \in \FF_{\wp}^{\times}} \chi(z)^{-1} \be_{\lambda}(az/\wp), & g_\lambda(\chi) &\assign g(\chi, \psi_\lambda;1). \notag
\end{align}
We will let $G_i(\chi;a) \assign G_{\pi_i}(\chi;a)$ and $g_i(\chi;a) \assign g_{\pi_i}(\chi;a)$ for $1 \leq i \leq r$. Necessarily we find for $\lambda = a_1\pi_1 + \dots + a_r \pi_r \in \Lambda_{\phi}$,
$G_{\lambda}(\chi) = \sum_{i=1}^r G_i(\chi; a_i)$ and 
$g_{\lambda}(\chi) = \sum_{i=1}^r g_i(\chi;a_i)$.
By Lemma~\ref{L:GTspaces}(c), $\cG(\phi,\chi)$ is generated over $\FF_{q^d}$ by $\{ G_{\lambda}(\chi) : \lambda \in \Lambda_{\phi} \}$.

We also define dual Gauss-Thakur sums for $\eta \in \Hom_A(\Lambda_{\phi},A)$ and $a \in A$,
\begin{align}
G_{\eta}^*(\chi; a) &\assign G^*(\chi,\psi_{\eta}^*; a) = -\sum_{z \in \FF_{\wp}^{\times}} \chi(z)^{-1} \otimes \be_{\eta}^*(az/\wp), &
G_\eta^*(\chi) &\assign G^*(\chi, \psi_\eta^*;1), \\
g_{\eta}^*(\chi; a) &\assign g^*(\chi, \psi_{\eta}^*; a) = -\sum_{z \in \FF_{\wp}^{\times}} \chi(z)^{-1} \be_{\eta}^*(az/\wp), &
g_\eta^*(\chi) &\assign g^*(\chi, \psi_\eta^*). \notag
\end{align}
As above, we let $G_i^*(\chi;a) \assign G_{\eta_i}^*(\chi;a)$ and $g_i^*(\chi;a) \assign g_{\eta_i}^*(\chi;a)$ for $1 \leq i \leq r$, and 
for $\eta = s_1 \eta_1 + \dots + s_r \eta_r \in \Hom(\Lambda_{\phi},A)$,
$G_{\eta}^*(\chi) = \sum_{i=1}^r G_i^*(\chi;s_i)$ and
$g_{\eta}^*(\chi) = \sum_{i=1}^r g_i^*(\chi; s_i)$.
By Lemma~\ref{L:GTspaces}(d), $\cG^*(\phi^*,\chi)$ is generated over $\FF_{q^d}$ by $\{ G_{\eta}^*(\chi) : \eta \in \Hom_A(\Lambda_{\phi},A) \}$.
\end{subsec}

\begin{subsec}{(Dual) Gauss-Thakur sums and evaluations of (dual) special functions} \label{SS:theorem5.2*}
We now investigate analogues of one of Gazda and Maurischat's main results \cite{GazdaMaurischat21}*{Thm.~5.2} in the context of dual Gauss-Thakur sums. Such results go back to Angl\`es and Pellarin~\cite{AnglesPellarin15}*{Thm.~2.9}, who proved identities between Thakur's Gauss sums for the Carlitz module and evaluations of the Anderson-Thakur special function (see \S\ref{S:CarlitzApp}).

For $\varepsilon \in \FF_{q^d}$ a fixed root of $\wp(t)=0$ we have the evaluative character $\chi_{\varepsilon}: A \to \FF_{q^d}$ as in \S\ref{SS:evalchars}. We also define $\bschi_{\varepsilon} : \C \hotimes A \to \FF_{q^d} \otimes_{\FF_q} \C$ by the composition
\[
\bschi_{\varepsilon} : \C \hotimes A \xrightarrow{1 \otimes \chi_{\varepsilon}} \C \otimes_{\FF_q} \FF_{q^d} \iso \FF_{q^d} \otimes_{\FF_q} \C,
\]
where the latter map is the canonical isomorphism that swaps the tensor product.
If we identify $\TT$ with $\C \hotimes A$ as in~\eqref{E:completeTate}, then multiplication composed with $\bschi_{\varepsilon}$,
\begin{equation} \label{E:chieval}
f \mapsto f(\varepsilon) : \TT \xrightarrow{\bschi_{\varepsilon}} \FF_{q^d} \otimes_{\FF_q} \C \to \C,
\end{equation}
is simply evaluation at $t=\varepsilon$. For $a \in A$, let $a'(\theta)$ be its derivative with respect to~$\theta$.
\end{subsec}

\begin{theorem} \label{T:DualGaussThakurSums/Zetas}
Fix an evaluative character $\chi_{\varepsilon} : A \to \FF_{q^d}$ modulo~$\wp$ as in \S\ref{SS:evalchars}.
\begin{alphenumerate}
\item \textup{(Gazda-Maurischat \cite{GazdaMaurischat21}*{Thm.~5.2})} For $\lambda \in \Lambda_{\phi}$,
\[
\bschi_{\varepsilon} \bigl( \omega_{\lambda}(t) \bigr) = G_{\lambda}(\chi_{\varepsilon}; \wp'(\theta)) \quad \text{and} \quad
\omega_{\lambda}(\varepsilon) = g_{\lambda}(\chi_{\varepsilon}; \wp'(\theta)).
\]
\item For $\eta \in \Hom_A(\Lambda_{\phi},A)$,
\[
\bschi_{\varepsilon} \bigl( \zeta_{\eta}(t) \bigr) = G_{\eta}^*(\chi_{\varepsilon}; \wp'(\theta))  \quad \text{and} \quad
\zeta_{\eta}(\varepsilon) = g_{\eta}^*(\chi_{\varepsilon}; \wp'(\theta)).
\]
\end{alphenumerate}
\end{theorem}

\begin{proof}
In principle the proof of (b) can follow similar lines to \cite{GazdaMaurischat21}*{Thm.~5.2}, though since we are restricting to the case $A = \FF_q[\theta]$, rather than the case of general~$A$ as in~\cite{GazdaMaurischat21}, we instead appeal to identities of Angl\`es and Pellarin \cites{AnglesPellarin14, AnglesPellarin15, Pellarin12} to obtain a more direct proof.

We first note that the second equality in (b) follows from the first. Next we consider
\[
1 \otimes \be_{\eta}^* : \FF_{q^d} \otimes_{\FF_q} k_{\infty} \to \FF_{q^d} \otimes_{\FF_q} \phi^*(\C),
\]
which is $\FF_{q^d} \otimes_{\FF_q} A$-linear and continuous with respect to the sup-norm. Then we have
\begin{equation} \label{E:chizetaexpand}
\bschi_{\varepsilon} \bigl( \zeta_\eta(t) \bigr)
= (1 \otimes \be_{\eta}^*) \biggl( \sum_{m=0}^{\infty} \varepsilon^m \otimes \theta^{-m-1} \biggr)
= (1 \otimes \be_{\eta}^*) \biggl( \frac{1}{1 \otimes \theta - \varepsilon \otimes 1} \biggr).
\end{equation}
On the other hand by~\eqref{E:Gcharab} and the definition of $\chi_{\varepsilon}$,
\begin{align} \label{E:Geta*expand}
G_{\eta}^*(\chi_{\varepsilon}; \wp'(\theta))
&= -(\wp'(\varepsilon) \otimes 1) \sum_{z \in \FF_{\wp}^{\times}} z(\varepsilon)^{-1} \otimes \be_{\eta}^*\bigl( z(\theta)/\wp(\theta) \bigr) \\
&= -(1 \otimes \be_{\eta}^*) \biggl( \bigl( \wp'(\varepsilon) \otimes \wp(\theta)^{-1} \bigr) \sumprime{\substack{z \in A \\ \deg z < d}} z(\varepsilon)^{-1} \otimes z(\theta) \biggr). \notag
\end{align}

Now for independent variables $x$, $y$ over $\FF_q$, it follows from \cite{AnglesPellarin14}*{Prop.~10} that, for $k \geq 0$, we have in $\FF_q(x)[y] = \FF_q(x) \otimes_{\FF_q} \FF_q[y]$,
\begin{equation} \label{E:PellarinSum}
\sum_{\substack{z \in A \\ \deg z = k}} z(x)^{-1} \otimes z(y)
= -\frac{(1 \otimes y - x \otimes 1)(1 \otimes y - x^q \otimes 1) \cdots (1 \otimes y - x^{q^{k-1}} \otimes 1)}{(x-x^q)(x-x^{q^2}) \cdots (x-x^{q^k}) \otimes 1}.
\end{equation}
Moreover, an induction implies that
\[
\sumprime{\substack{z \in A \\ \deg z \leq k}} z(x)^{-1} \otimes z(y)
= -\frac{(1 \otimes y - x^q \otimes 1)(1 \otimes y - x^{q^2} \otimes 1) \cdots (1 \otimes y - x^{q^{k}} \otimes 1)}{(x-x^q)(x-x^{q^2}) \cdots (x-x^{q^k}) \otimes 1}.
\]
(The reader will also observe that these are partial sums and partial products in the specialization of Pellarin's $L$-series at $s=1$ in~\cite{Pellarin12}*{Thm.~1}.) In our setting, this implies
\begin{align*}
\sumprime{\substack{z \in A \\ \deg z < d}} z(\varepsilon)^{-1} \otimes z(\theta)
&= -\frac{(1 \otimes \theta - \varepsilon^q \otimes 1)(1 \otimes \theta - \varepsilon^{q^2} \otimes 1) \cdots (1 \otimes \theta - \varepsilon^{q^{d-1}} \otimes 1)}{(\varepsilon - \varepsilon^q)(\varepsilon - \varepsilon^{q^2}) \cdots (\varepsilon - \varepsilon^{q^{d-1}}) \otimes 1} \\
&= -\frac{1 \otimes \wp(\theta)}{(1 \otimes \theta - \varepsilon \otimes 1)( \wp'(\varepsilon) \otimes 1)}.
\end{align*}
Combining this with~\eqref{E:chizetaexpand}--\eqref{E:Geta*expand}, we see that
\[
G_{\eta}^*(\chi_{\varepsilon}; \wp'(\theta))
= (1 \otimes \be_{\eta}^*) \biggl( \frac{1}{1 \otimes \theta - \varepsilon \otimes 1} \biggr) = \bschi_{\epsilon} \bigl( \zeta_{\eta}(t) \bigr).
\qedhere
\]
\end{proof}

\begin{subsec}{Twists of (dual) Gauss-Thakur sums} \label{ss:GaussSumTwists}
For $\lambda \in \Lambda_{\phi}$, $\eta \in \Hom_A(\Lambda_{\phi},A)$, $k\in \ZZ$, and $a\in A$, let 
\[
 g^{[k]}_\lambda(\chi;a) \assign -\sum_{z \in \FF_{\wp}^{\times}} \chi(z)^{-1} \be_{\lambda}(az/\wp)^{q^k}, \quad
 g_\eta^{*[k]}(\chi;a) \assign -\sum_{z \in \FF_{\wp}^{\times}} \chi(z)^{-1} \be_{\eta}^*(az/\wp)^{q^k}.
\]
This definition can be made for the tensored version of Gauss-Thakur sums as well, but we do not need this generality. The following corollary shows that twists of Gauss-Thakur sums satisfy a similar special function property to that of Gauss-Thakur sums in~\eqref{E:gcharab}.
\end{subsec}

\begin{corollary} \label{C:GaussSumsTwistsSomewhatSpecial}
Fix a multiplicative character $\chi \in \widehat{\FF}_{\wp}^{\times}$, $k \in \ZZ$, and $a \in A$. For $\lambda \in \Lambda_{\phi}$ and $\eta \in \Hom_A(\Lambda_{\phi},A)$, the following hold.
\begin{alphenumerate}
\item If $\chi$ is not evaluative as in~\S\ref{SS:evalchars}, then $g_{\lambda}^{[k]}(\chi;a)=0$ and $g_{\eta}^{*[k]}(\chi;a)=0$.
\item We have $g_{\lambda}^{[k]}(\chi; a)= \chi(a) g_{\lambda}^{[k]}(\chi;1)$ and
$g_{\eta}^{*[k]}(\chi; a)=\chi(a) g_{\eta}^{*[k]}(\chi;1)$.
\end{alphenumerate}   
\end{corollary}

\begin{proof}
We write only the dual case, as the argument is the same for both cases. Then
\begin{multline*}
g_{\eta}^{*[k]}(\chi; a) = -\sum_{z \in \FF_{\wp}^{\times}} \chi(z)^{-1} \be_{\eta}^*(az/\wp)^{q^k} = \biggl( -\sum_{z \in \FF_{\wp}^{\times}} \chi(z)^{-q^{-k}} \be_{\eta}^*(az/\wp) \biggr)^{q^k} \\
= g_{\eta}^*\bigl( \chi^{q^{-k}};a\bigr)^{q^k}
= \bigl( \chi^{q^{-k}}(a) g_{\eta}^*(\chi^{q^{-k}};1) \bigr)^{q^k}
= \chi(a) g_{\eta}^{*[k]}(\chi;1),
\end{multline*}
where the fourth equality follows from~\eqref{E:gcharab} and the last by essentially undoing the first and second. This proves (b). For (a), if $\chi$ is not evaluative, then neither is~$\chi^{q^{-k}}$, in which case
the above calculation and Proposition~\ref{P:vanishing} imply the desired vanishing.
\end{proof}

\begin{corollary} \label{C:GaussSumsTwistsAndZetas}
Fix an evaluative character $\chi_{\varepsilon} : A \to \FF_{q^d}$ modulo~$\wp$ as in \S\ref{SS:evalchars}. Then for $k \in \ZZ$ the following hold. For $\lambda \in \Lambda_{\phi}$ and $\eta \in \Hom_A(\Lambda_{\phi},A)$,
\[
\omega_{\lambda}^{(k)}(\varepsilon) = g_{\lambda}^{[k]}(\chi_{\varepsilon}; \wp'(\theta)), \quad
\zeta_{\eta}^{(k)}(\varepsilon) = g_{\eta}^{*[k]}(\chi_{\varepsilon}; \wp'(\theta)).
\]
\end{corollary}

\begin{proof}
Once again, the proof is similar for both cases so we proceed with only the dual case. Working as in the proof of the previous corollary,
\[
g_{\eta}^{*[k]}(\chi_{\varepsilon}; \wp'(\theta)) = g_{\eta}^{*}(\chi_{\varepsilon^{q^{-k}}}; \wp'(\theta))^{q^k}
= \bigl(\zeta_\eta(\varepsilon^{q^{-k}}) \bigr)^{q^k}
=\zeta_\eta ^{(k)}(\varepsilon).
\qedhere
\]
\end{proof}

\section{Duality of Gauss-Thakur sums} \label{S:MatrixGaussSection}

As one might expect the duality results of Ferraro~\cite{Ferraro25} in Theorems~\ref{T:FerraroUniv} and~\ref{T:Ferraropairing} induce a duality between Gauss-Thakur sums and dual Gauss-Thakur sums. This is also reflected in the interconnections between rigid analytic trivializations for the $t$-motive and $t$-comotive of a Drinfeld module, as in Lemma~\ref{L:PsiUpsilon} and Proposition~\ref{P:PsiNUpsilon11}. We explore this duality in detail in this section, and our first result is the following.

\begin{theorem} \label{T:GaussSumduality0}
Fix $\wp \in A_+$ irreducible, $\deg \wp =d$, and let $\varepsilon \in \FF_{q^d} \subseteq \C$ be a root of~$\wp$. For $1 \leq i \leq r$ the following hold.
\begin{alphenumerate}
\item We have
\[
\sum_{i=1}^r g_i^*(\chi_{\varepsilon}) g_i(\chi_{\varepsilon}) = \frac{1}{\wp'(\varepsilon)^2(\theta - \varepsilon)}.
\]
\item For $0 \leq j, k \leq r-1$ with $1 \leq j+k \leq r-1$,
\[
\sum_{i=1}^r g_i^{*[-j]}(\chi_\varepsilon) g_i^{[k]}(\chi_{\varepsilon}) = 0.
\]
\end{alphenumerate}
\end{theorem}

\begin{proof}
Let $\bsomega_{\phi} = (\omega_1, \dots, \omega_r)^{\tr}$ be the first column of $\Upsilon$ as in \eqref{E:Upsilondef} and $\bszeta_{\phi} = (\zeta_1, \dots, \zeta_r)$ be the first row of $\Psi$ as in \eqref{E:zetadef}. By Theorem~\ref{T:DualGaussThakurSums/Zetas} and Corollaries~\ref{C:GaussSumsTwistsSomewhatSpecial} and~\ref{C:GaussSumsTwistsAndZetas},
\[
\zeta_i^{(j)}(\varepsilon) = g_i^{*[j]}(\chi_{\varepsilon};\wp'(\theta)) = \wp'(\varepsilon)g_i^{*[j]}(\chi_{\varepsilon}), \quad
\omega_i^{(k)}(\varepsilon) = g_i^{[k]}(\chi_{\varepsilon};\wp'(\theta)) = \wp'(\varepsilon) g_i^{[k]}(\chi_{\varepsilon}),
\]
for $1 \leq i \leq r$ and $j$, $k \in \ZZ$. Therefore, for $j$, $k$, with $0 \leq j, k \leq r-1$, we have
\[
\sum_{i=1}^r g_i^{*[-j]} (\chi_{\varepsilon}) g_i^{[k]} (\chi_{\varepsilon}) = \frac{1}{\wp'(\varepsilon)^2} \cdot \bszeta_{\phi}^{(-j)}(\varepsilon) \cdot \bsomega_{\phi}^{(k)}(\varepsilon)
= \frac{1}{\wp'(\varepsilon)^2} \cdot \bigl[ \Psi \Upsilon \bigr]_{j+1,k+1} \Big|_{t=\varepsilon},
\]
and the result follows from Lemma~\ref{L:PsiUpsilon}.
\end{proof}

\begin{subsec}{Matrices of (dual) Gauss-Thakur sums} \label{SS:GTMatrices}
We next prove a more uniform result in Theorem~\ref{T:GaussSumduality1}, which does not depend on the choice of $\varepsilon$. For $\chi \in \widehat{\FF}_q^{\times}$ and $a \in A$, we recall the definitions of $\fG(\chi;a)$ and $\fG^*(\chi;a) \in \Mat_{r\times r}(\oK)$ from \eqref{E:frakG}--\eqref{E:frakG*}. For reasons of space we do not repeat their definitions here. Keeping with previous conventions, we set $\fG(\chi) \assign \fG(\chi;1)$ and $\fG^*(\chi) \assign \fG^*(\chi;1)$. We note as in \S\ref{SS:GTSMatrices}, Theorem~\ref{T:DualGaussThakurSums/Zetas} and Corollaries~\ref{C:GaussSumsTwistsSomewhatSpecial} and~\ref{C:GaussSumsTwistsAndZetas} imply
\begin{equation} \label{E:matrixGaussZetas}
\fG(\chi_{\varepsilon};\wp'(\theta)) = \wp'(\varepsilon)\fG(\chi_{\varepsilon}) = \Upsilon(\varepsilon), \quad
\fG^*(\chi_{\varepsilon};\wp'(\theta)) = \wp'(\varepsilon) \fG^*(\chi_{\varepsilon}) = J\Psi^{(r-1)}(\varepsilon).
\end{equation}
For $a \in A$, define $\Gamma(\phi;a) \in \Mat_{rd \times \NN}(\oK)$ and $\Gamma^*(\phi^*;a) \in \Mat_{\NN \times rd}(\oK)$ as in \eqref{E:Gammadef}--\eqref{E:Gamma*def}. Recall also the definitions of $\Gamma_{\ell}(\phi;a)$ and $\Gamma_{\ell}^*(\phi^*;a)$ from \S\ref{SS:GTSMatrices}, and set $\Gamma(\phi) \assign \Gamma(\phi;1)$ and $\Gamma^*(\phi^*) \assign \Gamma^*(\phi^*;1)$. Further for $\beta \in k$ that is regular at $\varepsilon$, recall
\[
\bQ(\beta) \assign \diag( \beta(\varepsilon) I_r, \dots, \beta(\varepsilon^{q^{d-1}}) I_r) \in \Mat_{rd \times rd}(\FF_{q^d}),
\]
as in~\eqref{E:bQintro}. For $a$, $b \in A$, we note that $\Gamma(\phi;ab) = \bQ(b)\Gamma(\phi;a)$ and $\Gamma^*(\phi^*;ab) =\Gamma^*(\phi^*;a) \bQ(b)$ by Corollary~\ref{C:GaussSumsTwistsSomewhatSpecial}.
Also, by Theorem~\ref{T:DualGaussThakurSums/Zetas} and Corollary~\ref{C:GaussSumsTwistsAndZetas},
\begin{equation} \label{E:matrixgaussandomega}
\Gamma(\phi;a) = \bQ(a/\wp') \Gamma(\phi;\wp'(\theta))
= \bQ(a/\wp') \begin{pNiceMatrix}
\Upsilon(\varepsilon) & \Upsilon^{(r)}(\varepsilon) & \Cdots[shorten-end=10pt] \\
\Vdots & \Vdots & \\
\Upsilon(\varepsilon^{q^{d-1}}) & \Upsilon^{(r)}(\varepsilon^{q^{d-1}}) & \Cdots[shorten-end=10pt]
\end{pNiceMatrix}
\end{equation}
and
\begin{equation} \label{E:matrixgaussandzeta}
\Gamma^*(\phi^*;a) = \Gamma^*(\phi^*;\wp'(\theta)) \bQ(a/\wp')
= \begin{pNiceMatrix}
J\Psi^{(r-1)}(\varepsilon) & \Cdots & J\Psi^{(r-1)}(\varepsilon^{q^{d-1}}) \\
J\Psi^{(2r-1)}(\varepsilon) & \Cdots & J\Psi^{(2r-1)}(\varepsilon^{q^{d-1}}) \\
\Vdots[shorten-end=5pt] &  & \Vdots[shorten-end=5pt]
\end{pNiceMatrix}
\bQ(a/\wp').
\end{equation}
Lastly, for $\beta \in k$ recall $\bP(\beta) \in \Mat_{\NN \times \NN}(\oK)$ as in~\eqref{E:bPintro}, which written more explicitly, for $i$, $j \geq 1$, satisfies
\begin{equation}
[\bP(\beta)]_{ij} = \begin{cases}
   {\displaystyle \brac{\beta}{i-j}_\phi^{q^{j-1}}} & \textup{if $i-j \geq 0$,} \\
   {\displaystyle -\brac{\beta}{j-i}_{\phi^*}^{q^{j-1}}} & \textup{if $i-j < 0$.}
\end{cases}
\end{equation}
Similarly, let $\bP_{\ell}(\beta) \in \Mat_{r\ell \times r\ell}(\oK)$ be the truncation of $\bP(\beta)$ after $r\ell$ rows and columns. 

The following theorem is the main result of this paper. It is similar to~\cite{Ferraro25}*{Thm.~7.18}, but made more explicit in this case where $A=\FF_q[\theta]$. Namely, it is a duality result between the Gauss-Thakur sums and their twists present in the $\Gamma(\phi)$ matrix and the dual Gauss-Thakur sums and their twists present in the $\Gamma^*(\phi)$ matrix.
Now, for $x = \sum_{j \geq j_0} c_j \theta^{-j} \in k_{\infty}$, we write $\langle x \rangle \assign \sum_{j \geq 1} c_j \theta^{-j} \in \theta^{-1} \power{\FF_q}{\theta^{-1}}$. For $a \in A$, we take $a \bmod \wp$ for the remainder of $a$ divided by $\wp$ with $\deg (a \bmod \wp) < \deg \wp$. Thus, $\langle a/\wp \rangle = (a \bmod \wp)/\wp$.
\end{subsec}

\begin{theorem} \label{T:GaussSumduality1}
With the matrices defined as above, we have the following identities.
\begin{alphenumerate}
\item For $a$, $b \in A$,
\[
\Gamma^*(\phi^*;a) \cdot \bQ(\wp') \cdot \Gamma(\phi;b)=\bP (\langle ab/\wp \rangle ).
\]
\item For $a$, $b \in A$ and $\ell \in \NN$,
\[
\Gamma_\ell^*(\phi^*;a) \cdot \bQ(\wp') \cdot \Gamma_\ell (\phi;b) = \bP_\ell (\langle ab/\wp \rangle ).
\]
\item In particular,
\[
\Gamma^*(\phi^*) \cdot \bQ(\wp') \cdot \Gamma(\phi) = \bP(1/\wp).
\]
\end{alphenumerate}
\end{theorem}

\begin{proof}
Parts (b) and (c) follow directly from (a). To prove (a), we fix $i$, $j \geq 1$ and note
\begin{equation} \label{E:productred}
\bigl[\Gamma^*(\phi^*;a) \bQ(\wp') \Gamma(\phi;b) \bigr]_{ij}
 = \bigl[ \Gamma^*(\phi^*;\wp'(\theta)) \bQ(ab/\wp') \Gamma(\phi;\wp'(\theta)) \bigr]_{ij}.
\end{equation}
We note that $J\Psi^{(r-1)} = ( \bszeta_{\phi}, \bszeta_{\phi}^{(1)}, \dots, \bszeta_{\phi}^{(r-1)})^{\tr}$,
and so~\eqref{E:matrixgaussandzeta} implies that the $i$-th row of $\Gamma^*(\phi^*;\wp'(\theta))$ is $( \bszeta_{\phi}^{(i-1)}(\varepsilon),\, \bszeta_{\phi}^{(i-1)}(\varepsilon^q),\, \dots,\, \bszeta_{\phi}^{(i-1)}(\varepsilon^{q^{d-1}}))$. Likewise, by~\eqref{E:matrixgaussandomega}, the $j$-th column of $\Gamma(\phi;\wp'(\theta))$ is $( \bsomega_{\phi}^{(j-1)}(\varepsilon), \bsomega_{\phi}^{(j-1)}(\varepsilon^q), \dots, \bsomega_{\phi}^{(j-1)}(\varepsilon^{q^{d-1}}) )^{\tr}$.
From the definition of $\bQ(\beta)$ as a diagonal matrix, it follows that
\begin{align} \label{E:firstred}
\bigl[ \Gamma^*(\phi^*;\wp'(\theta)) \cdot{} &\bQ(ab/\wp') \cdot \Gamma(\phi;\wp'(\theta)) \bigr]_{ij} \\
&= \sum_{k=1}^d \frac{a(\varepsilon^{q^{k-1}}) b(\varepsilon^{q^{k-1}})}{\wp'(\varepsilon^{q^{k-1}})} \cdot \bszeta_{\phi}^{(i-1)} (\varepsilon^{q^{k-1}}) \cdot \bsomega_{\phi}^{(j-1)} (\varepsilon^{q^{k-1}}), \notag \\
\intertext{and since $\bszeta_{\phi}^{(i-1)}\bsomega_{\phi}^{(j-1)}$ is the upper-left entry of $\Psi^{(i-1)} \Upsilon^{(j-1)}$,}
&= \sum_{k=1}^d \frac{a(\varepsilon^{q^{k-1}})b(\varepsilon^{q^{k-1}})}{\wp'(\varepsilon^{q^{k-1}})} \cdot
\bigl[ \Psi^{(i-1)}(\varepsilon^{q^{k-1}}) \cdot \Upsilon^{(j-1)}(\varepsilon^{q^{k-1}}) \bigr]_{11} \notag \\
&= \sum_{k=1}^d \frac{a(t)b(t)\bmod \wp(t)}{\wp'(t)} \cdot \bigl[ \Psi^{(i-j)}\Upsilon \bigr]^{(j-1)}_{11} \Big|_{t=\varepsilon^{q^{k-1}}}. \notag
\end{align}
By their definitions in \S\ref{SS:RATs}, the entries of $\Psi$ and $\Upsilon$ are in $\TT \cap \power{\oK}{t}$. Thus, $[\Psi^{(i-j)}\Upsilon]_{11}^{(j-1)} = \sum_{m \geq 0} c_m t^m$ with $c_m \in \oK$ and $\inorm{c_m} \to 0$. Letting $\rho(t) = a(t)b(t) \bmod \wp(t) \in \FF_q[t]$,
\begin{equation} \label{E:secondred}
\sum_{k=1}^d \frac{\rho(t)}{\wp'(t)} \cdot \bigl[ \Psi^{(i-j)}\Upsilon \bigr]^{(j-1)}_{11} \Big|_{t=\varepsilon^{q^{k-1}}}
= \sum_{m=0}^{\infty} c_m\sum_{k=1}^d \frac{\rho(t)}{\wp'(t)} t^m \Big|_{t=\varepsilon^{q^{k-1}}}.
\end{equation}
Now as a function on $\PP^1/\oFF_q$, $(\rho(t)/\wp(t))t^m$ can only have poles at $\infty$ and at the zeros of~$\wp(t)$. Hence, by~\eqref{E:residue} and the residue theorem,
\begin{equation} \label{E:residue1}
\rR_t \biggl( \frac{\rho(t)}{\wp(t)}t^m \biggr) = -\Res_\infty \biggl( \frac{\rho(t)}{\wp(t)}t^m \, dt \biggr)
= \sum_{k=1}^{d} \Res_{t=\varepsilon^{q^{k-1}}} \biggl( \frac{\rho(t)}{\wp(t)}t^m \, dt \biggr).
\end{equation}
As $(\rho(t)/\wp(t))t^m$ has at most simple poles at $\varepsilon^{q^{k-1}}$, $1 \leq k \leq d$, we have
\begin{equation} \label{E:residue2}
\Res_{t=\varepsilon^{q^{k-1}}} \biggl( \frac{\rho(t)}{\wp(t)} t^m\, dt \biggr)
= \frac{(t-\varepsilon^{q^{k-1}})\rho(t)}{\wp(t)}t^m \bigg|_{t=\varepsilon^{q^{k-1}}}
= \frac{\rho(t)}{\wp'(t)}t^m \bigg|_{t=\varepsilon^{q^{k-1}}}.
\end{equation}
Combining \eqref{E:productred}--\eqref{E:residue2}, we have
\[
\bigl[ \Gamma^*(\phi^*;a) \bQ(\wp') \Gamma(\phi;b) \bigr]_{ij}
= \sum_{m=0}^{\infty}c_m\rR_t \biggl( \frac{\rho(t)t^m}{\wp(t)} \biggr). 
\]
Note that $\rR_t((\rho(t)/\wp(t)) t^m) \in \FF_q$ is the coefficient of $t^{-m-1}$ in the $\laurent{\FF_q}{t^{-1}}$-expansion of $\rho(t)/\wp(t)$. Thus writing $\rho(t)/\wp(t) = \sum_{m=\delta}^{\infty} a_m t^{-m}$, with $\delta = -\deg(\rho(t)/\wp(t))$ and $a_m \in \FF_q$, and noting that $\delta >0$, we obtain
\[
\bigl[ \Gamma^*(\phi^*;a)\bQ(\wp')\Gamma(\phi;b) \bigr]_{ij}
= \sum_{m=\delta}^{\infty} a_{m} c_{m-1}
= \sum_{m=\delta}^{\infty} a_m \begin{cases}
{\displaystyle \brac{\theta^{-m}}{i-j}_{\phi}^{q^{j-1}},} & \textup{if $i-j \geq 0$,} \\
{\displaystyle -\brac{\theta^{-m}}{j-i}_{\phi^*}^{q^{j-1}},} & \textup{if $i-j < 0$,}
\end{cases}
\]
where the second equality follows from the definition of $c_m$ and Proposition~\ref{P:PsiNUpsilon11}. Because $\phi : k_{\infty} \to \laurent{\C}{\tau}$ and $\phi^* : k_{\infty} \to \laurent{\C}{\tau}$ have continuous coefficient functions,
\[
\bigl[ \Gamma^*(\phi^*;a)\bQ(\wp')\Gamma(\phi;b) \bigr]_{ij}
= \begin{cases}
{\displaystyle \brac{ \rho(\theta)/\wp(\theta)}{i-j}_{\phi}^{q^{j-1}},}& \textup{if $i-j \geq 0$,} \\
{\displaystyle -\brac{ \rho(\theta)/\wp(\theta)}{j-i}_{\phi^*}^{q^{j-1}},} & \textup{if $i-j < 0$,}
\end{cases}
\]
as desired.
\end{proof}

\begin{remark} \label{R:invertibleP}
Note that for $-r \deg \wp < i-j < 0$, $[\bP(1/\wp)]_{ij} = -\brac{1/\wp}{j-i}_{\phi^*}^{q^{j-1}}=0$ by conventions established in \S\ref{SS:phistarext}. This implies that for $\ell \leq \deg \wp$, and $a$, $b \in A$ such that $ab \bmod \wp=1$, the matrix
\[
\Gamma_{\ell}^*(\phi^*;a) \cdot \bQ(\wp') \cdot \Gamma_{\ell}(\phi,b) = \bP_{\ell} ( 1/\wp) \quad \in \GL_{r\ell}(\oK),
\]
is lower triangular with non-zero diagonal entries $\wp^{-q^{j-1}}$, and thus is invertible.
\end{remark}

\section{Carlitz module} \label{S:CarlitzApp}

In this section we restrict to the case of the Carlitz module, $\sC : A \to A[\tau]$ defined by $\sC_{\theta} = \theta + \tau$.
In this case the Gauss-Thakur sums $g(\chi_{\varepsilon};a)$ are precisely the sums defined and investigated by Thakur in~\cite{Thakur88}.

\begin{subsec}{Essentials of the Carlitz module}
The Carlitz module is rank~$1$, and its fundamental period, the \emph{Carlitz period}, is
\begin{equation} \label{E:pitilde}
\tpi = -(-\theta)^{q/(q-1)} \prod_{n=1}^{\infty} \Bigl( 1 - \theta^{1-q^n} \Bigr)^{-1} \quad \in k_{\infty}\bigl( (-\theta)^{1/(q-1)} \bigr).
\end{equation}
For more details about the Carlitz module, see \citelist{\cite{Goss}*{Ch.~3} \cite{Papikian}*{Ch.~3} \cite{Thakur}*{Ch.~2}}. The matrices $\Upsilon$ and $\Psi$ are then $1 \times 1$, given by
\begin{equation} \label{E:omegaCOmegaC}
\Upsilon = \omega_{\sC} = (-\theta)^{1/(q-1)} \prod_{n=0}^{\infty} \biggl( 1 - \frac{t}{\theta^{q^n}} \biggr)^{-1}, \quad
\Psi = -\Omega_{\sC} = -(-\theta)^{-q/(q-1)} \prod_{n=1}^{\infty} \biggl( 1 - \frac{t}{\theta^{q^n}} \biggr),
\end{equation}
for a fixed choice of $(-\theta)^{1/(q-1)}$.
The function $\omega_{\sC}$ is the Anderson-Thakur function defined in~\cite{AndThak90}*{\S 2.5}, and $\Omega_{\sC}$ was defined in~\cite{ABP04}*{\S 3.1}. In the context of \S\ref{SS:AGF} and \S\ref{SS:eetastar},
\[
\omega_{\tpi}(t) = \omega_{\sC}(t) = \sum_{m=0}^{\infty} \be_{\tpi}(1/\theta^{m+1}) t^m, \quad \zeta_{\teta}(t) = -\Omega_{\sC}(t) = \sum_{m=0}^{\infty} \be_{\teta}^*(1/\theta^{m+1}) t^m,
\]
where $\teta : \Lambda_{\sC} \to A$ is defined by $\teta(a\tpi) = a$ (cf.~\citelist{\cite{ABP04}*{\S 5.1} \cite{DavisP25}*{\S 2.5}}). We also have 
$\omega_{\sC}^{(1)} = (t-\theta) \omega_{\sC}$, $\Omega_{\sC}^{(-1)} = (t-\theta) \Omega_{\sC}$, and $\Omega_{\sC} = 1/\omega_{\sC}^{(1)}$,
which then provide directly Proposition~\ref{P:PsiNUpsilon11} in this case that $-\Omega_{\sC} \cdot \omega_{\sC} = 1/(\theta - t)$.
\end{subsec}

\begin{remark} \label{R:oldestar}
Throughout this section we let $\be = \be_{\tpi}$ and $\be^* = \be_{\teta}$, though as a warning to the reader, we note that this notation for $\be^*$ is different from that in~\cite{ABP04}*{\S 5.2}. Indeed if $\be_1^* : k_{\infty} \to \sC^*(\C)$ is the function defined in~\cite{ABP04}, then $\be^*(x) = -\be_1^*(x)^q$. The negative sign arises because $\zeta_{\teta} = -\Omega_{\sC}$, and so also $\be^*$ differs by a sign from $\tbe$ in~\cite{DavisP25}*{\S 3.1}.
\end{remark}

\begin{subsec}{Gauss-Thakur sums and their duals for $\sC$} \label{SS:GTsumsCarlitz}
Fix $\wp \in A_+$, irreducible of degree~$d$, and let $\varepsilon \in \FF_{q^d} \subseteq \C$ be one of its roots. As presented in \S\ref{SS:Gauss-ThakurSums}, we have Gauss-Thakur sums, for $\chi \in \widehat{\FF}_q^{\times}$ and $a \in A$,
\[
g(\chi;a) \assign g_{\tpi}(\chi;a), \quad
g^*(\chi;a) \assign g_{\teta}^*(\chi;a).
\]
As usual $g(\chi) = g(\chi;1)$ and $g^*(\chi) = g^*(\chi;1)$. By Theorem~\ref{T:DualGaussThakurSums/Zetas},
\begin{equation} \label{E:Carlitzomegaevals}
\wp'(\varepsilon) g(\chi_\varepsilon) = \omega_{\sC}(\varepsilon), \quad
\wp'(\varepsilon) g^*(\chi_\varepsilon) = -\Omega_{\sC}(\varepsilon),
\end{equation}
(the first originally due to Angl\`es and Pellarin~\cite{AnglesPellarin15}*{Thm.~2.9}).
Also by Theorem~\ref{T:GaussSumduality0}(a),
\begin{equation} \label{E:ggstaridentity}
g(\chi_{\varepsilon})g^*(\chi_{\varepsilon}) = \frac{1}{(\theta - \varepsilon) \wp'(\varepsilon)^2},
\end{equation}
from which we see in particular that $g(\chi_{\varepsilon})\neq 0$ and $g^*(\chi_{\varepsilon}) \neq 0$.

Now Thakur established several theorems for Gauss-Thakur sums in \cite{Thakur88}, which parallel results for classical Gauss sums. Let $K = k(\be(1/\wp))$, the $\wp$-th Carlitz cyclotomic field, and let $H=\FF_{q^d}K$. Then
$\Gal(H/k) \cong \Gal(\FF_{q^d}/\FF_q) \times \Gal(K/k)$, where $\Gal(\FF_{q^d}/\FF_q) = \langle \sigma \rangle$ is generated by the $q$-th power Frobenius and $\Gal(K/k) = \{ \rho_a : a \in (A/\wp)^{\times} \} \cong (A/\wp)^{\times}$. By \citelist{\cite{Goss95}*{Thm.~3.2} \cite{Goss}*{Thm.~1.7.11}}, we have $K = k(\be(1/\wp)) = k(\be^*(1/\wp))$, and moreover, by \citelist{\cite{Goss}*{Prop.~4.14.13} \cite{Papikian}*{Thm.~7.1.12} \cite{Poonen96}*{Cor.~11}}, for any $a$, $b \in (A/\wp)^{\times}$,
\begin{equation} \label{E:rhoCarlitz}
\rho_a ( \be(b/\wp) ) = \be(ab/\wp) = \sC_a(\be(b/\wp)), \quad
\rho_a^{-1} ( \be^*(b/\wp)) = \be^*(ab/\wp) = \sC_a^*( \be^*(b/\wp)),
\end{equation}
(see also \cite{DavisP25}*{Prop.~3.1.6}). Finally, for $0 \leq k \leq d-1$, let $\fp_k \subseteq \cO_H$ be the unique (totally ramified) prime ideal above $\theta - \varepsilon^{q^k}$ in $\FF_{q^d}[\theta]$, and normalize $\ord_{\fp_k}$ so that $\ord_{\fp_k}(\wp) = 1$.
\end{subsec}

\begin{theorem}[{Thakur~\cite{Thakur88}*{Thms.~I--IV}}] \label{T:ThakurI-IV}
Let $\wp \in A_+$ be irreducible of degree $d$, and fix a root $\varepsilon \in \FF_{q^d}$ of $\wp$. Let $K = k(\be(1/\wp))$ and $H = \FF_{q^d} K$.
\begin{alphenumerate}
\item Let $0 \leq k \leq d-1$. For $\sigma \in \Gal(H/k)$ as above, $\sigma(g(\chi_k)) = g(\chi_{k+1})$.
For $a \in (A/\wp)^\times$, $\rho_a(g(\chi_k)) = a(\varepsilon^{q^k}) g(\chi_k)$.
\item We have
\[
\prod_{k=0}^{d-1} g(\chi_k)^{q-1} = (-1)^d \wp.
\]
\item \textup{(Stickelberger's Theorem)} There is a unique root $(-\wp)^{1/(q^d-1)}$ in the completion $H_{\fp_0}$ such that $\ord_{\fp_0}( \be(1/\wp) - (-\wp)^{1/(q^d-1)}) > 1/(q^d-1)$. In this case, 
\[
\ord_{\fp_0} \biggl( \frac{g(\chi_{\varepsilon})}{(-\wp)^{1/(q^d-1)}} - 1 \biggr) > 0,
\]
and more generally, for $1 \leq k \leq d-1$,
\[
\ord_{\fp_0} \biggl( \frac{g(\chi_k)}{(-\wp)^{q^k/(q^d-1)}} - \frac{1}{(\varepsilon^{q^k} - \theta)(\varepsilon^{q^{k-1}} - \theta)^q \cdots (\varepsilon^q - \theta)^{q^{k-1}}} \biggr) > 0.
\]
\item As ideals in $\cO_H$, for $0 \leq k \leq d-1$,
\[
g(\chi_k) \cdot \cO_H = \prod_{j=0}^{d-1} \fp_j^{q^{(k-j) \bmod d}}.
\]
\item For $0 \leq k \leq d-1$, the degree of $g(\chi_k)$ at any infinite place of $H$ is $1/(q-1)$.
\end{alphenumerate}
\end{theorem}

The companion results for dual Gauss-Thakur sums for the Carlitz module are the following. They can be proved in the same manner as in \cite{Thakur88} directly from properties of $g^*(\chi_{\varepsilon})$, but in large part we instead apply Theorem~\ref{T:GaussSumduality0} via~\eqref{E:ggstaridentity}. The trickiest part to derive is the analogue of Stickelberger's theorem in~(c).

\begin{theorem} \label{T:dualThakurI-IV}
Let $\wp \in A_+$ be irreducible of degree $d$, and fix a root $\varepsilon \in \FF_{q^d}$ of $\wp$. Let $K = k(\be^*(1/\wp))=k(\be(1/\wp))$ and $H = \FF_{q^d} K$.
\begin{alphenumerate}
\item Let $0 \leq k \leq d-1$. For $\sigma \in \Gal(H/k)$ as above, $\sigma(g^*(\chi_k)) = g^*(\chi_{k+1})$.
For $a \in (A/\wp)^\times$, $\rho_a^{-1}(g^*(\chi_k)) = a(\varepsilon^{q^k}) g^*(\chi_k)$.
\item We have
\[
\prod_{k=0}^{d-1} g^*(\chi_k)^{q-1} = \frac{(-1)^d}{\wp^q}.
\]
\item \textup{(Stickelberger's Theorem)} For the same root $(-\wp)^{1/(q^d-1)}$ in Theorem \ref{T:ThakurI-IV}(c), we have $\ord_{\fp_0}(\be^*(1/\wp) + \wp'(\varepsilon)^{-1} (-\wp)^{-q^d/(q^d-1)} ) > -q^d/(q^d-1)$. In this case, 
\[
\ord_{\fp_0} \Bigl( g^*(\chi_{\varepsilon}) (-\wp'(\varepsilon))(-\wp)^{q^d/(q^d-1)} - 1 \Bigr) > 0,
\]
and more generally, for $1 \leq k \leq d-1$,
\[
\ord_{\fp_0} \biggl(\frac{g^*(\chi_k) (-\wp'(\varepsilon^{q^{k}})) (-\wp)^{q^{d+k}/(q^d-1)}} {(\varepsilon - \theta)^{q^k} (\varepsilon^q - \theta)^{q^{k-1}} \cdots (\varepsilon^{q^{k-1}} - \theta)^q} - 1 \biggr) > 0.
\]
\item As fractional ideals of $\cO_H$ in $H$, for $0 \leq k \leq d-1$,
\[
g^*(\chi_k) \cdot \cO_H = \prod_{j=0}^{d-1}\fp_j^{-q^{d - ((j-k)\bmod d)}}.
\]
\item For $0 \leq k \leq d-1$, the degree of $g^*(\chi_k)$ at any infinite place of $H$ is $-q/(q-1)$.
\end{alphenumerate}
\end{theorem}

\begin{proof}
Part~(a) follows from the definitions of $\sigma$, $\rho_a \in \Gal(H/k)$ and from~\eqref{E:rhoCarlitz}. For~(b), we combine~\eqref{E:ggstaridentity} with Theorem~\ref{T:ThakurI-IV}(b) and find
\[
\prod_{k=0}^{d-1} g^*(\chi_k)^{q-1} = \frac{1}{(-1)^d\wp} \cdot \prod_{k=0}^{d-1} \frac{1}{(\theta - \varepsilon^{q^k})^{q-1}} \cdot \frac{1}{\wp'(\varepsilon^{q^k})^{2(q-1)}}.
\]
In the product on the right, the first terms multiply together to give $\wp^{-(q-1)}$. The second terms multiply to give $\rN(\wp'(\varepsilon)^{2(q-1)})$, where $\rN = \rN^{\FF_{q^d}}_{\FF_q} : \FF_{q^d}^{\times} \to \FF_q^{\times}$ is the norm. But this is~$1$ as it is the $(q-1)$-st power of an element of $\FF_q^{\times}$.

For (c), substituting \eqref{E:ggstaridentity} into the first main inequality in Theorem~\ref{T:ThakurI-IV}(c), we obtain
\begin{equation} \label{E:partcsub}
\ord_{\fp_0} \biggl( \frac{1}{g^*(\chi_{\varepsilon}) (-\wp)^{1/(q^d-1)}(\theta-\varepsilon) \wp'(\varepsilon)^2} - 1 \biggr) > 0.
\end{equation}
We note that \eqref{E:ggstaridentity} implies $\ord_{\fp_0} (g(\chi_{\varepsilon})) + \ord_{\fp_0}(g^*(\chi_{\varepsilon})) = -\ord_{\fp_0}(\theta-\varepsilon) = -1$. Thus, as $\ord_{\fp_0}(g(\chi_{\varepsilon})) = 1/(q^d-1)$ by Theorem~\ref{T:ThakurI-IV}(c), $\ord_{\fp_0}(g^*(\chi_{\varepsilon})) = -q^d/(q^d-1)$.
Multiplying the expression in \eqref{E:partcsub} by $g^*(\chi_{\varepsilon}) (-\wp)^{q^d/(q^d-1)} \wp'(\varepsilon)$, which has $\fp_0$-adic valuation~$0$, yields
\[
\ord_{\fp_0} \biggl( g^*(\chi_{\varepsilon}) (-\wp'(\varepsilon))(-\wp)^{q^d/(q^d-1)} - \frac{\wp(\theta)}{(\theta-\varepsilon) \wp'(\varepsilon)} \biggr) > 0.
\]
Since $\wp(\theta) = \wp'(\varepsilon)(\theta - \varepsilon) + O((\theta - \varepsilon)^2)$ in $\FF_{q^d}[\theta]$, we then obtain
\begin{equation} \label{E:Stickelberger1}
\ord_{\fp_0} \Bigl( g^*(\chi_{\varepsilon}) (-\wp'(\varepsilon))(-\wp)^{q^d/(q^d-1)} - 1 \Bigr) > 0.
\end{equation}
Now for $k \in \ZZ$, since $\sC_{\theta}^* = \theta + \tau^{-1}$,
\begin{multline*}
\varepsilon^{q^k} g^*(\chi_k) = g^*(\chi_k; \theta) = -\sum_{z \in \FF_{\wp}^{\times}} z(\varepsilon^{q^k})^{-1} \be^*(\theta z/\wp) \\
= -\sum_{z \in \FF_{\wp}^{\times}} z( \varepsilon^{q^k})^{-1} \bigl( \theta \be^*(1/\wp) + \be^*(z/\wp)^{1/q} \bigr) 
= \theta g^*(\chi_k) + g^*(\chi_{k+1})^{1/q}.
\end{multline*}
Thus,
\begin{equation} \label{E:gstarCarltizrec}
g^*(\chi_k)^q = \frac{g^*(\chi_{k+1})}{(\varepsilon^{q^k} - \theta)^q}.
\end{equation}
Using the convention that an empty product is $1$, we see that~\eqref{E:Stickelberger1} is the $k=0$ version of the inequality in~(c). Thus proceeding by induction on $k$, assume the inequality in~(c) holds for $k< d-1$. Raising the expression in~(c) for $k$ to the $q$-th power and substituting in~\eqref{E:gstarCarltizrec}, we obtain
\[
\frac{g^*(\chi_k)^q (-\wp'(\varepsilon^{q^{k+1}})) (-\wp)^{q^{d+k+1}/(q^d-1)}}{(\varepsilon - \theta)^{q^{k+1}} (\varepsilon^q - \theta)^{q^{k}} \cdots (\varepsilon^{q^{k-1}} - \theta)^{q^2}}
= \frac{g^*(\chi_{k+1}) (-\wp'(\varepsilon^{q^{k+1}})) (-\wp)^{q^{d+k+1}/(q^d-1)}}{(\varepsilon - \theta)^{q^{k+1}} (\varepsilon^q - \theta)^{q^{k}} \cdots (\varepsilon^{q^k} - \theta)^q},
\]
which yields the desired inequality in~(c).

For (d), fix $0 \leq k \leq d-1$. For $0 \leq j \leq d-1$, we obtain from \eqref{E:ggstaridentity} that
\[
\ord_{\fp_j}(g^*(\chi_k)) + \ord_{\fp_j}(g(\chi_k))
= \begin{cases}
0 & \textup{if $j\neq k$,} \\
 -1 & \textup{if $j=k$,}
\end{cases}
\]
Using this together with Theorem~\ref{T:ThakurI-IV}(d), we find
\[
g^*(\chi_k) \cdot \cO_H
= \fp_k^{-q^d} \cdot \prod_{\substack{j=0 \\ j \neq k}}^{d-1} \fp_j^{-q^{(k-j) \bmod d}}=\prod_{j=0}^{d-1}\fp_j^{-q^{d - ((j-k)\bmod d)}}.
\]
Note that the powers on $\fp_j$ are affected by our normalization of the valuation $\ord_{\fp_j}$. Part~(e) follows directly from Theorem~\ref{T:ThakurI-IV}(e) and~\eqref{E:ggstaridentity}.
\end{proof}

\begin{remark}
One can show that with the chosen root $(-\wp)^{1/(q^d-1)}$ in Theorem~\ref{T:ThakurI-IV}(c),
\[
\ord_{\fp_0} \biggl( \be^*(1/\wp) + \frac{(-\wp)^{-q^d/(q^d-1)}}{\wp'(\varepsilon)}  \biggr) > -\frac{q^d}{q^d-1}.
\]
We have omitted the details, which are similar to those in~\cite{Thakur88}*{Lem.~II}, for brevity.
\end{remark}

\begin{subsec}{$\infty$-adic signs of Gauss-Thakur sums}
We have seen in Theorems~\ref{T:ThakurI-IV}(e) and~\ref{T:dualThakurI-IV}(e) that the $\infty$-adic degrees of Carlitz module Gauss-Thakur sums and their duals are explicit. These align with \cite{Thakur93a}*{Thm.~1.2}, where Thakur relates the $\infty$-adic size of $g(\chi_{\varepsilon})$ to the $\wp$-torsion on $\sC$ with largest norm. We can also develop information about the signs of $g(\chi_{\varepsilon})$ and $g^*(\chi_{\varepsilon})$. Suppose we have fixed a decomposition $\C^{\times} = \oFF_{q}^{\times} \times (1/\theta)^{\QQ} \times U_1(\C)$, where consistent roots of $\theta$ have been chosen and $U_1(\C)$ is the group of $1$-units. Projecting onto $\oFF_q^{\times}$ defines a sign function, $\sgn : \C^{\times} \to \oFF_q^{\times}$, that we fix throughout. In this way all rational powers of $\theta$ have sign~$1$, and moreover, the group $(1/\theta)^{\QQ} \times U_1(\C)$ of sign~$1$ elements is uniquely divisible. Also, by~\eqref{E:pitilde},
\begin{equation} \label{E:signpitilde}
\sgn(\tpi) = -(-1)^{q/(q-1)} = (-1)^{1/(q-1)},
\end{equation}
where throughout we fix this choice of root of $-1$ so that \eqref{E:pitilde}--\eqref{E:omegaCOmegaC} are consistent. Results on the signs of Gauss-Thakur sums are then due to Thakur.
\end{subsec}

\begin{proposition}[{Thakur~\cite{Thakur93a}*{Thm.~2.3}}] \label{P:gg*signs}
For $0 \leq k \leq d-1$,
\[
\sgn \bigl( g(\chi_k) \bigr)
= \frac{(-1)^{1/(q-1)}}{\wp'(\varepsilon^{q^k})},
\quad
\sgn \bigl( g^*(\chi_k) \bigr)
= \frac{1}{(-1)^{1/(q-1)} \wp'(\varepsilon^{q^k})}.
\]
\end{proposition}

\begin{proof}
For completeness we prove the case of $\sgn(g^*(\chi_{\varepsilon}))$, though it can be derived from~\cite{Thakur93a}*{Thm.~2.3} and~\eqref{E:ggstaridentity}.
The argument is similar in spirit to~\cite{Thakur93a}*{Thm.~2.3} (see also~\cite{Thakur}*{Thm.~4.4.1}). For $x \in k_{\infty}$ with $\inorm{x} < 1$, \cite{ABP04}*{Pf.~of Lem.~5.2.6} implies that $\inorm{\be^*(x)} = \inorm{\theta}^{-q^{-\deg x}/(q-1)}$,
(see also~\cite{DavisP25}*{Eq.~(2.5.9)}, keep in mind Remark~\ref{R:oldestar}). For $z \in A$, $\deg z < \deg \wp$, it follows that $\deg \be^*(z/\wp) = -q^{d-\deg z}/(q-1)$ is maximized when $\deg z = d-1$, and in that case
\begin{equation} \label{E:signtemp}
\sgn\bigl( \be^*(z/\wp) \bigr) = \sgn(z) \sgn \bigl( \be^*(\theta^{-1}) \bigr) = \sgn(z) \sgn \bigl( -(-\theta)^{-q/(q-1)} \bigr),
\end{equation}
by \eqref{E:omegaCOmegaC}. Therefore,
\[
g^*(\chi_{\varepsilon}) = - \sumprime{\substack{z \in A \\ \deg z < d}} z(\varepsilon)^{-1} \be^*(z/\wp)
= -\sum_{\substack{z \in A \\ \deg z = d-1}} z(\varepsilon)^{-1} \be^*(z/\wp) + E, \quad \deg E < -\frac{q}{q-1}.
\]
By Theorem~\ref{T:dualThakurI-IV}(e), this initial sum is non-zero, so by\eqref{E:signpitilde} and~\eqref{E:signtemp},
\[
\sgn\bigl( g^*(\chi_{\varepsilon}) \bigr)
= - \sgn \Biggl( \sum_{\substack{z \in A \\ \deg z = d-1}} z(\varepsilon)^{-1} \be^*(z/\wp) \Biggr)
= \frac{-1}{(-1)^{1/(q-1)}} \cdot \sgn \Biggl( \sum_{\substack{z \in A \\ \deg z = d-1}} \frac{z(\theta)}{z(\varepsilon)} \Biggr),
\]
By~\eqref{E:PellarinSum}, this latter sum is $-(\theta - \varepsilon) \cdots (\theta - \varepsilon^{q^{d-2}})/((\varepsilon - \varepsilon^q) \cdots (\varepsilon - \varepsilon^{q^{d-1}}))$, and so $\sgn(g^*(\chi_{\varepsilon})) = (-1)^{-1/(q-1)}/\wp'(\varepsilon)$. Replacing $\varepsilon$ by $\varepsilon^{q^k}$ yields the general result.
\end{proof}

\begin{subsec}{Matrices of (dual) Gauss-Thakur sums for the Carlitz module} \label{SS:GammaCGamma*C}
In the context of Theorem~\ref{T:GaussSumduality1}, consider the $d\times d$ matrices,
\[
\Gamma_d(\sC) = \begin{pNiceMatrix}
g(\chi_0) & \Cdots & g^{[d-1]}(\chi_0) \\
\Vdots & & \Vdots \\
g(\chi_{d-1}) & \Cdots & g^{[d-1]}(\chi_{d-1})
\end{pNiceMatrix},
\;
\Gamma_d^*(\sC^*) = \begin{pNiceMatrix}
g^*(\chi_0) & \Cdots & g^*(\chi_{d-1}) \\
\Vdots & & \Vdots \\
g^{*[d-1]}(\chi_0) & \Cdots & g^{*[d-1]}(\chi_{d-1})
\end{pNiceMatrix}.
\]
Theorem~\ref{T:GaussSumduality1}(c) and Remark~\ref{R:invertibleP} imply
\begin{equation} \label{E:Carlitz7.5}
\Gamma_d^*(\sC^*) \cdot \bQ(\wp') \cdot
\Gamma_d(\sC)
= \begin{pNiceMatrix}
1/\wp & 0 & 0 & \Cdots & 0 \\
\brac{1/\wp}{1}_{\sC} & 1/\wp^{q} & 0 & \Cdots & 0 \\
\Vdots & \Ddots & \Ddots[draw-first] & & \Vdots \\[5pt]
\brac{1/\wp}{d-2}_{\sC} & &  & 1/\wp^{q^{d-2}} & 0 \\
\brac{1/\wp}{d-1}_{\sC} & \brac{1/\wp}{d-2}_{\sC}^q & \Cdots & \brac{1/\wp}{1}_{\sC}^{q^{d-2}} & 1/\wp^{q^{d-1}}
\CodeAfter
\line{2-3}{4-5}
\end{pNiceMatrix}.
\end{equation}
A natural problem to consider is to evaluate the determinants of $\Gamma_d(\sC)$ and $\Gamma^*_d(\sC^*)$, which are nonzero by~\eqref{E:Carlitz7.5}. To proceed, since $\omega^{(1)}=(t-\theta)\omega$, we have
\[
\omega_{\sC}^{(j)} = \omega_{\sC} \prod_{\ell=0}^{j-1} (t - \theta^{q^{\ell}}), \quad j \geq 1. 
\]
Thus by Theorem~\ref{T:DualGaussThakurSums/Zetas} and Corollary~\ref{C:GaussSumsTwistsAndZetas}, as well as~\eqref{E:Carlitzomegaevals}, we find, by evaluating both sides of of this identity at $t=\varepsilon^{q^{k}}$, $k \in \ZZ$, and dividing by $\wp'(\varepsilon^{q^{k}})$, that
\[
g^{[j]}(\chi_k) = g(\chi_k) \prod_{\ell=0}^{j-1}(\varepsilon^{q^{k}} - \theta^{q^{\ell}}).
\]
Rewriting $\Gamma_d(\sC)$ using this identity, we can now pull out the common factor of $g(\chi_k)$, $0 \leq k \leq d-1$, in each row to obtain the determinant relation
\begin{equation} \label{E:detGamma}
\det \Gamma_d(\sC) = \prod_{k=0}^{d-1} g(\chi_k) \cdot \det
\begin{pNiceMatrix}
1 & \varepsilon -\theta & \Cdots & \prod_{\ell=0}^{d-1}(\varepsilon-\theta^{q^{\ell}}) \\
1 & \varepsilon^q -\theta & \Cdots & \prod_{\ell=0}^{d-1}(\varepsilon^q-\theta^{q^{\ell}}) \\
\Vdots & \Vdots &  & \Vdots \\
1 & \varepsilon^{q^{d-1}} -\theta & \Cdots & \prod_{\ell=0}^{d-1}(\varepsilon^{q^{d-1}}-\theta^{q^{\ell}}) \\
\end{pNiceMatrix}.
\end{equation}
By performing column manipulations of this determinant we find that
\begin{equation}\label{E:alphadef}
    \rV(\varepsilon) \assign
\det \begin{pNiceMatrix}
    {\scriptstyle 1} & {\scriptstyle \varepsilon -\theta} & \Cdots & {\scriptstyle \prod_{\ell=0}^{d-1}(\varepsilon-\theta^{q^{\ell}})} \\
    {\scriptstyle 1} & {\scriptstyle \varepsilon^q -\theta} & \Cdots & {\scriptstyle \prod_{\ell=0}^{d-1}(\varepsilon^q - \theta^{q^{\ell}})} \\
     \Vdots & \Vdots &  & \Vdots \\
    {\scriptstyle 1} & {\scriptstyle \varepsilon^{q^{d-1}} -\theta} & \Cdots & {\scriptstyle \prod_{\ell=0}^{d-1}(\varepsilon^{q^{d-1}} - \theta^{q^{\ell}})} \\
\end{pNiceMatrix}
= \det \begin{pNiceMatrix}
    1 & \varepsilon & \Cdots & \varepsilon^{d-1} \\
     1 & \varepsilon^q & \Cdots & \varepsilon^{(d-1)q}\\
     \Vdots & \Vdots & & \Vdots \\
     1 & \varepsilon^{q^{d-1}}  & \Cdots & \varepsilon^{(d-1)q^{d-1}}
\end{pNiceMatrix},
\end{equation}
which is the evaluation $V(\varepsilon, \dots, \varepsilon^{q^{d-1}})$ of the Vandermonde determinant $V(t_1, \dots, t_d)$. Notably $\rV(\varepsilon) \in \FF_{q^d}^{\times}$, and $\rV(\varepsilon)^2  = (-1)^{d(d-1)/2}\, \rN(\wp'(\varepsilon)) \in \FF_q^{\times}$. Since we have $\rV(\varepsilon) = \prod_{0\leq i<j \leq d-1}(\varepsilon^{q^j}-\varepsilon^{q^i})$, we observe that $\rV(\varepsilon)^q = (-1)^{d-1} \rV(\varepsilon)$. Thus right away, combining~\eqref{E:detGamma} with Theorem~\ref{T:ThakurI-IV}(b),
\begin{equation} \label{E:detGammaCq-1}
\det \Gamma_d(\sC)^{q-1} = (-1)^d \wp \cdot \rV(\varepsilon)^{q-1} = -\wp.
\end{equation}
We extract this root $\wp$-adically and $\infty$-adically in the following theorem.
\end{subsec}

\begin{theorem} \label{T:detGamma.Carlitz}
Let $\wp \in A_+$ be irreducible of degree $d$, and fix a root $\varepsilon \in \FF_{q^d}$ of $\wp$. Let $K = k(\be(1/\wp))$ and $H=\FF_{q^d}K$. 
Define $\Gamma_d(\sC)$, $\Gamma_d^*(\sC^*) \in \GL_d(H)$ as in~\S\ref{SS:GammaCGamma*C}.
\begin{alphenumerate}
\item We have in $H$,
\[
\det \Gamma_d(\sC)^{q-1} = -\wp \quad \text{and} \quad \det \Gamma_d^*(\sC^*)^{q-1} = -\frac{1}{\wp^{q^d}}.
\]
\item \textup{($\wp$-adic Considerations)} In $H_{\fp_0}$, by choosing the root $(-\wp)^{1/(q^d-1)}$ as in Theorem~\ref{T:ThakurI-IV}(c),
\[
\det \Gamma_d(\sC) = (-\wp)^{1/(q-1)} \quad \text{and} \quad
\det \Gamma_d^*(\sC^*) = \frac{(-1)^d}{\rN(\wp'(\varepsilon)) \cdot (-\wp)^{q^d/(q-1)}}.
\]
\item \textup{($\infty$-adic Considerations)} In $\C$, by choosing $(-1)^{1/(q-1)}$ as in~\eqref{E:signpitilde} and $\wp^{1/(q-1)}$ to have sign~$1$,
\[
\det \Gamma_d(\sC) = \frac{(-1)^{d(d-1)/2} (-1)^{d/(q-1)}}{\rV(\varepsilon)} \cdot \wp^{1/(q-1)},
\]
and
\[
\det \Gamma_d^*(\sC^*) = \frac{1}{(-1)^{d/(q-1)}\,
\rV(\varepsilon)} \cdot \frac{1}{\wp^{q^d/(q-1)}}.
\]
\end{alphenumerate}
\end{theorem}

\begin{proof}
The first identity in~(a) was proved in~\eqref{E:detGammaCq-1}. Now $\det \bQ(\wp') = \rN(\wp'(\varepsilon)) \in \FF_q^{\times}$,
and so taking determinants of~\eqref{E:Carlitz7.5},
\begin{equation} \label{E:detGammaCGamma*Csimpler}
\rN(\wp'(\varepsilon)) \cdot \det \Gamma_d^*(\sC^*) \cdot \det \Gamma_d(\sC) = \frac{1}{\wp^{(q^d-1)/(q-1)}}.
\end{equation}
Taking a $(q-1)$-st power and substituting the first part of~(a) yields the second part.

For the first part of (b), noting that $\theta \equiv \varepsilon \pmod{\fp_{0}}$, we rewrite the identity in Theorem~\ref{T:ThakurI-IV}(c) for $0 \leq k \leq d-1$ as
\[
\ord_{\fp_0} \biggl( \frac{g(\chi_k)}{(-\wp)^{q^k/(q^d-1)}} - \frac{1}{(\varepsilon^{q^k} - \varepsilon)(\varepsilon^{q^{k}} - \varepsilon^q) \cdots (\varepsilon^{q^k} - \varepsilon^{q^{k-1}})}\biggr) > 0.
\]
Multiplying over all $k$, we obtain
\[
\ord_{\fp_0} \biggl( \prod_{k=0}^{d-1} \frac{g(\chi_k)}{(-\wp)^{q^k/(q^d-1)}} - \frac{1}{\rV(\varepsilon)} \biggr) > 0
\quad \stackrel{\eqref{E:detGamma}}{\Rightarrow} \quad 
\ord_{\fp_0} \biggl( \frac{\det \Gamma_d(\sC)}{(-\wp)^{1/(q-1)}} - 1 \biggr) > 0.
\]
Thus in $H_{\fp_0}$, $\det \Gamma_d(\sC) = (-\wp)^{1/(q-1)}$, with root chosen as in Theorem~\ref{T:ThakurI-IV}(c).
Substituting this into~\eqref{E:detGammaCGamma*Csimpler} yields the second part of~(b).

For the first part of~(c), it suffices from~(a) to verify that sign of $\det \Gamma_d(\sC)$ is correct. Applying Proposition~\ref{P:gg*signs} to~\eqref{E:detGamma},
\[
\sgn \bigl( \det \Gamma_d(\sC) \bigr) = \rV(\varepsilon) \prod_{k=0}^{d-1} \sgn \bigl( g( \chi_k) \bigr) = \frac{(-1)^{d/(q-1)} \rV(\varepsilon)}{\rN(\wp'(\varepsilon))},
\]
which simplifies to the desired expression as  $\rV(\varepsilon)^2  = (-1)^{d(d-1)/2}\, \rN(\wp'(\varepsilon))$. The expression for $\det \Gamma^*_d(\sC)$ then follows from~\eqref{E:detGammaCGamma*Csimpler}.
\end{proof}

\begin{subsec}{Summation identities for Dirichlet $L$-series and Gauss-Thakur sums} \label{SS:LgTrace}
In~\cite{GazdaMaurischat21}, Gazda and Maurischat analyzed the Dirichlet $L$-series of Goss type,
\[
L(\chi_{\varepsilon},s) = \sum_{a \in A_+} \frac{\chi_{\varepsilon}(a)}{a^s},
\]
especially in the case $s = q^i$, $i \geq 0$, where \cite{GazdaMaurischat21}*{Eq.~(6.1.1) and Pf.~of Thm.~6.1} imply
\begin{equation} \label{E:GM-Lvalue}
L(\chi_{\varepsilon},q^i) = -\frac{\tpi^{q^i}}{(\varepsilon - \theta^{q^i}) \cdots (\varepsilon-\theta^q)(\varepsilon - \theta) \wp'(\varepsilon) g(\chi_{\varepsilon})}, \quad i \geq 0.
\end{equation}
See also~\citelist{\cite{AnglesPellarin14}*{p.~225} \cite{Pellarin12}*{Cor.~5}}.
In light of our previous calculations, we can express $L(\chi_{\varepsilon},q^i)$ in terms of dual Gauss-Thakur sums.

Applying Theorem~\ref{T:GaussSumduality1}(a) for the Carlitz module, we have for $b \in A$ and  $j$, $k \in \ZZ$,
\begin{equation} \label{E:sumg*gtwistsCarlitz}
\sum_{k=0}^{d-1} \wp'(\varepsilon^{q^k}) g^{*[i]}(\chi_k) g^{[j]}(\chi_k;b) =
\begin{cases}
{\displaystyle \brac{ \langle b/\wp \rangle }{i-j}_{\sC}^{q^{j}}} & \textup{if $i-j \geq 0$,} \\[10pt]
{\displaystyle -\brac{ \langle b/\wp \rangle }{j-i}_{\sC^*}^{q^{j}}} & \textup{if $i-j < 0$.}
\end{cases}
\end{equation}
Now by Pellarin~\cite{Pellarin12}*{Thm.~1},
\[
\sum_{a \in A_+} \frac{a(t)}{a} = -\frac{\tpi}{(t-\theta)\omega_{\sC}(t)} =  -\tpi \Omega_{\sC}(t),
\]
and so by Corollary~\ref{C:GaussSumsTwistsAndZetas}, we have a dual Gauss-Thakur sum version of~\eqref{E:GM-Lvalue} for $j \geq 0$,
\begin{equation}
L(\chi_\varepsilon,q^i) = \sum_{a\in A_{+}} \frac{a(\varepsilon)}{a^{q^i}} =
-\tpi^{q^i} \Omega_{\sC}^{(i)}(\varepsilon) = \wp'(\varepsilon) g^{*[i]}(\chi_{\varepsilon}) \tpi^{q^i}.   
\end{equation}
Combining this with~\eqref{E:sumg*gtwistsCarlitz} we have the following.
\end{subsec}

\begin{proposition} \label{P:Lgtrace}
Let $\wp \in A_+$ be irreducible of degree $d$, and fix a root $\varepsilon \in \FF_{q^d}$. Then for $b \in A$ and $i$, $j \in \ZZ$,
\[
\sum_{k=0}^{d-1} L(\chi_k,q^{i}) g^{[j]}(\chi_k;b) =
\begin{cases}
{\displaystyle \tpi^{q^{i}}\brac{ \langle b/\wp \rangle }{i-j}_{\sC}^{q^{j}}} & \textup{if $i-j \geq 0$,} \\[10pt]
{\displaystyle -\tpi^{q^{i}}\brac{ \langle b/\wp \rangle }{j-i}_{\sC^*}^{q^{j}}} & \textup{if $i-j < 0$.}
\end{cases}
\]
\end{proposition}

As a corollary we evaluate some log-algebraic series introduced by Anderson~\cite{And96}*{\S 4.1},
\[
\lambda_m(x) \assign \sum_{a \in A_+} \frac{\be(ax)^m}{a}, \quad x \in k_{\infty},\ m \geq 0.
\]
By Anderson's log-algebraicity theorem for the Carlitz module~\cite{And96}*{Prop.~8}, for each $m \geq 0$ there is $S_m(t) \in A[t]$ such that for $b \in A$,
\begin{equation} \label{E:logalg}
\exp_{\sC} \bigl( \lambda_m(b/\wp) \bigr) = S_m \bigl( \be(b/\wp) \bigr).
\end{equation}
Anderson related $\lambda_m(b/\wp)$ to Gauss-Thakur sums~\cite{And96}*{\S 4.6} (see also~\cite{LutesP13}*{Eq.~(3.10)}, but there is missing negative sign there) by
\[
\lambda_m(b/\wp) = -\sum_{\chi \in \widehat{\FF}_q^{\times}} \biggl( \sum_{a \in \FF_{\wp}^{\times}} \chi^{-1}(a) \be(ab/\wp)^m \biggr) L(1,\chi).
\]
When $m = q^j$, the inner sum is simply $g^{[j]}(\chi;b)$, which vanishes if $\chi$ is not evaluative.

\begin{corollary} \label{C:logalg}
With notation as above, for $j \geq 0$ and $b \in A$ with $\deg b < \deg \wp$,
\[
\lambda_{q^j}(b/\wp) = \sum_{k=0}^{d-1} L(\chi_k,1) g^{[j]}(\chi_k;b) =
\begin{cases}
{\displaystyle \tpi b/\wp }& \textup{if $j=0$,} \\[10pt]
{\displaystyle -\tpi\brac{b/\wp}{j}_{\sC^*}^{q^j}} & \textup{if $j \geq 1$.}
\end{cases}
\]
\end{corollary}

\begin{remark}
By \cite{And96}*{Prop.~8(II)}, $S_1(t) = t$, and so for $j=0$, the formula for $\lambda_1(b/\wp)$ is immediately consistent with~\eqref{E:logalg}. When $j\geq 1$, \cite{And96}*{Prop.~8(VII)} implies that $S_{q^j}(t) = 0$ identically. The formula for $\lambda_{q^j}(b/\wp)$ is also consistent with~\eqref{E:logalg} by Lemma~\ref{L:qPowerBracketStarPoly}.
\end{remark}

\begin{remark}
One should compare these results with work of Angl\`es and Tavares Ribeiro \cite{AT17}*{\S 4.2}, where they consider special values of $\wp$-adic Goss-Dirichlet $L$-functions, also in terms of Gauss-Thakur sums and log-algebraicity.
\end{remark}

\section{Applications to general Drinfeld modules} \label{S:DrinfeldApps}

\begin{subsec}{The determinants of $\Gamma_d(\phi)$ and $\Gamma^*_d(\phi^*)$}
For an intermediate field $k \subseteq K \subseteq \C$, fix a rank $r$ Drinfeld module $\phi : A \to K[\tau]$ as in~\eqref{E:rhodef}. We wish to compute the determinants,
\[
\det \Gamma_d(\phi)
= \det \bQ(1/\wp')
\det\begin{pNiceMatrix}
\Upsilon(\varepsilon) & \Cdots & \Upsilon^{(r(d-1))}(\varepsilon)\\
\Vdots &  &\Vdots\\
\Upsilon(\varepsilon^{q^{d-1}}) & \Cdots & \Upsilon^{(r(d-1))}(\varepsilon^{q^{d-1}})
\end{pNiceMatrix}
\]
and
\[
\det \Gamma_d^*(\phi^*)
= \det\bQ(1/\wp') \det \begin{pNiceMatrix}
J\Psi^{(r-1)}(\varepsilon) & \Cdots & J\Psi^{(r-1)}(\varepsilon^{q^{d-1}}) \\
\Vdots & & \Vdots \\
J\Psi^{(dr-1)}(\varepsilon) & \Cdots & J\Psi^{(dr-1)}(\varepsilon^{q^{d-1}})
\end{pNiceMatrix},
\]
which we recall from \eqref{E:matrixgaussandomega}--\eqref{E:matrixgaussandzeta}.
By \eqref{E:Upsilontwist}, $\det\Upsilon^{(1)} = \det\Upsilon \det\Theta^{\tr} = (-1)^{r-1} (t-\theta)/b_r \cdot \det\Upsilon$. Therefore,
\begin{equation} \label{E:detUpsilon}
\det \Upsilon = \mu \omega_{\sC}
\end{equation}
for some $\mu \in \oK$ with $\mu^{q-1}=(-1)^{r-1}/b_r$. Indeed, both $\omega_{\sC}$ and $(\det \Upsilon)/\mu$ satisfy the same difference equation, $f^{(1)} = (t-\theta)f$, and since $\Upsilon \in \GL_r(\TT)$ by~\cite{Pellarin08}*{\S 4.2}, the two functions can differ by only an $\FF_q^{\times}$-multiple (e.g., see~\cite{AndThak90}*{Lem.~2.5.4}).
Also by Lemma~\ref{L:PsiUpsilon},
\begin{equation} \label{E:detPsi}
\det \Psi = \frac{(-1)^{r(r-1)/2}\, b_r^{(1-q^{1-r})/(1-q)}}{(\theta-t)\det \Upsilon}
= \frac{(-1)^{r(r+1)/2}\, \Omega_\sC}{\mu^{q^{1/(r-1)}}}.
\end{equation}
Our main result on $\det \Gamma_d(\phi)$ and $\det \Gamma^*_d(\phi^*)$ is Theorem~\ref{T:detGamma.General} below.
\end{subsec}

\begin{subsec}{The determinants of $\fG_{\phi}(\chi_{\varepsilon})$ and $\fG^*_{\phi^*}(\chi_{\varepsilon})$}
Recall $\fG_{\phi}(\chi_{\varepsilon})$ and $\fG_{\phi^*}^*(\chi_{\varepsilon})$ from~\eqref{E:frakG}--\eqref{E:frakG*}, where we have included $\phi$ and $\phi^*$ in the notation to indicate the dependence on~$\phi$, and let $g_{\sC}(\chi_{\varepsilon})$ and $g_{\sC^*}^*(\chi_{\varepsilon})$ denote the Gauss-Thakur sums for the Carlitz module as in~\S\ref{SS:GTsumsCarlitz}. Then Corollary~\ref{C:GaussSumsTwistsAndZetas}, together with~\eqref{E:matrixGaussZetas} and~\eqref{E:detUpsilon}--\eqref{E:detPsi}, implies for $0 \leq j \leq d-1$,
\begin{equation} \label{E:detGprep}
\det \fG_{\phi}(\chi_j)
= \frac{\det \Upsilon(\varepsilon^{q^j})}{\wp'(\varepsilon^{q^j})^r}
= \frac{\mu g_\sC(\chi_j)}{\wp'(\varepsilon^{q^j})^{r-1}},   
\end{equation}
and
\begin{equation} \label{E:detG*prep}
\det \fG_{\phi^*}^*(\chi_j)
= \frac{(-1)^{r(r-1)/2} \det \Psi^{(r-1)}(\varepsilon^{q^j})}{\wp'(\varepsilon^{q^j})^r}
= \frac{(-1)^{r-1} g^{*[r-1]}_{\sC^*}(\chi_j)}{\mu\wp'(\varepsilon^{q^j})^{r-1}}.
\end{equation}
\end{subsec}

\begin{proposition} \label{P:detfGq-1}
With notation as above,
\[
\prod_{j=0}^{d-1}\det \fG_\phi(\chi_j)^{q-1} = \frac{(-1)^{rd}\wp}{b_r^d}
\quad \text{and} \quad
\prod_{j=0}^{d-1}\det \fG^*_{\phi^*}(\chi_j)^{q-1} = \frac{(-1)^{rd}b_r^{d}}{\wp^{q^r}}.
\]
\end{proposition}

\begin{proof}
Multiplying~\eqref{E:detGprep} over the roots of~$\wp$ and raising to the $(q-1)$-st power, we obtain
\[
\prod_{j=0}^{d-1}\det \fG_\phi(\chi_j)^{q-1}
= \frac{\mu^{(q-1)d}}{\rN(\wp'(\varepsilon))^{(q-1)(r-1)}} \prod_{j=0}^{d-1} g_{\sC}(\chi_j)^{q-1}
=\frac{(-1)^{rd}\wp}{b_r^d},
\]
where the final equality follows from Theorem~\ref{T:ThakurI-IV}(b) and that $\mu^{q-1}=(-1)^{r-1}/b_r$ and $\rN(\wp'(\varepsilon)) \in \FF_q^{\times}$.
For the $\fG^*_{\phi^*}$ case, note that for $0\leq j\leq d-1$, $g_{\sC^*}^{*[r-1]}(\chi_{\varepsilon^{q^j}}) = g_{\sC^*}^*(\chi_{\varepsilon^{q^{j-r+1}}})^{q^{r-1}}$, as in the proof of Corollary~\ref{C:GaussSumsTwistsSomewhatSpecial}. Then using~\eqref{E:detG*prep},
\[
\prod_{j=0}^{d-1}\det \fG^*_{\phi^*}(\chi_j)^{q-1}
= \frac{1}{\mu^{(q-1)d}\, \rN(\wp'(\varepsilon))^{(q-1)(r-1)}} \prod_{j=0}^{d-1} g_{\sC^*}^{*}(\chi_j)^{(q-1)q^{r-1}}
=\frac{(-1)^{rd}b_r^{d}}{\wp^{q^r}},
\]
where the last equality follows from Theorem~\ref{T:dualThakurI-IV}(b).
\end{proof}

We now analyze $\det \Gamma_d(\phi)$ in a manner similar to Theorem~\ref{T:detGamma.Carlitz}. By \eqref{E:Upsilontwist}, we find
\begin{equation} \label{E:UpsilonTwist/Upsilon}
 \Upsilon^{(k)}=\Upsilon\Theta^{\tr} \bigl( \Theta^{\tr} \bigr)^{(1)} \cdots \bigl( \Theta^{\tr} \bigr)^{(k-1)}, \quad k \geq 0.  
\end{equation}

\begin{lemma} \label{L:Theta-t-linear}
As matrices in $\Mat_r(K[t])$, we have
\[
\Theta^{\tr} \bigl( \Theta^{\tr} \bigr)^{(1)} \cdots \bigl( \Theta^{\tr} \bigr)^{(r-1)} 
= \bB_1 \cdot t + \bC_1,
\]
where $\bB_1$, $\bC_1 \in \Mat_{r}(K)$ and 
$\bB_1$ is upper triangular with $[\bB_1]_{ii} = 1/b_r^{q^{i-1}}$ for $1 \leq i \leq r$. Moreover, for $k \geq 1$,
\[
\Theta^{\tr} \bigl( \Theta^{\tr} \bigr)^{(1)} \cdots \bigl( \Theta^{\tr} \bigr)^{(rk-1)}
= \bB_k \cdot t^k + \bC_k(t),
\]
where $\bB_k = \bB_1 \bB_1^{(r)} \cdots \bB_1^{(r(k-1))}$ is upper triangular with $[\bB_k]_{ii} = b_r^{-(q^{kr+i-1}-q^{i-1})/(q^r-1)}$ for $1 \leq i \leq r$ and where $\bC_k(t) \in \Mat_r(K[t])$ has entries with degree $< k$ in~$t$.
\end{lemma}

\begin{proof}
(cf.~\cite{HartlJuschka20}*{Ex.~2.5.16})
We claim that for $0 \leq \ell \leq r-1$,
\[
\Theta^{\tr} \bigl( \Theta^{\tr} \bigr)^{(1)} \cdots \bigl( \Theta^{\tr} \bigr)^{(\ell)} = \bV_\ell + \bW_\ell,
\]
where $\bV_\ell$, $\bW_\ell \in \Mat_r(K[t])$ have the form
\begin{center}
${\displaystyle
\bV_\ell = \begin{pNiceMatrix}
0 & & &  t/b_r & & \\
& \Ddots & & & \Ddots & \\
 & & & & & \smash{t/b_r^{q^\ell}} \\
 1 & & & & & \\
 & \Ddots & & & & \\
 & & 1 & & & 0
\end{pNiceMatrix}, \quad
\bW_\ell = \begin{pNiceMatrix}[first-row]
\rule[-8pt]{0pt}{8pt}
\Block{1-3}{r-\ell-1} & & & \Block{1-4}{\ell+1}\\
0 & \Cdots & 0 & * & * & \Cdots & * \\
\Vdots & & \Vdots & \Vdots & \Ddots & \Ddots & \Vdots \\
 & &  &  & &  & * \\
 & &  &  & &  & * \\
  & &  &  & &  & \Vdots \\
0 & \Cdots & 0 & * & \Cdots & & * \\
\CodeAfter%
\OverBrace[shorten,yshift=2pt]{1-1}{1-3}{}
\OverBrace[shorten,yshift=2pt]{1-4}{1-7}{}
\end{pNiceMatrix},
}$
\end{center}
and where the entries of $\bW_\ell$ have degree at most~$1$ in $t$ and if $j-i \leq r-1-\ell$ then $[\bW_\ell]_{ij} \in K$. That is, the entries in the upper triangle in the diagram above for $\bW_\ell$ can be linear in~$t$, and the ones in the lower trapezoid are all in~$K$.
Note that in $\bV_\ell$ the non-zero super-diagonal starts in the $r-\ell$ entry of the first row, and the non-zero sub-diagonal starts in the $\ell+2$ entry of the first column. Thus $\bV_{r-1}$ is a diagonal matrix with diagonal entries $t/b_r$, $t/b_r^q, \dots, t/b_r^{q^{r-1}}$. Now since
\[
\Theta^{\tr} = \begin{pNiceMatrix}
0 & \Cdots & 0 & (t-\theta)/b_r \\
1 & \Cdots & 0 & -b_1/b_r \\
  \Vdots & \Ddots & \Vdots & \Vdots \\
 0 & \Cdots & 1 & -b_{r-1}/b_r 
\end{pNiceMatrix}
= \begin{pNiceMatrix}
& & & t/b_r \\
1 & & & \\
& \Ddots & & \\
& & 1 &
\end{pNiceMatrix}
+ \begin{pNiceMatrix}
0 & \Cdots & 0 & -\theta/b_r \\
\Vdots & & \Vdots & \Vdots \\
& & & \\
0 & \Cdots & 0 & -b_{r-1}/b_r
\end{pNiceMatrix},
\]
we see that the claim is true for $\ell = 0$. Suppose that it is true for $\ell < r-1$. Then
\[
\Theta^{\tr} \bigl( \Theta^{\tr} \bigr)^{(1)} \cdots \bigl( \Theta^{\tr} \bigr)^{(\ell+1)} = (\bV_\ell + \bW_\ell)( \bV_0^{(\ell+1)} + \bW_0^{(\ell+1)}).
\]
One checks that $\bV_\ell \bV_0^{(\ell+1)} = \bV_{\ell+1}$. Letting $\bW_{\ell+1} \assign \bW_\ell \bV_0^{(\ell+1)} + \bV_\ell\bW_0^{(\ell+1)} + \bW_\ell \bW_0^{(\ell+1)}$, it is straightforward to verify using the induction hypothesis that each of $\bW_\ell \bV_0^{(\ell+1)}$, $\bV_\ell \bW_0^{(\ell+1)}$, and $\bW_\ell \bW_0^{(\ell+1)}$ satisfy the defining conditions for $\bW_{\ell+1}$, and the claim is proved. Therefore, $\Theta^{\tr} (\Theta^{\tr})^{(1)} \cdots (\Theta^{\tr})^{(r-1)} = \bV_{r-1} + \bW_{r-1} \rassign \bB_1 \cdot t + \bC_1$, where $\bB_1$, $\bC_1 \in \Mat_r(K)$ satisfy the desired conditions in the lemma. The second part of the lemma then follows from the first by taking $r$-twists of $\Theta^{\tr} (\Theta^{\tr})^{(1)} \cdots (\Theta^{\tr})^{(r-1)}$ and multiplying.
\end{proof}

Using the matrices defined in Lemma~\ref{L:Theta-t-linear}, for variables $t_1, \dots, t_d$, let
\[
\bD(t_1, \dots, t_d) \assign
\begin{pNiceMatrix}
I_r  & \bB_1 t_1+\bC_1 & \Cdots & \bB_{d-1} t_1^{d-1}+\bC_{d-1}(t_1)\\
I_r  & \bB_1 t_2+\bC_1 &  \Cdots & \bB_{d-1} t_2^{d-1}+\bC_{d-1}(t_2)\\
\Vdots  &\Vdots  &  & \Vdots \\
I_r & \bB_1 t_d+\bC_1 & \Cdots & \bB_{d-1} t_d^{d-1}+\bC_{d-1}(t_d)
\end{pNiceMatrix}
\in \Mat_{rd}(K[t_1, \dots, t_d]),
\]
and let
\[
\Delta = \Delta(t_1, \dots ,t_d)\assign \det \bD(t_1, \dots ,t_d).
\]
Similar to the simplifications in~\eqref{E:alphadef}, as each $\bB_k$, $1 \leq k \leq d-1$, is upper triangular with non-zero diagonal entries, we can perform column manipulations of this determinant to eliminate the $\bC_k(t_i)$ terms as their entries have degrees $< k$ in $t_i$, and we find
\[
\Delta = \det \left( \begin{smallmatrix}
I_r  & \bB_1 t_1 & \cdots & \bB_{d-1} t_1^{d-1}\\
I_r  & \bB_1 t_2 & \cdots & \bB_{d-1} t_2^{d-1}\\
\vdots & \vdots  &  & \vdots \\
I_r  & \bB_1 t_d & \cdots & \bB_{d-1} t_d^{d-1}
\end{smallmatrix} \right)
= \det \left( \begin{smallmatrix}
I_r & t_1 I_r  & \cdots & t_1^{d-1}I_r \\
I_r & t_2 I_r& \cdots & t_2^{d-1}I_r\\
\vdots & \vdots & & \vdots \\
I_r& t_dI_r  & \cdots & t_d^{d-1}I_r 
\end{smallmatrix} \right)
\det \left( \begin{smallmatrix}
I_r & 0 & \cdots & 0 \\
0 & \bB_1 & \cdots & 0 \\
\vdots & \vdots & \ddots & \vdots \\
0 & 0 & \cdots & \bB_{d-1}
\end{smallmatrix} \right).
\]
After performing an equal number of column and row swaps on the first determinant after the first equality, we then obtain
\[
\Delta = V(t_1, \dots, t_d)^r \det \bB_1 \cdots \det \bB_{d-1}
= b_r^{-[(q^{rd}-1)/(q^r-1) - d]/(q-1)} \cdot V(t_1, \dots, t_d)^r,
\]
where $V(t_1, \dots, t_d)$ is the Vandermonde determinant in $t_1, \dots, t_d$, and where we have used that $\det \bB_k = b_r^{-(q^{rk}-1)/(q-1)}$ by Lemma~\ref{L:Theta-t-linear}. Substituting $b_r = (-1)^{r-1}/\mu^{q-1}$ (and keeping track of signs), we obtain
\begin{equation} \label{E:Delta}
\Delta = \mu^{((q^{rd}-1)/(q^r-1) - d)} \cdot V(t_1, \dots, t_d)^r.
\end{equation}
For $0 \leq j \leq d-1$, once we divide~\eqref{E:UpsilonTwist/Upsilon} by $\wp'(\varepsilon^{q^j})$, we obtain from Corollary~\ref{C:GaussSumsTwistsAndZetas},
\[
\fG^{[k]}_\phi(\chi_j) = \fG_{\phi}(\chi_j) \cdot \Theta^{\tr} \bigl( \Theta^{\tr} \bigr)^{(1)} \cdots \bigl( \Theta^{\tr} \bigr)^{(k-1)} \Big|_{t=\varepsilon^{q^j}}, \quad k \geq 0.
\]
Taking $\Gamma_d(\phi)$ written in block form as in \eqref{E:Gammadef}, this then implies
\begin{equation} \label{E:MatrixPrep}
\Gamma_d(\phi) = \begin{pNiceMatrix}
\fG_{\phi}(\chi_0) & \Cdots & 0 \\
\Vdots & \Ddots & \Vdots \\
0 & \Cdots & \fG_{\phi}(\chi_{d-1})
\end{pNiceMatrix}
\bD(\varepsilon, \dots, \varepsilon^{q^{d-1}}).
\end{equation}

\begin{theorem} \label{T:detGamma.General}
Let $\wp \in A_+$ be irreducible of degree $d$, and fix a root $\varepsilon \in \FF_{q^d}$ of $\wp$. Define $\Gamma_d(\phi)$, $\Gamma_d^*(\phi^*) \in \GL_d(\oK)$ as in~\S\ref{SS:GTMatrices}.
\begin{alphenumerate}
\item We have
\[
\det \Gamma_d(\phi)^{q-1} = \frac{(-1)^{r}\wp}{b_r^{(q^{rd}-1)/(q^r-1)}}
\quad \text{and} \quad
\det \Gamma_d^*(\phi^*)^{q-1} = \frac{(-1)^{r}b_r^{(q^{rd}-1)/(q^r-1)}}{\wp^{q^{rd}}}.
\]
\item By choosing $(-1)^{1/(q-1)}$ as in~\eqref{E:signpitilde} and $\wp^{1/(q-1)}$ to have sign~$1$,
\[
\det \Gamma_d(\phi) = \frac{(-1)^{rd(d-1)/2} \cdot (-1)^{d/(q-1)} \cdot \mu^{(q^{rd}-1)/(q^r-1)}}{\rV(\varepsilon)^r} \cdot \wp^{1/(q-1)}
\]
and
\[
\det \Gamma_d^*(\phi^*) = \frac{1}{(-1)^{d/(q-1)} \cdot \mu^{(q^{rd}-1)/(q^r-1)} \cdot \rV(\varepsilon)^r}\cdot \frac{1}{\wp^{q^{rd}/(q-1)}}.
\]
\end{alphenumerate}
\end{theorem}

\begin{proof}
Let $\rV(\varepsilon)$ be defined as in~\eqref{E:alphadef}. Then combining \eqref{E:Delta}--\eqref{E:MatrixPrep}, we have
\begin{equation} \label{E:detGammadphi}
\det \Gamma_d(\phi) = \mu^{((q^{rd}-1)/(q^r-1) - d)}\, \rV(\varepsilon)^r \prod_{j=0}^{d-1} \det \fG_{\phi}(\chi_j).
\end{equation}
We also recall from Theorem~\ref{T:GaussSumduality1}(c) that
\begin{equation} \label{E:detGammaphiduality}
\det \Gamma_d(\phi) \cdot \det \Gamma_d^*(\phi^*) \cdot \rN(\wp'(\varepsilon))^r = \det \bP_d(1/\wp) = \wp^{-(q^{rd}-1)/(q-1)}.
\end{equation}
Though (a) follows from (b), we can prove it simply by raising~\eqref{E:detGammadphi} to the $(q-1)$-st power and applying Proposition~\ref{P:detfGq-1}. Recalling that $\mu = (-1)^{r-1}/b_r$ and $\rV(\varepsilon)^{q-1} = (-1)^{d-1}$, a straightforward calculation yields the identity for $\det \Gamma_d(\phi)^{q-1}$ in~(a). The identity for $\det \Gamma_d^*(\phi^*)^{q-1}$ in~(a) follows from \eqref{E:detGammaphiduality}.
For (b), we combine~\eqref{E:detGprep} with~\eqref{E:detGammadphi} and find
\begin{align*}
\det \Gamma_d(\phi) &= \frac{\mu^{(q^{rd}-1)/(q^r-1)} \, \rV(\varepsilon)^r}{\rN(\wp'(\varepsilon))^{r-1}} \prod_{j=0}^{d-1} g_{\sC}(\chi_j), \\
\intertext{and then after applying~\eqref{E:detGamma} and Theorem~\ref{T:detGamma.Carlitz}(c),}
&= (-1)^{d(d-1)/2} \cdot (-1)^{d/(q-1)} \cdot \mu^{(q^{rd}-1)/(q^r-1)} \cdot \frac{\rV(\varepsilon)^{r-2}}{\rN(\wp'(\varepsilon))^{r-1}} \cdot \wp^{1/(q-1)}.
\end{align*}
Using that $\rN(\wp'(\varepsilon)) = (-1)^{d(d-1)/2} \rV(\varepsilon)^2$ completes the proof for~$\det \Gamma_d(\phi)$, and then the identity for $\det \Gamma_d^*(\phi^*)$ follows from~\eqref{E:detGammaphiduality}.
\end{proof}

\begin{subsec}{Special values of Drinfeld modular forms in Tate algebras} \label{SS:DMFspecialvalues}
We explore connections between Drinfeld modular forms of arbitrary rank with values in~$\TT$ and dual Gauss-Thakur sums. We hope the methods of the present section shed additional light on certain modular form identities. For precise definitions of Drinfeld modular forms of higher rank, see~\cites{BassonBreuer17, BassonBreuerPink24, ChenGezmis22, Gekeler17, GezmisUlkem27}, and for modular forms in Tate algebras~\cites{Pellarin12, Pellarin25, PellarinPerkins18, Perkins14}.

Throughout we primarily adopt the conventions and notation of Chen and Gezmi\c{s} in~\cite{ChenGezmis22} for Drinfeld modular forms. For $r \geq 2$, let
$\HH^r = \PP^{r-1}(\C) \setminus \{k_{\infty}$-rational hyperplanes$\}$ be the Drinfeld upper half space. For $\bz = (z_1, \dots, z_r)^{\tr} \in \HH^r$, $z_1, \dots, z_r$ are $k_{\infty}$-linearly independent, which we normalize so that $z_r=1$. For $\bz \in \HH^r$, define $\Lambda_{\bz} = A z_1 + \dots + A z_r \subseteq \C$, and let $\phi^{\bz}$ be the Drinfeld module associated to $\Lambda_{\bz}$,
\[
\phi^{\bz}_\theta \assign \theta + b_1(\bz)\tau + \cdots + b_r(\bz)\tau^r.
\]
As functions on $\HH^r$, the coefficients $b_i(\bz)$ are then Drinfeld modular forms of rank~$r$, weight $q^i-1$, and type~$0$ for $\GL_r(A)$ \citelist{\cite{BassonBreuer17}*{Prop.~3.4} \cite{ChenGezmis22}*{Ex.~2.8} \cite{Gekeler17}*{Ex.~1.12}}. The leading coefficient $b_r(\bz)$ is the discriminant form $\Delta_r(\bz)$ of Gekeler as in~\cite{Gekeler17}*{\S 3}. Gekeler~\citelist{\cite{Gekeler89}*{\S 4} \cite{Gekeler17}*{Thm.~3.8}} further defined a rank~$r$ modular form $h_r(\bz)$, which has weight $(q^r-1)/(q-1)$ and type~$1$, such that (see also~\citelist{\cite{ChenGezmis22}*{Ex.~2.8(iii)} \cite{BassonBreuerPink24}*{Thm.~5.3(2)}})
\begin{equation} \label{E:brz}
\Delta_r(\bz) = b_r(\bz) = (-1)^{r-1} \tpi^{q^r-1} h_r(\bz)^{q-1}. 
\end{equation}
We have Eisenstein series of weight~$k$ and type~$0$ (see \citelist{\cite{BassonBreuer17}*{Ex.~3.1} \cite{ChenGezmis22}*{Ex.~2.8(i)} \cite{Gekeler17}*{\S 1}}),
\[
E_k(\bz) = \sumprime{(a_1, \dots, a_r) \in A^r} \frac{1}{(a_1 z_1 + \dots + a_r z_r)^k}.
\]
Right away we see that if we take the dual special functions $\zeta_i(\bz,t) \in \TT$, $1 \leq i \leq r$, for~$\phi^{\bz}$ as in Corollary~\ref{C:samezetas}, then
\[
E_{q^j-1}(\bz) = -\sum_{i=1}^r z_i \bigl( \zeta_i(\bz,t) \bigr)^{(j)} \big|_{t=\theta}, \quad j \geq 0.
\]

Chen and Gezmi\c{s}~\cite{ChenGezmis22} investigate Eisenstein series taking values in $\TT$, originally constructed by Pellarin~\cite{Pellarin12}*{\S 2} and Perkins~\cite{Perkins14}*{\S 4}. For $k \geq 0$, set
\[
\cE_k(\bz,t) \assign \Biggl(\, \sumprime{(a_1, \dots, a_r) \in A^r} \frac{a_1(t)}{(a_1 z_1 + \dots + a_r z_r)^k}, \ \dots, \sumprime{(a_1, \dots, a_r) \in A^r} \frac{a_r(t)}{(a_1 z_1 + \dots + a_r z_r)^k} \Biggr)^{\tr},
\]
which, for fixed $\bz$, is a vector in $\TT^r$. Let $\Upsilon(\bz,t)$, $\Psi(\bz,t) \in \GL_r(\TT)$ be defined with respect to $\phi^{\bz}$ as in \eqref{E:Upsilondef}, \eqref{E:Psidef}, and let $\bsomega(\bz,t)$ and $\bszeta(\bz,t)$ be their first column and row respectively. Then again by Corollary~\ref{C:samezetas},
\begin{equation} \label{E:cEqjzeta}
\cE_{q^j}(\bz,t) = \bigl( -\bszeta(\bz,t)^{(j)} \bigr)^{\tr}.
\end{equation}
Also, we note that $\Upsilon(\bz,t)$ is exactly the same as $\cF(\bz,t)$ in~\cite{ChenGezmis22}*{Eq.~(3.1)}, and moreover,
\begin{equation} \label{E:CGProp3.4}
\det \Upsilon(\bz,t) = \tpi^{-(q^r-1)/(q-1)}h_r(\bz)^{-1}\, \omega_{\sC}(t),
\end{equation}
by \cite{ChenGezmis22}*{Prop.~3.4}. Thus, we see that
\begin{equation} \label{E:muz}
\mu(\bz) \assign \tpi^{-(q^r-1)/(q-1)} h_r(\bz)^{-1}
\end{equation}
is the value of $\mu$ in~\eqref{E:detUpsilon} and elsewhere, as a function of~$\bz$. Let $\cC(\bz,t)$ be the cofactor matrix of $\Upsilon(\bz,t)$, so that $\cC(\bz,t)^{\tr} \Upsilon(\bz,t) = \det \Upsilon(\bz,t)\, I_r$.
\end{subsec}

\begin{theorem} \label{T:ChenGezmis11}
For fixed $\bz \in \HH^r$ the following identities hold in $\TT^r$.
\begin{alphenumerate}
\item \textup{(Pellarin~\cite{Pellarin12}*{Thm.~8}, case $r=2$):} For $z \in \C \setminus k_{\infty}$,
\[
\Biggl(\, \sumprime{(a,b) \in A^2} \frac{a(t)}{az  + b},\ \sumprime{(a,b) \in A^2} \frac{b(t)}{az  + b} \Biggr) = \frac{\tpi^{q+1}\,h_2(z,1)}{(t-\theta)\omega_{\sC}(t)} \cdot \bigl( \omega_1^{(1)}(t), -\omega_z^{(1)}(t) \bigr).
\]
\item \textup{(Chen-Gezmi\c{s}~\cite{ChenGezmis22}*{Thm.~3.2}, case $r \geq 2$):} For $\bz \in \HH^r$,
\[
\cE_1(\bz,t)^{\tr} = \frac{\tpi^{(q^r-1)/(q-1)}\,h_r(\bz)}{(t-\theta) \omega_{\sC}(t)} \cdot \bigl( \cC(\bz,t)_{11}, \dots, \cC(\bz,t)_{r1} \bigr).
\]
\end{alphenumerate}
\end{theorem}

\begin{remark}
Once we know how $\mu(\bz)$ is related to $h_r(\bz)$ from~\cite{ChenGezmis22}*{Prop.~3.4}, as in~\eqref{E:muz}, we see below these results can also be proved by combining information about rigid analytic trivializations from Lemma~\ref{L:PsiUpsilon} and Corollary~\ref{C:samezetas} (relying on results of~\cite{Ferraro25} stated in Theorems~\ref{T:FerraroUniv} and~\ref{T:Ferraropairing}), without appealing further to the theory of Drinfeld modular forms. Various extensions of this theorem can be found in \citelist{\cite{ChenGezmis26}*{Lem.~5.5} \cite{GezmisUlkem27}*{Thm.~1.3(i)}}, and it would be interesting to investigate the extent to which the methods of the present paper could prove them from this standpoint (e.g., cf.\ Lemma~\ref{L:PsiUpsilon} and~\cite{GezmisUlkem27}*{Eq.~(3.8)}).
\end{remark}

\begin{proof}[Alternative proof of Theorem~\ref{T:ChenGezmis11}]
Part (a) is simply a special case of (b). For~(b) we need~\cite{ChenGezmis22}*{Prop.~3.4} to pin down the sign of $\det \Upsilon(\bz,t)$ as a function of $\bz$ as in~\eqref{E:CGProp3.4} (cf.~\eqref{E:detUpsilon}, \eqref{E:muz}). In any case, similar to the proof of Corollary~\ref{C:samezetas},
\[
\cE_1(\bz,t)^{\tr} = -\bszeta(\bz,t) = (-1,0, \dots, 0) \Psi(\bz,t)
= \biggl( \frac{1}{t-\theta},0,\dots, 0 \biggr) \cC(\bz,t)^{\tr} \det \Upsilon(\bz,t)^{-1},
\]
using Lemma~\ref{L:PsiUpsilon} at the end (cf.~\cite{ChenGezmis22}*{Prop.~3.18}). The result follows from~\eqref{E:CGProp3.4}.
\end{proof}

As usual, fix irreducible $\wp \in A_+$ of degree~$d$ with root $\varepsilon \in \FF_{q^d}$. For a character $\chi \in \widehat{\FF}_{\wp}^{\times}$, $1 \leq i \leq r$ and $\bz \in \HH^r$, we will denote $g_i(\bz,\chi)$, $g_i^*(\bz,\chi)$, $\fG(\bz,\chi)$, $\fG^*(\bz,\chi)$, $\Gamma(\bz)$, $\Gamma^*(\bz)$, etc. for the corresponding constructions on Gauss-Thakur sums and their duals as in \S\S\ref{S:GTsums}--\ref{S:MatrixGaussSection}, with respect to the Drinfeld modules $\phi^{\bz}$ and their adjoints~$(\phi^{\bz})^*$. The connection with dual Gauss-Thakur sums arrives through evaluation, as by Corollary~\ref{C:GaussSumsTwistsAndZetas} and~\eqref{E:cEqjzeta},
\begin{equation} \label{E:VDMF-dualGTsums}
\cE_{q^j}(\bz,\varepsilon) = -\biggl( \frac{g_1^{*[j]}(\bz,\chi_{\varepsilon})}{\wp'(\varepsilon)}, \dots, \frac{g_r^{*[j]}(\bz,\chi_{\varepsilon})}{\wp'(\varepsilon)}  \biggr)^{\tr}.
\end{equation}
From \eqref{E:detG*prep}, Proposition~\ref{P:detfGq-1}, and Theorem~\ref{T:detGamma.General}, we immediately obtain identities for $h_r(\bz)$ and $\Delta_r(\bz)$ in terms of dual Gauss-Thakur sums.

\begin{corollary} \label{C:DMFs}
Fix $\wp \in A_+$ irreducible of degree $d$ and root $\varepsilon \in \FF_{q^d}$. For $\bz \in \HH^r$,
\begin{alphenumerate}
\item $\begin{aligned}[t]
&h_r(\bz) = \frac{(-1)^{r-1} \wp'(\varepsilon)^{r-1}}{\tpi^{(q^r-1)/(q-1)} g_{\sC}^{*[r-1]}(\chi_\varepsilon)} \cdot \det \fG^*(\bz,\chi_{\varepsilon}), \\
&h_r(\bz)^{(q^{rd}-1)/(q^r-1)} = \frac{(-1)^{d/(q-1)} \rV(\varepsilon)^r \wp^{q^{rd}/(q-1)}}{\tpi^{(q^{rd}-1)/(q-1)}} \cdot \det \Gamma_d^*(\bz);
\end{aligned}$
\item $\begin{aligned}[t]
&\Delta_r(\bz)^d
= (-1)^{rd} \wp^{q^r} \prod_{j=0}^{d-1} \det \fG^*(\bz,\chi_j)^{q-1}, \\
&\Delta_r(\bz)^{(q^{rd}-1)/(q^r-1)}
= (-1)^{r} \wp^{q^{rd}} \cdot \det \Gamma_d^*(\bz)^{q-1}.
\end{aligned}$
\end{alphenumerate}
\end{corollary}

\end{document}